\documentclass[leqno, myheadings, twoside]{amsart}
\usepackage{amsfonts}
\usepackage{amsmath, amsthm, amssymb, amscd, amsxtra,graphicx}
\usepackage{latexsym, amsfonts}
\usepackage{url}
\usepackage{texdraw}
\usepackage{epsfig}
\def\1-1{(1.1)}

\makeatletter \@addtoreset{equation}{section}

\makeatletter \renewcommand{\@biblabel}[1]{#1.}

\theoremstyle{remark}

\begin{document}
\title [The Weyl-type asymptotic formula for biharmonic Steklov eigenvalues]
{The Weyl-type asymptotic formula for biharmonic Steklov eigenvalues
 with Dirichlet boundary condition on Riemannian manifolds}
\author{Genqian Liu}


\maketitle Department of Mathematics, Beijing Institute of
Technology,
 Beijing, the People's Republic of China.
 \ \
E-mail address:  liugqz@bit.edu.cn

\vskip 0.46 true cm

\vskip 15 true cm

\begin{abstract}   Let $\Omega$ be a bounded domain with $C^2$-smooth
  boundary in an $n$-dimensional
  oriented Riemannian manifold. It is well-known that for the bi-harmonic equation
   $\Delta^2 u=0$ in $\Omega$ with the $0$-Dirichlet boundary
 condition, there exists an infinite set $\{u_k\}$ of biharmonic functions in
 $\Omega$ with positive eigenvalues $\{\lambda_k\}$ satisfying
 $\Delta u_k+ \lambda_k \varrho \frac{\partial u_k}{\partial \nu}=0$ on
  the boundary $\partial \Omega$.
 In this paper, by a new method we establish the Weyl-type asymptotic formula
    for the counting function of the biharmonic Steklov
  eigenvalues $\lambda_k$.
   \end{abstract}

\vskip 1.39 true cm

\section{ Introduction}

\vskip 0.45 true cm

  Spectral asymptotics for partial differential operators have been
the subject of extensive research for over a century. It has
attracted the attention of many outstanding mathematicians and
physicists. Beyond the beautiful asymptotic formulas that are
intimately related to the geometric properties of the domain and its
boundary, a sustaining force has been its important role in
 mathematics, mechanics and theoretical physics (see, for
example, \cite{Ch1}, \cite{Ch2}, \cite{CH}, \cite{CLN}, \cite{ES},
\cite{Ho}, \cite{Iv}, \cite{Ka}, \cite{MS}, \cite{Pa1}, \cite{Pa2},
\cite{Po}, \cite{Sa}, \cite{Sar}, \cite{St}, \cite{SV}, \cite{We4}).

Let $(\mathcal{M}, g)$ be an oriented Riemannian manifold of
dimension $n$
 with a positive definite metric tensor $g$, and
let $\Omega\subset {\mathcal{M}}$ be a bounded domain with
$C^2$-smooth boundary $\partial \Omega$. Assume $\varrho$ is a
non-negative bounded function defined on $\partial \Omega$.
   We consider the following classical biharmonic Steklov eigenvalue problem:
        \begin{eqnarray} \label{1-1}  \left\{\begin{array}{ll}\triangle^2_g u=0
      \quad \;\; &\mbox{in}\;\; \Omega,\\
  u=0 \;\; & \mbox{on}\;\; \partial \Omega,\\
   \triangle_g u+\lambda \varrho \frac{\partial u}{\partial \nu}=0
   \;\;& \mbox{on}\;\; \partial \Omega, \end{array} \right. \end{eqnarray}
   where $\nu$ denotes the inward unit
   normal vector to $\partial \Omega$,
   and $\triangle_g$ is the Laplace-Beltrami operator
  defined in local coordinates by the expression,
   \begin{eqnarray*} \triangle_g =\frac{1}{\sqrt{|g|}}\sum_{i,j=1}^n
   \frac{\partial}{\partial x_i} \left( \sqrt{|g|}\,g^{ij}
   \frac{\partial}{\partial x_j}\right).\end{eqnarray*}
Here $| g | : =  det(g_{ij})$ is the determinant of the metric
tensor, and $g^{ij}$ are the components of the inverse of the metric
tensor $g$.

\vskip 0.12 true cm

 The problem (\ref{1-1}) has nontrivial solutions $u$ only for a discrete set
 of $\lambda = \lambda_k$, which are called biharmonic Steklov
  eigenvalues (see \cite{FGW}, \cite{KS}, \cite{Pa1} or \cite{WX}).
  Let us enumerate the eigenvalues in increasing order:
\begin{eqnarray*}  & 0<\lambda_1 \le \lambda_2 \le  \cdots \le
   \lambda_k\le\cdots,
   \end{eqnarray*}
  where each eigenvalue is counted as many times as its multiplicity.
     The corresponding eigenfunctions $\frac{\partial u_1}{\partial \nu},
     \frac{\partial  u_2}{\partial \nu},
   \cdots, \frac{\partial u_k}{\partial \nu}, \cdots$
   form a complete orthonormal basis in $L_\varrho^2(\partial \Omega)$
   (see, Proposition 3.5).
It is clear that $\lambda_k$ can be characterized variationally as
\begin{eqnarray*}
 \lambda_1=\frac{\int_\Omega |\triangle_g u_1|^2 dR} {
\int_{\partial \Omega} \varrho \big(\frac{\partial u_1}{\partial
\nu}\big)^2 ds} =
  \inf_{\underset{
  0\ne \frac{\partial v}{\partial \nu}\in L^2(\partial \Omega)}{v\in H_0^1 (\Omega)\cap H^2 (\Omega)
  }}\;\; \frac{\int_\Omega |\triangle_g v|^2 dR} {
\int_{\partial \Omega} \varrho \big(\frac{\partial v}{\partial \nu}\big)^2 ds},\;\;\qquad \qquad \\
 \lambda_k=\frac{\int_\Omega |\triangle_g u_k|^2 dR} {
\int_{\partial \Omega} \varrho \big(\frac{\partial u_k}{\partial
\nu}\big)^2 ds} =
  \max_{\underset{codim(\mathcal {F})=k-1}
  {\mathcal{F}\subset H_0^1 (\Omega)\cap H^2 (\Omega)}
}  \,\,\inf_{\underset{0\ne \frac{\partial v}{\partial\nu}\in
L^2(\partial \Omega)}{v\in \mathcal{F}}}\;\; \frac{\int_\Omega
|\triangle_g v|^2 dR} { \int_{\partial \Omega} \varrho
\big(\frac{\partial v}{\partial \nu}\big)^2 ds}, \quad \; k=2,3, 4,
\cdots
\end{eqnarray*} where
 $H^m(\Omega)$ is the Sobolev space, and where $dR$ and $ds$ are the Riemannian
elements of volume and area on $\Omega$ and $\partial \Omega$,
respectively.

In elastic mechanics, when the weight of the body $\Omega$
is the only body force, the stress function $u$
 must satisfy the equation $\Delta^2 u=0$ in $\Omega$ (see, p.$\,$32 of \cite{TG}).
 In addition, the boundary
 condition in (\ref{1-1}) has an interesting
interpretation in theory of elasticity.
    Consider the model problem
(see \cite{FGW}):
 \begin{eqnarray} \label{1-2}  \left\{\begin{array}{ll}\triangle^2 u=f
      \quad \;\; &\mbox{in}\;\; \Omega,\\
  u=0,\;\; \triangle u+(1-\sigma) \iota \,\frac{\partial u}{\partial \nu}=0 &
  \mbox{on}\;\; \partial
  \Omega, \end{array} \right. \end{eqnarray}
where $\Omega\subset {\Bbb R}^2$ is an open bounded domain with
smooth boundary, $\sigma\in (-1, 1/2)$ is the Poisson ratio and
$\iota$ is the mean  curvature of the boundary $\partial \Omega$.
Problem (\ref{1-2}) describes the deformation $u$ of the linear
elastic supported plate $\Omega$ under the action of the transversal
exterior force $f = f(x)$, $x\in \Omega$. The Poisson ratio $\sigma$
 of an elastic material is the negative transverse strain divided by
the axial strain in the direction of the stretching force. In other
words, this parameter measures the transverse expansion
(respectively, contraction) if $\sigma>0$ (respectively, $\sigma<0$)
when the material is compressed by an external force. We refer to
   \cite{La}, \cite{Vi} for more details. The restriction on the
   Poisson ratio is
due to thermodynamic considerations of strain energy in the theory
of elasticity. As shown in \cite{La}, there exist materials for
which the Poisson ratio is negative and the limit case $\sigma=-1$
  corresponds to materials with an infinite flexural rigidity (see,
 p.$\,$456 of \cite{St}). This limit value for $\sigma$ is strictly related to the
eigenvalue problem (\ref{1-1}). Hence, the limit value $\sigma=-1$,
which is not allowed from a physical point of view, also changes the
structure of the stationary problem (\ref{1-2}): For example (see
\cite{FGW}), when $\Omega$ is the unit disk and $\lambda_1=
(1-\sigma)\iota= 1-\sigma=2$,  $\,$(\ref{1-2}) either admits an
infinite number of solutions or it admits no solutions at all,
depending on $f$.

 Problem (\ref{1-1}) is also important in conductivity
 and biharmonic analysis because the related problem was initially studied by Calder\'{o}n
 (cf.$\,$\cite{Cal}). This connection arises because the set of the eigenvalues
 for the biharmonic Steklov problem is the same as the set
 of eigenvalues of the well-known ``Neumann-to-Laplacian'' map for biharmonic equation
 (This map associates each normal derivative $\partial u/\partial \nu$
 defined on the boundary
 $\partial \Omega$ to the restriction $(\triangle u)\big|_{\partial
 \Omega}$ of the Laplacian of $u$
 for the biharmonic function
 $u$ on $\Omega$, where the biharmonic function $u$ is uniquely
 determined by $u\big|_{\partial \Omega}=0$  and
  $(\partial u/\partial \nu)\big|_{\partial \Omega}$).

In the general case the eigenvalues $\lambda_k$ can not be evaluated
explicitly. In particular, for large $k$ it is difficult to
calculate them numerically. In view of the important applications,
one is interested in finding the asymptotic formulas for $\lambda_k$
as $k\to \infty$. However, for a number of reasons it is traditional
in such problems to deal with the matter the other way round, i.e.,
to study the sequential number $k$ as a function of $\tau$. Namely,
let us introduce the counting function $A(\tau)$ defined as the
number of eigenvalues $\lambda_k$ less than or equal to a given
$\tau$. Then our asymptotic problem is reformulated as the study of
the asymptotic behavior of
 $A(\tau)$ as $\tau\to +\infty$.

In order to better understand our problem (\ref{1-1}) and its
asymptotic behavior, let us mention the Steklov eigenvalue problem
for the harmonic equation
\begin{eqnarray} \label{1-3}  \left\{\begin{array}{ll} \triangle_g v=0&
\quad \, \mbox{in}\;\; \Omega,\\
\frac{\partial v}{\partial \nu}+ \eta \varrho v =0 &\quad \,
\mbox{on}\;\; \partial \Omega,
\end{array} \right.\end{eqnarray}
 where $\eta$ is a real number. This problem was first introduced by
 V. A. Steklov for bounded domains in the plane in \cite{St} (The reader should be aware
 that ``Steklov'' is also often transliterated as ``Stekloff''.) His motivation
 came from physics. The function $v$ represents
 the steady state temperature on $\Omega$ such that
 the flux on the boundary is proportional
 to the temperature (In two dimensions, it can also be interpreted as
 a membrane with whole mass concentrated on the boundary).
 For the harmonic Steklov eigenvalue problem (\ref{1-3}), in a special case
 in two dimensions,
 {\AA}. Pleijel \cite{Pl} outlined an investigation of the
 asymptotic behavior of both eigenvalues and the eigenfunctions.
 In 1955, L. Sandgren \cite{Sa} established the asymptotic formula of
 the counting function $B(\tau)=
 \#\{\eta_k\big|\eta_k \le \tau\}$:
   \begin{eqnarray} \label{1-4}   B(\tau)\sim
   \frac{\omega_{n-1}\tau^{n-1}}{(2\pi)^{n-1}}
 \int_{\partial \Omega} \varrho^{n-1} ds \quad \;\mbox{as}\;\;
 \tau\to +\infty,\end{eqnarray}
 i.e.,
 \begin{eqnarray*}  \lim_{\tau \to +\infty} \; \frac{B(\tau)}{\tau^{n-1}}=
 \frac{\omega_{n-1}}{(2\pi)^{n-1}}
 \int_{\partial \Omega} \varrho^{n-1} ds,\end{eqnarray*}
 where $\omega_{n-1}$ is the volume of the  unit
 ball of ${\Bbb R}^{n-1}$,  and the integral is
 over the boundary $\partial \Omega$.
This asymptotic behavior is motivated by the similar one for the
eigenvalues of the Dirichlet Laplacian. The classical result for the
Dirichlet (or Neumann) eigenvalues of the Laplacian on a smooth
 bounded domain is Weyl's formula (see \cite{Ch1}, \cite{CH} or
\cite{We2}):
\begin{eqnarray}  \label{1-5}   N(\tau,\Omega)\sim
 \frac{\omega_n}{(2\pi)^{n}} \big(\mbox{vol}(\Omega)\big) \tau^{n/2}
\quad \mbox{as}\;\; \tau\to +\infty,\end{eqnarray}
where $N(\tau, \Omega)$ is the number of the Dirichlet (or Neumann) eigenvalues of domain $\Omega$
 less than or equal to a given $\tau$.
 In the case of two-dimensional Euclidean space, Pleijel \cite{Ple} in 1950 proved an
  asymptotic formula for the
  eigenvalues $\Xi_k^2$ of a clamped plate problem:
\begin{eqnarray}   \left\{ \begin{array}{ll} \triangle^2 u
  -\Xi^2 u=0 \quad \;\; &\mbox{in}\;\; \Omega,\\
  u=\frac{\partial u}{\partial \nu}=0\quad \;\;  &  \mbox{on}\;\;
  \partial \Omega. \end{array} \right. \end{eqnarray}
   Grub\cite{Gr} and Ashbaugh, Gesztesy, Mitrea and
 Teschl \cite{AGMT} obtained Weyl's asymptotic formula for the eigenvalues $\Lambda_k$ of the
 buckling problem in
${\Bbb R}^n$:
\begin{eqnarray} \left\{ \begin{array}{ll} \triangle^2 u+\Lambda
  \triangle u=0 \quad \;\; &\mbox{in}\;\; \Omega,\\
  u=\frac{\partial u}{\partial \nu}=0 \quad \;\;  &  \mbox{on}\;\;
   \partial \Omega. \end{array} \right. \end{eqnarray}
 Note that for the Dirichlet eigenvalues, the Neumann
 eigenvalues,  the buckling eigenvalues and
 the square root of the clamped plate eigenvalues
  in a fixed domain, their counting functions
  have the same asymptotic
 formula (\ref{1-5}) (see, for example,
  \cite{Ch1}, \cite{Liu}, \cite{Gr}, \cite{We1} and \cite{We2}).

The study of asymptotic behavior of the biharmonic Steklov
eigenvalues is much more difficult than that of the harmonic
 Steklov eigenvalues. It had been a challenging problem in the
 past 50 years. The main stumbling
 block that lies in the way is the estimates for the distribution
  of the boundary eigenvalues for bi-harmonic
     equations with suitable boundary conditions. Some important works have
 contributed to the research of this problem, for example, L. E. Payne \cite{Pa1},
J. R. Kuttler and V. G. Sigillito \cite {KS}, A. Ferrero, F. Gazzola
and T. Weth \cite{FGW},  Q. Wang and C. Xia \cite{WX}, and others.

 \vskip 0.16 true cm

In this paper, by a new method we establish the Weyl-type asymptotic
formula for the counting function of the biharmonic Steklov
eigenvalues. The main result is the following:

\vskip 0.25 true cm

 \noindent  {\bf Theorem 1.1.} \ \  {\it
Let $(\mathcal{M},g)$ be an $n$-dimensional oriented Riemannian
manifold, and let $\Omega\subset \mathcal{M}$ be a bounded domain
with $C^{2}$-smooth boundary $\partial \Omega$. Then
\begin{eqnarray} \label{1-6} A(\tau)\sim
\frac{\omega_{n-1}\tau^{n-1}}{(4\pi)^{(n-1)}}
 \int_{\partial \Omega} \varrho^{n-1} ds \quad \, \; \mbox{as}\;\; \tau\to +\infty,\end{eqnarray}
where $A(\tau)$ is defined as before}.

\vskip 0.25 true cm

 \noindent  {\bf Corollary 1.2.} \ \  {\it
Let $(\mathcal{M},g)$ be an $n$-dimensional oriented Riemannian
manifold, and
 let $\Omega\subset \mathcal{M}$ be a bounded domain with
  $C^{2}$-smooth boundary $\partial \Omega$.  If, in problem
(\ref{1-1}), $\varrho\equiv 1$ on $\partial \Omega$,
 then \begin{eqnarray} \label{1-7}  \lambda_k \sim
 (4\pi) \left(\frac{k}{\omega_{n-1} (\mbox{vol}(\partial \Omega))}\right)^{1/(n-1)}
 \quad \, \mbox{as}\;\; k\to +\infty.\end{eqnarray}}

\vskip 0.16 true cm

We outline the idea of the proof of Theorem 1.1.
 First, we make a division of $\bar \Omega$ into subdomains (by dividing $\partial \Omega$
 into sufficiently small parts,
 then taking a depth $\sigma>0$ (small enough) in the direction
 of inner normal of
 $\partial \Omega$ to form a finite number of $n$-dimensional subdomains).
   Under a sufficiently fine division of $\partial \Omega$ (also $\sigma$
sufficiently small),  $g^{ik}$ and $\varrho$ can be replaced by
constants because their variant will be small, so that the
corresponding subdomains whose partial boundaries are situated at
the $\partial \Omega$ can be approximated by Euclidean cylinders.
 Next, we construct three Hilbert spaces of
 functions and their self-adjoint linear transformations whose
 eigenfunctions are just the Steklov eigenfunctions with
  corresponding  boundary conditions. It can be
 shown that these Steklov eigenvalue problems have the same boundary
 conditions on the base of each cylinder as the original one in problem
 (\ref{1-1}) but they have relevant boundary conditions on the
 other parts of a cylinder.
 In particular, on each cylindrical surface, these boundary conditions will be
  one of the three forms $u=\Delta_g u=0$,
 $\Delta_g u=\frac{\partial (\Delta_g u)}{\partial \nu}=0$ and
 $\frac{\partial u}{\partial \nu}=\frac{\partial (\Delta_g u)}{\partial \nu}=0$.
  The main purpose of constructing such Steklov problems is that when putting
   together such cylinders, we can obtain global
  upper and lower estimates for the counting function
  $A(\tau)$ of the original Steklov problem
   (i.e., $A^0(\tau)\le A(\tau)\le A^d(\tau) \le A^f(\tau)$ for all
   $\tau>0$, see Sections 3, 6). For each Euclidean cylinder, by using a cubical net
    we can divide the base of
    the cylinder into $(n-1)$-dimensional cubes and some smaller
   parts which intersect boundary of the base, so that we get $n$-dimensional
   parallelepipeds and
   some smaller $n$-dimensional cylinders.
 As for the $n$-dimensional parallelepiped, we can explicitly calculate the
 Steklov eigenfunctions and eigenvalues by
 separating variables, and then we can compute the asymptotic distribution of eigenvalues by means of
 the well-known variational
 methods used by H. Weyl \cite{We4}, R. Courant and D. Hilbert \cite{CH}
   in the case of the membranes. Meanwhile, for each small $n$-dimensional cylinder,
   by introducing a
 nice transformation we may map it into a special cylinder whose counting functions
 of Steklov eigenvalues can also be estimated.
 Finally, applying normal coordinates system at a fixed point of each subdomain of a division
   and combining
 these estimates, we establish the desired asymptotic
 formula for $A(\tau)$.
 Note that
 the Holmgren uniqueness theorem for the solutions of elliptic equations
  plays a crucial role in this paper.

 This paper is organized as follows. In Section 2, we prove two
 compact trace lemmas for bounded domains with piecewise smooth
 boundaries. In Section 3, we define various self-adjoint transformations
  on the associated Hilbert spaces of functions, and give the connections between
  the eigenfunctions of self-adjoint transformations and the
  Steklov eigenfunctions (corresponding to different kinds of boundary
  conditions).
 Section 4 is dedicated to deriving the explicit
 formulas for the
 biharmonic Steklov eigenvalues and eigenfunctions in an $n$-dimensional rectangular parallelepiped
 of ${\Bbb R}^n$, which depends on a key calculation for the
 solutions of biharmonic equations.
  The counting functions of Steklov eigenvalues for general cylinder of the Euclidean space
  are dealt with in Section 5. In the final section,
  we prove Theorem 1.1 and Corollary 1.2 on Riemannian manifolds.

\vskip 1.39 true cm

\section{Compact trace Lemmas}

\vskip 0.45 true cm

An $n$-dimensional cube in ${\Bbb R}^n$ is the set $\{x\in {\Bbb
R}^n\big| 0\le x_i \le a, \, i=1, \cdots, n\}$. \vskip 0.16 true cm

Let $f$ be a real-valued function defined in an open set $\Omega$ in
${\Bbb R}^{n}$ ($n\ge 1$). For $y\in \Omega$ we call $f$ {\it real
analytic at $y$} if there exist $a_\beta \in {\Bbb R}^1$ and a
neighborhood $U$ of $y$ (all depending on $y$) such that $$ f(x)=
\sum_\beta a_\beta (x-y)^\beta$$ for all $x$ in $U$. We say $f$ is
{\it real analytic in $\Omega$}, if $f$ is real analytic at each
$y\in \Omega$.

Let $\Omega$ together with its boundary  be transformed pointwise
into the domain $\Omega'$ together with its boundary by equations of
the form
\begin{eqnarray} \label {02-001} x'_i= x_i+ f_i(x_1, \cdots,
x_{n}), \quad \; i=1,2, \cdots, n. \end{eqnarray}
 where the functions $f_i$ and their first order derivatives are
 Lipschitz continuous throughout the domain, and they are less in absolute
 value than a small positive number $\epsilon$.
 Then we say that the domain $\Omega$ is approximated by the domain $\Omega'$ with the
 degree of accuracy $\epsilon$.

 Let $(\mathcal {M}, g)$ be a Riemannian manifold. A subset
 $\Gamma$ of $(\mathcal{M},g)$ is said to be an
$(n-1)$-dimensional smooth (respectively,
  real analytic) surface if $\Gamma$ is nonempty and if for every
point $x$ in $\Gamma$, there is a smooth (respectively,  real
analytic) diffeomorphism of the open unit ball $B(0,1)$ in ${\Bbb
R}^n$ onto an open neighborhood $U$ of $x$ such that $B(0,1)\cap
\{x\in {\Bbb R}^n\big|x_n=0\}$ maps onto $U\cap \Gamma$.

  An $(n-1)$-dimensional surface $\Gamma$ in $(\mathcal {M},g)$ is said to
  be piecewise smooth (respectively, piecewise real analytic) if there exist
  a finite number of $(n-2)$-dimensional
smooth surfaces, by which $\Gamma$ can be divided into
  a finite number of  $(n-1)$-dimensional smooth (respectively, real analytic) surfaces.

\vskip 0.18 true cm

A subset $\mathfrak{F}$ of $L^2 (\Gamma)$ is called precompact if
any infinite sequence $\{u_k\}$ of elements of $\mathfrak{F}$
contains a Cauchy subsequence $\{u_{k'}\}$, i.e., one for which
\begin{eqnarray} \int_\Gamma (u_{k'}- u_{l'})^2 ds \to 0 \quad
\;\;\mbox{as}\;\; k', l'\to \infty.\end{eqnarray}

\vskip 0.16 true cm

 From here up to Section 5, let $\mathcal{M}$ be an
 $n$-dimensional Riemannian manifold with real analytic metric tensor $g$.

\vskip 0.25 true cm

 \noindent  {\bf Lemma 2.1.} \ \  {\it Let $D\subset (\mathcal{M}, g)$
 be a bounded domain with piecewise smooth boundary.
     Assume that $\mathfrak{M}$ is a set of functions $u$ in
 $H_0^1 (D) \cap  H^2(D)$ for which
  \begin{eqnarray} \int_D |\triangle_g u|^2 dR \end{eqnarray}
  is uniformly bounded. Then the set $\{\frac{\partial u}{\partial \nu} \big|u\in \mathfrak{M}\}$
  is precompact in $L^2(\partial D)$. }

\vskip 0.28 true cm

\noindent  {\bf Proof.}
 Put
 \begin{eqnarray} \label{2-13}\Lambda_1(D)= \inf_{u\in H_0^1(D)\cap  H^2(D)}
\;\frac{\int_D|\triangle_g u|^2 dR}{\int_D |\nabla_g
u|^2dR}.\end{eqnarray}
 We claim that $\Lambda_1(D)>0$.
 In fact, by applying Green's formula (see, for example, \cite{Ch1} or \cite{Sa})
  and Schwarz's inequality
  we see that for any $u\in H_0^1(D)\cap H^2(D)$,
  $$\left(\int_D |\nabla_g u|^2 dR\right)^2=\bigg|\int_D -u(\triangle_g u)dR
   \bigg|^2\le \left(\int_D u^2 dR\right)
  \left(\int_D |\triangle_g u|^2 dR\right), $$
i.e., \begin{eqnarray} \label{2-14} \frac{\int_D|\nabla_g u|^2 dR
}{\int_D|u|^2 dR}\le \frac{\int_D |\triangle_g u|^2 dR}{\int_D
|\nabla_g u|^2 dR},\end{eqnarray} where \begin{eqnarray*} \int_D
|\nabla_g u|^2 dR = \int_{D} g^{ik} (x)
 \frac{\partial u}{\partial x_i}\,\frac{\partial u}{\partial x_k}\, \sqrt{|g|}dx.\end{eqnarray*}
Since the first Dirichlet
 eigenvalue $\lambda_1 (D)$ is positive for the bounded domain
 $D$, $\,$
  i.e.,  \begin{eqnarray}\label{2-15}
 0<\lambda_1 (D)=\inf_{u\in H_0^1(D)}\frac{\int_D
|\nabla_g u|^2dR}{\int_D|u|^2dR},\end{eqnarray}
 we find by (\ref{2-14}) and (\ref{2-15}) that $\Lambda_1(D)>0$,
 and the claim is proved.

 From (\ref{2-15}) and (\ref{2-13}) we obtain that
\begin{eqnarray} \label {02.19} \int_D |u|^2 dR \le
  \frac{1}{\lambda_1 (D)}\int_D|\nabla_g u|^2 dR\quad \,
\mbox{for all}\;\;  u\in H_0^1(D)
  \end{eqnarray}
and
   \begin{eqnarray} \label {2-16} \qquad \quad\int_D |\nabla_g u|^2 dR \le
  \frac{1}{\Lambda_1(D)}\int_D|\triangle_g u|^2 dR\quad\,
 \mbox{for all}\;\;  u\in H_0^1(D)\cap  H^2(D).\end{eqnarray}

 Since $\partial D$ is piecewise smooth, we can write
  $\partial D=\cup_{i=1}^m \Gamma_i$, where $\Gamma_i$ is an
  $(n-1)$-dimensional surface. For each fixed $i$, $(i=1,\cdots, m)$,
 we choose a smooth $(n-1)$-dimensional surface $\Gamma'\subset\subset  D$
  such that $\partial \Gamma'_i=\partial \Gamma_i$ and
 $\Gamma_i\cup \Gamma'_i$ bounds an $n$-dimensional
 Lipschitz domain $D'_i$ satisfying
 $D'_i\subset\subset D\cup \Gamma_i$.
 Note that $u=0$ on $\Gamma_i$ for $u\in H_0^1(D)\cap H^2(D)$
 (see, for example,
p.$\,$62 of \cite{LM} or Corollary 6.2.43 of \cite{Ha}).
 It follows from the {\it a priori}
 estimate of the elliptic operators (see, for example,
  Theorem 9.13 of
 \cite{GT}) that there exists a constant $C_i>0$ depending only on $n, \Gamma_i,  D'_i$ and $D$
  such that
  \begin{eqnarray} \label {2-17} \| {u}\|_{H^2({D'}_i)}
  \le C_i(\|\triangle {u}\|_{L^2(D)}
  +\|{u}\|_{L^2(D)}).\end{eqnarray}
  By assumption, we have $\int_D |\triangle u|^2 dR\le {\tilde C}$ for all $u\in \mathfrak{M}$,
  where ${\tilde C}>0$ is a constant.
  According to (\ref{02.19}), (\ref{2-16}) and (\ref{2-17}), we
  see
  that for every $u\in {\mathfrak{M}}$,
\begin{eqnarray} \label {2-18} \| u\|_{H^2(D'_i)} \le C''_i,\end{eqnarray}
where $C''_i>0$ is a constant depending only on $n,\Gamma_i, D'_i$,
$D$ and $\tilde C$.
  Since  $D'_i$ is a domain with Lipschitz boundary in $(\mathcal {M}, g)$, it
   follows from the Neumann trace theorem (see, p.$\,$16 of \cite{AGMT},
  p.$\,$127 of \cite{Mc}, \cite{Gri} or Chs V, VI of \cite{EE}) that
     \begin{eqnarray*} \frac{\partial}{\partial \nu}\bigg|_{\Gamma_i} =\nu \cdot \nabla_g : \,\,
      \mathfrak{M}\to L^2(\Gamma_i)\end{eqnarray*}
     is precompact for each $i$ ($i=1,\cdots, m$).  Consequently, we obtain that
     $\{\frac{\partial u}{\partial \nu}\big| u\in \mathfrak{M}\}$
     is precompact in $L^2(\partial D)$.
     \quad $\square$

\vskip 0.25 true cm

 \noindent  {\bf Lemma 2.2.} \ \  {\it Let $(\mathcal{M},g)$ be a real
 analytic Riemannian manifold, and let  $D\subset (\mathcal{M}, g)$ be a
  bounded domain with piecewise
 smooth boundary. Suppose $\Gamma_1$ is a domain in $\partial D$
 with $\partial D-\bar \Gamma_1\ne \emptyset$ and assume that $\Gamma_2$ is
     an $(n-1)$-dimensional real analytic surface in $\partial D$
    satisfying $\bar \Gamma_2\subset \subset \partial D-\bar\Gamma_1$.
      Assume $\mathfrak{E}$ is a set of functions $u$ in
 $K^d(D)=\{u\big| u\in H^2(D),\, u=0 \;\;\mbox{on}\;\;\Gamma_1,\,
 \, u=\frac{\partial u}{\partial \nu}=0 \,\, \mbox{on}\;\; \Gamma_2\}$ for which
  \begin{eqnarray} \int_D |\triangle_g u|^2 dR \end{eqnarray}
  is uniformly bounded. Then the set $\{\frac{\partial u}{\partial \nu}\big|_{\Gamma_1} :
  u\in \mathfrak{E}\}$
  is precompact in $L^2(\Gamma_1)$. }

\vskip 0.28 true cm

\noindent  {\bf Proof.} \   Since $\partial D$ is  piecewise smooth,
it follows that $\Gamma_1$ can be divided into a finite number of
smooth $(n-1)$ dimensional surfaces. Without loss of generality, we
  let
 $\Gamma_1$ itself be a smooth $(n-1)$ dimensional surface. Put
 \begin{eqnarray} \label{2-19}\lambda_{\Gamma_1} (D)= \inf_{v\in K^d(D),\; \int_{D} |v|^2 dR =1}
 \; \frac{\int_{D}|\triangle_g v|^2 dR}{\int_{D} |v|^2 dR}.\end{eqnarray}
  In order to prove the existence of a minimizer to (\ref{2-19}),
  consider a minimizing sequence $v_m$ in $K^d(D)$,
 i.e.,
    \begin{eqnarray*}  \int_{D}|\triangle_g v_m|^2 dR\to
    \lambda_{\Gamma_1} (D)=0\quad \mbox{as}\;\, m\to +\infty\end{eqnarray*}
 with $\int_{D}  |v_m|^2 dR=1$.
  Then, there is a constance $C>0$ such that
  \begin{eqnarray} \label {2-21} \|\triangle_g v_m\|_{L^2(D)}\le C, \quad \, \|v_m\|_{L^2(D)}\le C \quad \;
  \mbox{for all}\;\; m.\end{eqnarray}
 Let $\{D_l\}$  be a sequence of Lipschitz domains such that
   $D_1\subset D_2\subset \cdots \subset D_l\subset \cdots \subset
   \subset D\cup \Gamma_1\cup \Gamma_2$, $\;\cup_{l=1}^\infty D_l= D$,
   and $\Gamma_1\cup \Gamma_2 \subset \partial D_l$ for
   all $l$.
It follows from
 the {\it a priori} estimate for elliptic
equations  (see, for example, Theorem 9.13 of
 \cite{GT})  that there exists a constant $C'_l>0$ depending only on
 $n, D_l, D, \Gamma_1$ and $\Gamma_2,$
  such that
  \begin{eqnarray} \label {2-22} \| v_m\|_{H^2(D_l)}
  \le C'_l(\|\triangle_g  v_m\|_{L^2(D)}
  +\|v_m\|_{L^2(D)}).\end{eqnarray}
From this and (\ref{2-21}), we see  that
\begin{eqnarray*} \|v_m\|_{H^2(D_l)}\le C''_l \;\;\mbox{for all}\;\; m,\end{eqnarray*}
where $C''_l$ is a constant depending only on $n, D_l, D, \Gamma_1,
\Gamma_2$ and $C$. For each $l$, by the Banach-Alaoglu theorem we
can then extract a subsequence $\{v_{l,m}\}_{m=1}^\infty$ of
$\{v_m\}$, which converges weakly in $H^2(D_l)$ to a limit $u$, and
converges strongly in $L^2(D_l)$ to $u$. We may assume that
$\{v_{l+1,m}\}_{m=1}^\infty$ is a subsequence of
$\{v_{l,m}\}_{m=1}^\infty$ for every $l$. Then, the diagonal
sequence $\{v_{l,l}\}_{l=1}^\infty$ converges weakly in $H^2$ to
$u$, and strongly converges to $u$ in $L^2$, in every compact subset
$E$ of $D$.  It is obvious that $\|u\|_{L^2 (D)}=1$. Since the
functional $\int_{D_l} |\triangle_g u|^2 dR$ is lower semicontinuous
in the weak $H^2(D_l)$ topology, we have
$$\int_{D_l} |\triangle_g u|^2 dR\le
 \underset {k\to \infty} {\underline{\lim}} \int_{D_l} |\triangle_g v_{k,k}|^2 dR,$$
 so that
 \begin{eqnarray*} \int_D |\triangle_g u|^2 dR &=&\lim_{l\to \infty} \int_{D_l}
 |\triangle_g u|^2 dR \le \lim_{l\to \infty} \left(\underset {k\to
 \infty}{\underline{\lim}}
 \int_{D_l} |\triangle_g v_{k,k}|^2 dR\right)\\
 &\le& \lim_{l\to \infty} \left(\underset {k\to
 \infty}{\underline{\lim}}
 \int_{D} |\triangle_g v_{k,k}|^2 dR\right)
 =\lambda_{\Gamma_1}(D).\end{eqnarray*}
For each fixed $l$, since $v_{k,k} \to u$ weakly in $H^2(D_l)$, we
get that $v_{k,k}\to u$ strongly in $H^r (D_l)$ for any $0<r<2$.
Note that $\frac{\partial v_{k,k}}{\partial \nu}\big|_{\Gamma_2}=0$
and $v_{k,k}\big|_{\Gamma_i}=0$ for $i=1,2$. It follows that
$u\big|_{\Gamma_1}=0$ and $u\big|_{\Gamma_2}=\frac{\partial
u}{\partial \nu}\big|_{\Gamma_2}=0$.
  Therefore $u\in K^d(D)$ is a minimizer.

 We claim that  $\lambda_{\Gamma_1}(D)>0$.
Suppose by contradiction that $\lambda_{\Gamma_1}(D)=
\frac{\int_{D}|\triangle_g u|^2 dR}{\int_{D} |u|^2 dR}= 0$. Then
$\triangle_g u=0$ in $D$.
 Since the coefficients of the Laplacian are real analytic in $D$,
  and since $\Gamma_2$ is a real analytic surface, we find with the aid of the
  regularity for elliptic equations (see,
  Theorem A of \cite{MN}, \cite{Mo} or \cite{ADN}) that $u$ is real analytic
  up to the partial boundary $\Gamma_2$.
 Note that $u=\frac{\partial
  u}{\partial \nu}=0$ on $\Gamma_2$.
 Applying Holmgren's uniqueness theorem (see, Corollary 5 of p.$\,$39 in \cite{Ra}
 or p.$\,$433 of \cite{Ta1})
  for the real analytic
 elliptic equation $\triangle_g u=0$ in $D$, we get $u\equiv 0$ in
 $D$. This
   contradicts the fact $\int_D |u|^2dR =1$, and the claim is proved.
Therefore we have
    \begin{eqnarray} \label {2-23} \int_{D} |u|^2 dR \le
  \frac{1}{\lambda_{\Gamma_1}(D)}\int_{D}|\triangle_g u|^2 dR
  \quad \;\mbox{for}\;\; u\in K^d(D).\end{eqnarray}
According to the assumption,
  there is a constant $C''$ such that  \begin{eqnarray}
  \label{2-24}\|\triangle_g u
  \|_{L^2(D)}  \le C''\quad \; \mbox{for all}\;\; u \in
  \mathfrak{E}.\end{eqnarray}
  Again, applying the {\it a priori} estimate  for the
  elliptic equations in some (fixed) subdomain
  $D_l\subset\subset D\cup \Gamma_1\cup \Gamma_2$ (see, Theorem 9.13 of
  \cite{GT}), we obtain that
  \begin{eqnarray} \label {2-26} \| u\|_{H^2(D_l)}
  \le C'_l (\|\triangle_g  u\|_{L^2(D)}
  +\|u\|_{L^2(D)}),\end{eqnarray}
  where the constant $C'_l$ is as in (\ref{2-22}). By (\ref{2-23})---(\ref{2-26}), we get
  that for every $u\in {\mathfrak{E}}$,
  \begin{eqnarray*} \label {2-27} \|u\|_{H^2(D_l)} \le C''',\end{eqnarray*}
  where $C'''>0$ is a constant depending only on $n,D_l, D, \Gamma_1$,
  $\Gamma_2$ and $C''$.
  It follows from the Neumann trace theorem (see, p.$\,$16 of
  \cite{AGMT},
   \cite{Gri} or \cite{MMS} ) that
  $\{\frac{\partial u}{\partial \nu}\big|_{\Gamma_1} : u\in {\mathfrak{E}}\}$
  is precompact in $L^2(\Gamma_1)$. \ \ $\square$

\vskip 0.32 true cm

The following two results will be needed later:

\vskip 0.2 true cm

 \noindent  {\bf Proposition 2.3 (see, p.$\;$12 of \cite{Sa}).} \ \  {\it
 Let $\Pi^0$ be an  isometric transformation which maps a Hilbert
 space ${\mathcal{H}}^0$ onto a subspace $\Pi^0{\mathcal{H}}^0$ of
 another Hilbert space $\mathcal{H}$, so that
 \begin{eqnarray*} \langle u^0, v^0\rangle^0 =\langle \Pi^0 u^0,
 \Pi^0 v^0\rangle \quad \;\; \mbox{for all}\;\; u^0, v^0\in
 {\mathcal{H}}^0.\end{eqnarray*}
 Suppose that $G^0$ and $G$ are two non-negative, self-adjoint, completely
 continuous linear transformations on ${\mathcal{H}}^0$ and
 $\mathcal{H}$ respectively, such that
 \begin{eqnarray*}
  \langle G^0 u^0, v^0\rangle^0 = \langle G\Pi^0 u^0, \Pi^0 v^0\rangle
  \quad \;\mbox{for all}\;\; u^0, v^0\in
  {\mathcal{H}}^0.\end{eqnarray*}
  Then
  \begin{eqnarray*}  \mu_k^0 \le \mu_k \quad \;\; \mbox{for}\;\;
  k=1,2,3,\cdots,\end{eqnarray*}
 where $\{\mu_k^0\}$ and $\{\mu_k\}$ are the
 eigenvalues of $G^0$ and $G$, respectively.}

\vskip 0.25 true cm

\noindent  {\bf Proposition 2.4 (see, p.$\;$13 of \cite{Sa}).} \ \
{\it
  Assume that $\mathcal{H}$ is a direct sum of $p$ Hilbert spaces
  ${\mathcal{H}}_j$
  \begin{eqnarray*}
  \mathcal{H}={\mathcal{H}}_1\oplus{\mathcal{H}}_2\oplus
 \cdots  \oplus {\mathcal{H}}_p\end{eqnarray*}
 and that the self-adjoint, completely continuous linear
 transformation $G$ maps every ${\mathcal{H}}_j$ into itself,
 $$G{\mathcal{H}}_j\subset {\mathcal{H}}_j, \quad \;\; j=1,2,3,
 \cdots,p.$$
 Denote by $G_j$ the restriction of $G$ to ${\mathcal{H}}_j$.
 Then the set of eigenvalues of the transformation $G$ (each
 eigenvalue repeated according to its multiplicity) is identical to
 the union of the sets of eigenvalues of $G_1, \cdots, D_p$.}

\vskip 1.39 true cm

\section{Completely continuous transformations and eigenvalues}

\vskip 0.45 true cm

  Let $(\mathcal{M}, g)$ be an $n$-dimensional real analytic Riemannian
  manifold and let $D\subset \mathcal{M}$ be a bounded
  domain with piecewise smooth boundary $\Gamma$. Suppose that $\varrho$
  is a non-negative bounded function defined on
  $\Gamma$ or only on a portion $\Gamma_\varrho$ of $\Gamma$
  (measure $\Gamma_\varrho=\int_{\Gamma_\varrho}
  ds >0$) and assume that
  $\int_{\Gamma_\varrho}\varrho\,ds>0$.
  In  case $\Gamma_\varrho \ne \Gamma$ we denote $\Gamma_0=\Gamma-\bar \Gamma_\varrho$, and
  assume that  $\Gamma_{00}$ is a  real analytic $(n-1)$-dimensional surface in $\Gamma_{0}$.

  If
$\Gamma_\varrho \ne \Gamma$ (measure $\Gamma_0>0$), we denote
\begin{eqnarray*} & K(D)=\{u\big|u\in H_0^1(D) \cap  H^2(D),\,\,\mbox{and}\;\;
\frac{\partial u}{\partial \nu}=0 \;\; \mbox{on}\;\;
\Gamma_{00}\}\\
&  K^d(D)=\{u\big|u\in  H^2(D),\,\, u=0\;\; \mbox{on}\;\;
\Gamma_\varrho, \,\,
 \mbox{and}\,\, u=\frac{\partial u}{\partial \nu}=0 \;\; \mbox{on}\;\; \Gamma_{00}\}.\end{eqnarray*}

If $\Gamma_\varrho =\Gamma$, we denote
\begin{eqnarray*} N(D)=\{u\big|u\in H_0^1(D)\cap  H^2(D)\}.\end{eqnarray*}

It follows from the property of $H_0^1(\Omega)$  (see, for example,
p.$\,$62 of \cite{LM} or Corollary 6.2.43 of \cite{Ha} or \cite{Mc})
that $u=0$ on $\partial D$ for any $u\in H^1_0(\Omega)$ (Therefore,
we always have that $u=0$ on $\Gamma$ for any $u\in K(D)$ or
$N(D)$).

 We shall also use the notation
\begin{eqnarray*}    \langle u,v\rangle^{\star} =\int_D (\triangle_g u)(\triangle_g v) dR,
 \quad \, u, v \in K(D) \;\; \mbox{or}\;\; K^d(D) \;\;\mbox{or}\;\;
 N(D).\end{eqnarray*}
 The bilinear functional $\langle u,v\rangle^{\star}$ can be used
as an inner product in each of the spaces $K(D)$, $K^d(D)$ and
$N(D)$.
  In fact, $\langle u,v\rangle^{\star}$ is a positive, symmetric, bilinear functional.
  In addition, if $\langle u,u\rangle^{\star}=0$, then $\triangle_g u=0$ in $D$. In the case
  $u\in K(D)$ or $N(D)$, by
   applying
 the maximum principle, we have
 $u\equiv 0$ in $D$. In the case  $u\in K^d (D)$, since $u=\frac{\partial
 u}{\partial \nu}=0$ on $\Gamma_{00}$, we find by Holmgren's uniqueness theorem
 (see, Corollary 5 of p.$\,$39 in \cite{Ra}) that
  $u\equiv 0$ in $D$.
  Closing $K(D)$, $K^d(D)$ and $N(D)$ with respect to  the norm
  $\|u\|^\star=\sqrt{\langle u, u\rangle^{\star}}$,
  we get the Hilbert spaces $(\mathcal{K}, \|\cdot\|^\star)$,
   $({\mathcal{K}}^d, \|\cdot\|^\star)$
  and $(\mathcal{N}, \|\cdot\|^\star)$, respectively.

  Next, we consider two linear
functionals
\begin{eqnarray*} [u,v]=\int_{\Gamma_\varrho} \varrho \,\frac{\partial u}{\partial \nu}
\, \frac{\partial v}{\partial \nu}
  \, ds
\end{eqnarray*}
and \begin{eqnarray} \label{3.0}\langle u,v\rangle
 =\langle u,v\rangle^{\star}  +[u,v], \end{eqnarray} where $u,v\in K(D)$ or
 $u,v\in K^d(D)$ or $u,v\in N(D)$.  It is clear that $\langle
u,v\rangle$ is an inner product in each of the spaces $K(D)$, $K^d
(D)$ and $N(D)$.

\vskip 0.25 true cm

 \noindent  {\bf Lemma 3.1.} \ \  {\it The norm
 \begin{eqnarray*}  \|u\|^\star=\sqrt{\langle u,u\rangle^{\star} }\end{eqnarray*}
 and \begin{eqnarray*} \|u\|=\sqrt{\langle u,u\rangle}\end{eqnarray*}
are equivalent in $K(D)$, $K^d(D)$ and $N(D)$.}

\vskip 0.2 true cm

 \noindent  {\bf Proof.} \ \  Obviously, $\|u\|^\star\le \|u\|$ for all $u$
 in $K(D)$ or $K^d(D)$ or $N(D)$.
  In order to prove the equivalence of the two norms,
 we first consider the case in linear space $N(D)$.
 It suffices to show that
$\|u\|$ is bounded when $u$ belongs to the set
\begin{eqnarray*} \mathfrak{M}=\{u\big|u\in N(D), \|u\|^\star\le
1\}.\end{eqnarray*}  It follows from Lemma 2.1 that
${\mathfrak{M}}_\Gamma:=\{\frac{\partial u}{\partial \nu}\big| u\in
\mathfrak{M}\}$ is precompact in $L^2(\Gamma)$. This implies that
there exists a constant $C>0$ such that $\int_\Gamma
\big(\frac{\partial u}{\partial \nu}\big)^2 ds \le C$ for all $u\in
\mathfrak{M}$. Therefore, $[u,u]=\int_{\Gamma} \varrho
\left(\frac{\partial u}{\partial \nu}\right)^2 ds$ is bounded in
$\mathfrak{M}$, and so is $\|u\|^2=\langle u,u\rangle^{\star}
+[u,u]$.
 Similarly, applying Lemmas 2.1, 2.2 we can prove the
 corresponding results for the spaces $K(D)$ and $K^d(D)$.
  \ \  $\square$

\vskip 0.22 true cm

From Lemmas 2.1, 2.2, it follows that
\begin{eqnarray*}  \qquad  |[u,u]|=\bigg|\int_{\Gamma_\varrho} \varrho
\left(\frac{\partial u}{\partial \nu}\right)^2  ds \bigg|\le C
\langle u,u\rangle^{\star} \; \; \mbox{for all}\;\; u\;\;
\mbox{in}\;\; K(D) \;\;\mbox{or}\;\; K^d(D) \;\; \mbox{or}\;\;
N(D).\end{eqnarray*} Therefore, $[u,v]$  is a  bounded,  symmetric,
  bilinear functional in
$(K(D)$,$\langle\cdot$,$\cdot\rangle^{\star})$,
$(K^d(D)$,$\langle\cdot$,$\cdot\rangle^{\star})$ and $(N(D)$,
$\langle \cdot$,$\cdot\rangle^{\star})$. Since it is densely defined
in $(\mathcal{K}$,$\langle\cdot,\cdot\rangle^{\star})$,
$({\mathcal{K}}^d$, $\langle\cdot,\cdot\rangle^{\star})$ and
$\,(\mathcal{N}, \langle\cdot,\cdot\rangle^{\star})$, respectively,\
it can immediately be extended to
$\,(\mathcal{K}$,$\langle\cdot,\cdot\rangle^{\star})$,
$({\mathcal{K}}^d$,$\langle\cdot,\cdot\rangle^{\star})$ and
$(\mathcal{N}, \langle\cdot,\cdot\rangle^{\star})$. We still use
$[u,v]$ to express the extended functional. Then there is a bounded
linear transformation $G_{\mathcal{K}}^{(\star)}$ of
$(\mathcal{K},\langle\cdot,\cdot\rangle^{\star})$ into
$(\mathcal{K},\langle\cdot,\cdot\rangle^{\star})$
 (respectively, $G_{{\mathcal{K}}^d}^{(\star)}$ of
 $({\mathcal{K}}^d,\langle\cdot,\cdot\rangle^{\star})$
into $({\mathcal{K}}^d,\langle\cdot,\cdot\rangle^{\star})$,
    $G_{\mathcal{N}}^{(\star)}$ of
$(\mathcal{N},\langle\cdot,\cdot\rangle^{\star})$ into
$(\mathcal{N},\langle\cdot,\cdot\rangle^{\star})$) such that
   \begin{eqnarray} \label {3-1} [u,v]=\langle G_{\mathcal{K}}^{(\star)}
    u, v\rangle^{\star} \quad \,\;
   \mbox{for all}\;\; u \;\; \mbox{and}\;\; v\;\;\mbox{in}\;\;
   \mathcal{K} \end{eqnarray}
  (respectively, \begin{eqnarray} \label {3.1} [u,v]=
  \langle G_{{\mathcal{K}}^d}^{(\star)} u, v\rangle^{\star} \quad \,\;
   \mbox{for all}\;\; u \;\; \mbox{and}\;\; v\;\;\mbox{in}\;\;
   {\mathcal{K}}^d, \end{eqnarray}
    \begin{eqnarray} \label{3.2} [u,v]= \langle G_{\mathcal{N}}^{(\star)}
    u, v\rangle^{\star} \quad \;\,
   \mbox{for all}\;\; u \;\; \mbox{and}\;\; v\;\;\mbox{in}\;\;
   \mathcal{N}).\end{eqnarray}

 \vskip 0.25 true cm

 \noindent  {\bf Lemma 3.2.} \ \  {\it The transformations
 $G_{\mathcal{K}}^{(\star)}$, $G_{{\mathcal{K}}^d}^{(\star)}$
 and $G_{\mathcal{N}}^{(\star)}$ are self-adjoint and compact.}

\vskip 0.2 true cm

 \noindent  {\bf Proof.} \  Since $[u,v]$ is symmetric, we immediately get
 that the transformation $G_{\mathcal{K}}^{(\star)}$, $G_{{\mathcal{K}}^d}^{(\star)}$ and
$G_{\mathcal{N}}^{(\star)}$ are all self-adjoint. For the
compactness,
 we only discuss the case for the transformation
$G_{\mathcal{K}}^{(\star)}$. It suffices to
  show (see, p.$\,$204 of \cite{RN}):

  From every sequence $\{u_m\}$ in $K(D)$ which
is bounded
\begin{eqnarray} \label{3-4} \|u_m\|^\star \le constant, \; \, m=1,2, 3, \cdots, \end{eqnarray}
   we can pick out a subsequence $\{u_{m'}\}$ such that
   \begin{eqnarray} \label{3-5} \langle G_{\mathcal{K}}^{(\star)} (u_{m'}-u_{l'}),
   (u_{m'}-u_{l'})\rangle^{\star}
   \to 0 \,
   \, \mbox{when}\;\; m', l'\to \infty.\end{eqnarray}

   Applying Lemma 2.1 with the aid of (\ref{3-4}), we find that the
   sequence $\{\frac{\partial u_m}{\partial \nu}\}$ is
   precompact in $L^2(\Gamma_\varrho)$, so that
  there is a subsequence $\{u_{m'}\}$ such that
  \begin{eqnarray*} \int_{\Gamma_\varrho} \left(\frac{\partial (u_{m'}
  -u_{l'})}{\partial \nu}\right)^2 ds \to 0 \;\;
  \mbox{as} \;\;
    m', l'\to \infty. \end{eqnarray*}
   Therefore
 \begin{eqnarray*} [u_{m'}- u_{l'}, u_{m'} -u_{l'}] =
 \int_{\Gamma_\varrho} \varrho \left(\frac{\partial (u_{m'}
  -u_{l'})}{\partial \nu}\right)^2 ds \to 0 \;\;
  \mbox{as} \;\;
    m', l'\to \infty, \end{eqnarray*}
which implies (\ref{3-5}). This proves the compactness of
$G_{\mathcal{K}}^{(\star)}$.
   \ \  $\square$

\vskip  0.25 true cm

Except for the transformations $G^{(\star)}_{\mathcal{K}}$,
$G^{(\star)}_{{\mathcal {K}}^d}$ and $G^{(\star)}_{\mathcal{N}}$,
 we need introduce corresponding transformations
 $G_{\mathcal{K}}$,
  $G_{{\mathcal{K}}^d}$
   and  $G_{\mathcal{N}}$
  by the inner product $\langle \cdot, \cdot \rangle$.
  Since
\begin{eqnarray} 0\le [u,u]\le \langle u,u\rangle \quad \;\; \mbox{for all } u
\;\;\mbox{in}\;\; K(D)\;\;\mbox{or}\;\; K^d(D) \;\; \mbox{or}\;\;
 N(D),\end{eqnarray}  there is a bounded linear self-adjoint
transformation $G_{\mathcal{K}}$ of $(\mathcal{K}, \langle \cdot,
\cdot\rangle)$
 (respectively,  $G_{{\mathcal{K}}^d}$ of $({\mathcal{K}}^d, \langle \cdot,
\cdot\rangle)$,  $\,G_{\mathcal{N}}$ of $(\mathcal{N}, \langle
\cdot, \cdot \rangle)$)
 such that
\begin{eqnarray} \label{3-6} [u,v]=\langle G_{\mathcal{K}} u,v\rangle
\quad \;\; \mbox{for all} \;\; u\;\; \mbox{and}\;\; v\;\;
\mbox{in}\;\; \mathcal{K}\end{eqnarray}
 (respectively,
\begin{eqnarray} \label{3-7} [u,v]=\langle G_{{\mathcal{K}}^d}u,v\rangle
\quad \;\; \mbox{for all} \;\; u\;\; \mbox{and}\;\; v\;\;
\mbox{in}\;\; {\mathcal{K}}^d, \end{eqnarray}   \begin{eqnarray}
\label{3-8} [u,v]=\langle G_{\mathcal{N}}u,v\rangle
 \quad \;\; \mbox{for all} \;\; u\;\;
\mbox{and}\;\; v\;\; \mbox{in}\;\; \mathcal{N}).\end{eqnarray}

\vskip 0.25 true cm

 \noindent  {\bf Lemma 3.3.} \ \  {\it The transformations
 $G_{\mathcal{K}}$, $G_{{\mathcal{K}}^d}$
 and  $G_{\mathcal{N}}$ are positive and compact.}

\vskip 0.2 true cm

\noindent  {\bf Proof.} \ From $[u,u]\ge 0$  for any $u\in
{\mathcal{K}} $ or ${\mathcal{K}}^d$ or
 ${\mathcal{N}}$, we immediately know that
 $G_{\mathcal{K}}$, $G_{{\mathcal{K}}^d}$
 and  $G_{\mathcal{N}}$ are positive. The proof of the compactness is completely
 similar to that of Lemma 3.2. \ \   $\square$

 It follows from Lemma 3.3 that
 $G_{\mathcal{K}}$ (respectively, $G_{{\mathcal{K}}^d}$,  $G_{\mathcal{N}}$)
  has only non-negative eigenvalues and that the positive
 eigenvalues form an enumerable sequence  $\{\mu_{\mathcal{K}}\}$
 (respectively, $\{\mu_{{\mathcal{K}^d}}\}$,
   $\{\mu_{\mathcal{N}}\}$) with $0$ as the
 only limit point.

\vskip 0.25 true cm

 \noindent  {\bf Theorem 3.4.} \ \  {\it The transformations
 $G_{\mathcal{K}}^{(\star)}$ and
  $G_{\mathcal{K}}$ (respectively,  $G_{{\mathcal{K}}^d}^{(\star)}$ and
  $G_{{\mathcal{K}}^d}$,  $\,G_{\mathcal{N}}^{(\star)}$ and
  $G_{\mathcal{N}}$) have the same eigenfunctions. If $\mu_{\mathcal{K}}^{\star}$
 and $\mu_{\mathcal{K}}$ (respectively,
  $\mu_{{\mathcal{K}}^d}^{\star}$
 and $\mu_{{\mathcal{K}}^d}$,
    $\, \mu_{\mathcal{N}}^{\star}$
 and $\mu_{\mathcal{N}}$)
 are eigenvalues corresponding to
 the same eigenfunction we have
 \begin{eqnarray} \label{3.10}\mu_{\mathcal{K}}=
 \frac{\mu_{\mathcal{K}}^{\star}}{1+\mu_{\mathcal{K}}^{\star}}\end{eqnarray}
(respectively,   \begin{eqnarray} \mu_{{\mathcal{K}}^d}=
 \frac{\mu_{{\mathcal{K}}^d}^{\star}}{1+\mu_{{\mathcal{K}}^d}^{\star}}, \end{eqnarray}
   \begin{eqnarray} \label {3.12}\mu_{\mathcal{N}}=
 \frac{\mu_{\mathcal{N}}^{\star}}{1+\mu_{\mathcal{N}}^{\star}}).\end{eqnarray}
 }

\vskip 0.2 true cm

 \noindent  {\bf Proof.} \  We only prove the case
 for the $G_{\mathcal{K}}$
 (a similar argument will work for $G_{{\mathcal{K}}^d}$
 and $G_{\mathcal{N}}$). Since
  $G_{\mathcal{K}}^{(\star)}$ is positive, we can easily conclude that the inverse
 $(1+G_{\mathcal{K}}^{(\star)})^{-1}$ exists and is a bounded self-adjoint
 transformation.
   By virtue of  (\ref{3-1}), (\ref{3-6}) and (\ref{3.0}), we have
\begin{eqnarray} \label{3-9}  \langle G_{\mathcal{K}}^{(\star)} u, v\rangle^{\star}
 &=&[u,v] =  \langle G_{\mathcal{K}}u, v\rangle \\
 &=&  \langle G_{\mathcal{K}}u, v\rangle^{\star}
   + \langle G_{\mathcal {K}} G_{\mathcal{K}}^{(\star)} u, v\rangle^{\star}, \quad \;
 (u,v\in \mathcal{K}).\nonumber \end{eqnarray}
 It follows that
 \begin{eqnarray} \label {3-10} G_{\mathcal{K}}=
  G_{\mathcal{K}}^{(\star)}(1+G_{\mathcal{K}}^{(\star)})^{-1},\end{eqnarray}
   from which the desired result follows immediately.
\ $\square$

\vskip 0.25 true cm

 \noindent  {\bf Proposition 3.5.} \   {\it Let $u$ and $v$ be two eigenfunctions in
 $({\mathcal{K}},\langle \cdot, \cdot\rangle)$ (respectively,
 $({\mathcal{K}}^d,\langle \cdot, \cdot\rangle)$,
   $({\mathcal{N}},\langle \cdot, \cdot\rangle)$)
 of the transformation $G_{\mathcal{K}}$ (respectively,
 $G_{{\mathcal{K}}^d}$,
  $G_{\mathcal{N}}$) at least one of which corresponds to a non-vanishing eigenvalue.
 Then $u$ and $v$ are orthogonal if and only if the
 $\frac{\partial u}{\partial \nu}\big|_{\Gamma_\varrho}$
 and $\frac{\partial v}{\partial\nu}\big|_{\Gamma_\varrho}$
 are orthogonal in $L_\varrho^2 (\Gamma_\varrho)$, that is,
 \begin{eqnarray} [u,v]=\int_{\Gamma_\varrho} \varrho
 \frac{\partial u}{\partial\nu} \,\frac{\partial v}{\partial \nu}\,ds=0.\end{eqnarray}}

 \vskip 0.15 true cm

 \noindent  {\bf Proof.} \ \  Without loss of generality,
 we suppose that $u$ is the eigenfunction corresponding
 to the eigenvalue $\mu \ne 0$. Then
 \begin{eqnarray*} [u,v] = \langle G_{\mathcal{K}}u, v\rangle
 =\mu
 \langle u, v\rangle, \end{eqnarray*}
   which implies the desired result.
    $\;\; \square$

\vskip 0.25 true cm

We can now prove  \vskip 0.20 true cm

 \noindent  {\bf Theorem 3.6.} \ \  {\it  Let  $D\subset (\mathcal{M},g)$
  be a bounded domain with piecewise smooth boundary $\Gamma$.
 Assume that $\Gamma_{00}$ is an $(n-1)$-dimensional
 surface in $\Gamma-\bar \Gamma_\varrho$.
   If $u$ is an eigenfunction of
 the transformations $G_{\mathcal{K}}^{(\star)}$
  or $G_{\mathcal{N}}^{(\star)}$ with eigenvalue $\mu^\star \ne 0$,
  then
  $u$ has derivatives of any order in $D$ and is such that
  \begin{eqnarray}  \label{3-15} \left\{\begin{array}
  {ll}\triangle_g^2 u=0 \quad \mbox{in}\;\; D,  \\
   u=0 \; \; \mbox{on}\;\; \Gamma,\\
 \frac{\partial u}{\partial \nu}=0 \;\; \mbox{on}\;\;
 \Gamma_{00},\;\;  \Delta_g u=0 \;\; \mbox{on}\;\;
 \Gamma-(\Gamma_\varrho\cup \Gamma_{00}),\\
    \triangle_g u+ \gamma \, \varrho \,
 \frac{\partial u}{\partial \nu}=0 \;\; \mbox{on}\;\; \Gamma_\varrho,
 \; \quad \; \mbox{with}\;\; \gamma=\frac{1}{\mu^\star}. \end{array}\right.\end{eqnarray}}

\vskip 0.2 true cm

 \noindent  {\bf Proof.} \ \  Let $\{u_j\}$ be a sequence
 of functions in $K(D)$ such that $\| u_j-u \|^\star \to 0$ as $j\to \infty$.
   We first claim that \begin{eqnarray} \label{3-18}  u_j\to u\;\; \mbox{in}\;\;
   L^2(D).\end{eqnarray}
    In fact, since $u_j\in H_0^1(D)\cap  H^2(D)$,
    it follows from (\ref{02.19}) and (\ref{2-16})
 that \begin{eqnarray} \label{3-19}
  \int_{D} |u_j-u_l|^2 dR \le \frac{1}{\lambda_1(D)}
  \int_D |\nabla_g (u_j-u_l)|^2 dR \quad \, \mbox{for any}\;\;
  j\;\; \mbox{and}\;\;l,\end{eqnarray}
    \begin{eqnarray} \label{3-20}\qquad \;\; \int_D
  |\nabla_g (u_j- u_l)|^2 dR \le \frac{1}{\Lambda_1 (D)}
 \int_D |\triangle_g (u_j-u_l)|^2 dR \quad \, \mbox{for any}\;\;
  j\,\,\mbox{and}\,\,l,
 \end{eqnarray}
where  $\lambda_1(D)$ and $\Lambda_1(D)$ are the first
  Dirichlet and buckling eigenvalues for $D$, respectively.
   Since $\int_D
  |\triangle_g (u_j-u_l)|^2 dR\to 0\;\; \mbox{as}\;\; j, l\to +\infty$,
  we find by (\ref{3-19}) and (\ref{3-20}) that
$\int_D
  |u_j-u_l|^2 dR \to 0\;\; \mbox{as}\;\; j, l\to +\infty$.
  Therefore the claim is proved.

 For any point $p$ in $D$,
  let $U$ be a coordinate neighborhood of $p$, and let $E\ni p$ be a
  bounded domain with smooth boundary such that
$\bar E \subset U\cap D$. Let $f$ be a function
  in $C_0^4 (E)$. Then, by Green's formula (see, for example, p.$\,$6 of \cite{Ch1}), we have
  \begin{eqnarray*}  \label {3-21} \langle u_n, f\rangle^\star =
  \int_E (\triangle_g u_n)(\triangle_g f) dR = \int_E
  u_n (\triangle_g^2 f) dR,\end{eqnarray*}
   so that
\begin{eqnarray}  \label {3-22}  \langle u, f\rangle^\star =\int_E
  u (\triangle_g^2 f) dR.\end{eqnarray}
Now by assumption, $G_{\mathcal{K}}^{(\star)} u=\mu^{\star} u$
 with $\mu^{\star} \ne 0$ and $\frac{\partial f}{\partial \nu}=0$ on $\Gamma$, and hence we have
\begin{eqnarray}  \label {3-23} \mu^{\star} \langle u, f\rangle^{\star}  =
     \langle G_{\mathcal{K}}^{(\star)} u, f\rangle^{\star}
   =[u,f]=0.\end{eqnarray}
 Since $p$ is arbitrary in $D$, it follows from (\ref{3-22}) and
 (\ref{3-23}) that
  \begin{eqnarray*} \int_D u (\triangle^2_g f) dR =0 \quad \;\;
  \mbox{for all}\;\; f\in C_0^4(D).\end{eqnarray*}
 By applying Green's formula again, we get
 \begin{eqnarray} \int_D (\triangle_g u) (\triangle_g f) \, dR =0 \quad \;\;
  \mbox{for all}\;\; f\in C_0^4(D),\end{eqnarray}
 i.e., $u$ is a weak solution of $\triangle_g^2 u=0$ in $D$ (see \cite{GT}).
 It follows from the interior regularity of elliptic equations that $u\in C^\infty (D)$,
 and in the classic sense
   \begin{eqnarray}\label{3-24}  \triangle _g^2 u=0 \;\; \mbox{in}\;\; D.\end{eqnarray}
 In exactly the same way, the corresponding result can be proved for $G_{\mathcal{N}}^{(\star)}$.

 Next, suppose
 that  $\varrho$ is continuous.
  That the boundary conditions of (\ref{3-15}) hold follows from Lemma 2.1 and Green's
formula. In fact, if \begin{eqnarray*} G_{\mathcal{K}}^{(\star)} u=
\mu^{\star} u,\end{eqnarray*} then $u\big|_{\Gamma}=0$  and
$\frac{\partial u}{\partial \nu}\big|_{\Gamma_{00}}=0$, and that
 \begin{eqnarray*}  \label{03-03}\int_{\Gamma_\varrho} \varrho \frac{\partial u}{\partial \nu}\,
 \frac{\partial v}{\partial \nu} \, ds= \mu^{\star} \int_D
 (\triangle_g u)(\triangle v) dR \quad \, \mbox{for all} \;\;
 v\in {\mathcal{K}} (D). \end{eqnarray*}
 By this  and Green's formula (see, p.$\,$114-120 of \cite{LM}, \cite{Lo} and \cite{Fe}),
 we obtain that
\begin{eqnarray*} \frac{1}{\mu^{\star}} \int_{\Gamma_\varrho} \varrho \frac{\partial u}{\partial \nu}\,
 \frac{\partial v}{\partial \nu} \, ds=
  \int_D
 (\triangle_g^2 u) v\, dR - \int_{\Gamma} (\triangle_g u)\frac{\partial v}{\partial \nu}\,ds
 +\int_{\Gamma} \frac{\partial (\triangle_g u)}{\partial \nu} \,v \, ds. \end{eqnarray*}
 for all $v\in {\mathcal {K}}(D)$,
where $\frac{\partial (\triangle_g u)}{\partial \nu}\in H^{-3/2}
(\Gamma)$ (see \cite{AGMT}).  Thus
 \begin{eqnarray}   && \int_D  (\triangle_g^2 u)v\, dR -
 \int_{\Gamma_\varrho} \left(\triangle_g u +\frac{1}{\mu^{\star}}
  \varrho \frac{\partial u}{\partial \nu}\right)  \frac{\partial
  v}{\partial \nu} \, ds  \\
  && \quad \;\; + \int_{\Gamma -(\Gamma_\varrho\cup
  \Gamma_{00})} (\Delta_g u)\frac{\partial v}{\partial \nu}\, ds +
\int_\Gamma  \frac{\partial (\triangle_g u)}{\partial \nu} \;v\, ds
  =0\end{eqnarray}
  for all $v\in {\mathcal K}(D)$.
 Note that
   $v\big|_{\Gamma}=0$  and $\frac{\partial v}{\partial \nu}\big|_{\Gamma_{00}}=0$,
 and that  $\frac{\partial v}{\partial \nu}\big|_{\Gamma_{\varrho}}$
 and $\frac{\partial v}{\partial \nu}\big|_{\Gamma-(\Gamma_\varrho\cup \Gamma_{00})}$
 run throughout
   space $L^2(\Gamma_\varrho)$ and
  $L^2 (\Gamma-(\Gamma_\varrho\cup \Gamma_{00}))$, respectively,
   when $v$ runs throughout space $K(D)$. This implies that
    \begin{eqnarray*}  \Delta_g u=0 \;\;\mbox{on}\;\;
   \Gamma-(\Gamma_\varrho\cup \Gamma_{00}), \;\; \mbox{and}\;\;
   \triangle_g u +\frac{1}{\mu^{\star}} \varrho
  \frac{\partial u}{\partial \nu}
  =0 \;\; \mbox{on}\;\; \Gamma_\varrho.\end{eqnarray*}
 Therefore, (\ref{3-15}) holds.  In a similar way, we can prove
 the desired result for $G_{\mathcal{N}}$.
 \ \ $\square$

\vskip 0.24 true cm

  \noindent  {\bf Theorem 3.7.} \  {\it  Let $(\mathcal{M},g)$ be a real
 analytic  Riemannian manifold, and let $D\subset (\mathcal{M},g)$
  be a bounded domain with piecewise smooth boundary $\Gamma$.
 Assume that $ \Gamma_{00}$ is a real analytic $(n-1)$-dimensional
 surface in $\Gamma-\bar \Gamma_\varrho$.
   If $u$ is an eigenfunction of
 the transformations $G_{{\mathcal{K}}^d}^{(\star)}$
  with eigenvalue $\mu^\star\ne 0$,
  then $u$ has derivatives of any order in $D$ and is such that
  \begin{eqnarray}  \label{3--30} \left\{\begin{array} {ll} \triangle_g^2 u=0 \quad \mbox{in}\;\; D,\\
      u=0 \; \; \mbox{on}\;\; \Gamma_\varrho,\\
      u=\frac{\partial u}{\partial \nu}=0 \quad \mbox{on}\;\;
      \Gamma_{00}, \\
 \Delta_g u=0 \;\; \mbox{and}\;\;  \frac{\partial (\triangle_g u)}{\partial \nu}=0 \;\; \mbox{on}\;\;
  \Gamma -(\Gamma_\varrho \cup \Gamma_{00}),\\
   \triangle_g u+ \kappa \, \varrho \,
 \frac{\partial u}{\partial \nu}=0 \;\; \mbox{on}\;\; \Gamma_{\varrho}, \quad \; \mbox{with}\;\;
   \kappa=\frac{1}{\mu^\star}. \end{array}\right.\end{eqnarray}}

\vskip 0.2 true cm

 \noindent  {\bf Proof.} \
 If $G_{{\mathcal{K}}^d}^{(\star)} u= \mu^\star
u$, then we have that $u=0$ on $\Gamma_\varrho$ and
$u=\frac{\partial u}{\partial \nu}=0$ on $\Gamma_{00}$, and that
 \begin{eqnarray} \label{3/25}
 \int_{\Gamma_\varrho} \varrho \frac{\partial u}{\partial \nu}\,
 \frac{\partial v}{\partial \nu}\, ds =
   \mu^\star \int_D (\triangle_g u)(\triangle_g v) dR \quad \; \mbox{for all}\;\;
 v\in {\mathcal{K}}^d (D).\end{eqnarray}
  Applying Green's formula on the right-hand side of (\ref{3/25}), we get
   that \begin{eqnarray} \label{3..27}
    \int_D  (\triangle_g^2 u) v\, dR+
 \int_{\partial D} \frac{\partial (\triangle_g u)}{\partial \nu} \; v\, ds
 -\int_{\Gamma-(\Gamma_\varrho\cup
 \Gamma_{00})} (\triangle_g u)\frac{\partial v}{\partial \nu}\, ds \\
   -
\int_{\Gamma_\varrho} \left(\triangle_g u+\frac{1}{\mu^\star}
\varrho \frac{\partial u}{\partial \nu} \right) \frac{\partial
v}{\partial \nu} \, ds  =0 \quad \; \mbox{for all}\;\; v\in {K}^d
(D).\nonumber
  \end{eqnarray}
  By taking all $v\in C^\infty_0(D)$, we obtain $\triangle_g^2 u=0$ in
  $D$.
  Note that $v\big|_{\Gamma_\varrho}=0$ and $v\big|_{\Gamma_{00}}=
  \frac{\partial v}{\partial \nu}\big|_{\Gamma_{00}}=0$,
  and that $v\big|_{\Gamma -(\Gamma_\varrho\cup \Gamma_{00})}$
  and $\frac{\partial v}{\partial \nu}\big|_{\Gamma -\Gamma_{00}}$
   run throughout
  the spaces  $L^2(\Gamma -(\Gamma_{\varrho}\cup \Gamma_{00}))$
  and $L^2(\Gamma -\Gamma_{00})$, respectively,
  when $v$ runs throughout the space $K^d(D)$.
    Thus we have
\begin{eqnarray*} \label{3..28}
    & \Delta_g u=0 \;\; \mbox{and}\;\;
    \frac{\partial (\triangle_g u)}{\partial \nu}=0 \;\; \mbox{on}\;\;
  \Gamma -(\Gamma_\varrho\cup \Gamma_{00}),\\
  & \quad\quad \qquad \qquad \qquad \qquad \qquad \triangle u+\frac{1}{\mu^\star} \varrho
\frac{\partial u}{\partial \nu} =0 \quad \,\mbox{on}\;\;
\Gamma_\varrho. \qquad \qquad \qquad \qquad \qquad  \quad   \square
\end{eqnarray*}

\vskip 0.25 true cm

  \noindent  {\bf Theorem 3.8.} \ \  {\it  Let $(\mathcal{M},g)$, $D$ and $\Gamma_{00}$ be as
  in Theorem 3.7.
  Assume that $\varsigma_k$ and $\kappa_k$ are the $k$-th Steklov eigenvalues of the following problems:
  \begin{eqnarray}  \label{3-30} \left\{\begin{array} {ll} \triangle_g^2 u=0 \quad \mbox{in}\;\; D,\\
      u=0 \; \; \mbox{on}\;\; \Gamma_\varrho,\\
      u=\frac{\partial u}{\partial \nu}=0 \quad \mbox{on}\;\;
      \Gamma_{00}, \\
   \frac{\partial u}{\partial \nu}=0\;\; \mbox{and}\;\;
    \frac{\partial (\triangle_g u)}{\partial \nu}=0 \;\; \mbox{on}\;\;
  \Gamma -(\Gamma_\varrho \cup \Gamma_{00}),\\
   \triangle_g u+ \varsigma \, \varrho \,
 \frac{\partial u}{\partial \nu}=0 \;\; \mbox{on}\;\; \Gamma_{\varrho} \end{array}\right.\end{eqnarray}
  and
  \begin{eqnarray}  \label{3-31} \left\{\begin{array} {ll} \triangle_g^2 u=0 \quad \mbox{in}\;\; D,\\
      u=0 \; \; \mbox{on}\;\; \Gamma_\varrho,\\
      u=\frac{\partial u}{\partial \nu}=0 \quad \mbox{on}\;\;
      \Gamma_{00}, \\
  \Delta_g u=0\;\; \mbox{and}\;\; \frac{\partial (\triangle_g u)}{\partial \nu}=0
   \;\;  \mbox{on}\;\;
  \Gamma -(\Gamma_\varrho \cup \Gamma_{00}),\\
   \triangle_g u+ \kappa \, \varrho \,
 \frac{\partial u}{\partial \nu}=0 \;\; \mbox{on}\;\; \Gamma_{\varrho}, \end{array}\right.\end{eqnarray}
   respectively.
  Then $\varsigma_k\le \kappa_k$ for all $k\ge 1$.}

\vskip 0.26 true cm

 \noindent  {\bf Proof.} \  \ For $0<\alpha<1$, let $u_k=u_k(\alpha, x)$  be
 the normalized eigenfunction corresponding to the $k$-th Steklov eigenvalue $\lambda_k$ for the
 following problem:
   \begin{eqnarray} \label{3---31}  \left\{\begin{array} {ll} \triangle_g^2 u_k=0 \quad \mbox{in}\;\; D,\\
      u_k=0 \; \; \mbox{on}\;\; \Gamma_\varrho,\\
      u_k=\frac{\partial u_k}{\partial \nu}=0 \quad \mbox{on}\;\;
      \Gamma_{00}, \\
  \alpha \Delta_g u_k+(1-\alpha) \frac{\partial u_k}{\partial \nu}=0\;\; \mbox{and}\;\;
  \frac{\partial (\Delta_g u_k)}{\partial \nu}=0
   \;\;  \mbox{on}\;\;
  \Gamma -(\Gamma_\varrho \cup \Gamma_{00}),\\
   \triangle_g u_k+ \lambda \, \varrho \,
 \frac{\partial u_k}{\partial \nu}=0 \;\; \mbox{on}\;\; \Gamma_{\varrho}. \end{array}\right.\end{eqnarray}
   It is easy to verify (cf, p.$\;$410 or Theorem 9 of p.$\;$419 in \cite{CH})
   that the $k$-th Steklov eigenvalue
    $\lambda_k=\lambda_k(\alpha)$ is continuous on the closed interval $[0,1]$
   and differentiable in the open interval $(0,1)$, and that
   $u_k(\alpha, x)$ is also differentiable with respect to $\alpha$ in
   $(0,1)$ (see, \cite{Fri}).
   \  We will denote by $'$ the derivative with respect to $\alpha$.
    Then
 \begin{eqnarray} \label{32}  \quad \;\; \quad \left\{\begin{array} {ll}
  \triangle_g^2 u'_k=0 \quad \mbox{in}\;\; D,\\
      u'_k=0 \; \; \mbox{on}\;\; \Gamma_\varrho,\\
      u'_k=\frac{\partial u'_k}{\partial \nu}=0 \quad \mbox{on}\;\;
      \Gamma_{00}, \\
   \Delta_g u_k+ \alpha \triangle_g u'_k -\frac{\partial u_k}{\partial
   \nu} +  (1-\alpha) \frac{\partial u'_k}{\partial \nu}=0\,\, \mbox{and}\,\,
  \frac{\partial (\Delta_g u'_k)}{\partial \nu}=0
   \,\,  \mbox{on}\,\,
  \Gamma -(\Gamma_\varrho \cup \Gamma_{00})\\
   \triangle_g u'_k+ \lambda' \, \varrho \,
 \frac{\partial u_k}{\partial \nu}+\lambda \varrho \frac{\partial u'_k}{\partial \nu}=0
 \;\; \mbox{on}\;\; \Gamma_{\varrho}. \end{array}\right. \end{eqnarray}
  Multiplying  (\ref{32}) by $u_k$,  integrating the product over
  $D$, and then applying Green's formula,
 we get
\begin{eqnarray*} 0&=&\int_D (\Delta^2_g u'_k)u_k \,dR =
\int_D (\Delta_g^2 u_k)u'_k\, dR  - \int_{\partial D} (\Delta_g
  u_k)\frac{\partial u'_k}{\partial \nu}\, ds \\
  && +\int_{\partial D} u'_k
  \frac{\partial (\Delta_g u_k)}{\partial \nu}\, ds -\int_{\partial D}
  u_k\frac{\partial (\Delta_g u'_k)}{\partial \nu}\, ds +\int_{\partial
  D} (\Delta_g u'_k) \frac{\partial u_k}{\partial \nu}\,ds \nonumber \\
  &=& -\left [ \int_{\Gamma_\varrho} (\Delta_g u_k)\frac{\partial
  u'_k}{\partial \nu} \, ds + \int_{\Gamma -(\Gamma_\varrho \cup
  \Gamma_{00})} (\Delta_g u_k)\frac{\partial
  u'_k}{\partial \nu} \, ds \right]\nonumber \\
  && +  \left [ \int_{\Gamma_\varrho} (\Delta_g u'_k)\frac{\partial
  u_k}{\partial \nu} \, ds + \int_{\Gamma -(\Gamma_\varrho \cup
  \Gamma_{00})} (\Delta_g u'_k)\frac{\partial
  u_k}{\partial \nu} \, ds \right]\nonumber \\&=&
  \left [ \int_{\Gamma_\varrho}
  \left (\lambda \varrho \frac{\partial u_k}{\partial \nu} \right)\frac{\partial
  u'_k}{\partial \nu} \, ds
   + \int_{\Gamma -(\Gamma_\varrho \cup
  \Gamma_{00})} \left(\frac{1-\alpha}{\alpha} \,\frac{\partial u_k}{\partial \nu}\right)\frac{\partial
  u'_k}{\partial \nu} \, ds \right]\nonumber \\
  && +   \int_{\Gamma_\varrho} \left(-\lambda' \varrho \frac{\partial
  u_k}{\partial \nu} -\lambda \varrho \frac{\partial u'_k}{\partial \nu}\,
\right) \frac{\partial
  u_k}{\partial \nu} \, ds \nonumber \\
  && + \int_{\Gamma -(\Gamma_\varrho \cup
  \Gamma_{00})} \left( -\frac{1}{\alpha}\, \Delta_g u_k
  +\frac{1}{\alpha}\,
  \frac{\partial
  u_k}{\partial \nu} -\frac{1-\alpha}{\alpha} \, \frac{\partial u'_k}{\partial \nu} \right)
  \frac{\partial u_k}{\partial \nu} \, ds \nonumber \\
&=& -\lambda' \int_{\Gamma_\varrho} \varrho \left(\frac{\partial
u_k}{\partial \nu}\right)^2 ds +
 \int_{\Gamma -(\Gamma_\varrho \cup
  \Gamma_{00})} \left[\left(\frac{1-\alpha}{\alpha^2}\right) \frac{\partial u_k}{\partial \nu}
  +\frac{1}{\alpha} \,\frac{\partial u_k}{\partial \nu}\right] \frac{\partial u_k}{\partial \nu}\,
  ds \nonumber \\
  &=& -\lambda' \int_{\Gamma_\varrho} \varrho \left(\frac{\partial
u_k}{\partial \nu}\right)^2 ds +
 \int_{\Gamma -(\Gamma_\varrho \cup
  \Gamma_{00})} \bigg(\frac{1}{\alpha} \, \frac{\partial u_k}{\partial
  \nu} \bigg)^2 ds,\nonumber \end{eqnarray*}
 i.e., $$\lambda'_k(\alpha)= \frac{\int_{\Gamma -(\Gamma_\varrho
\cup
  \Gamma_{00})} \left(\frac{1}{\alpha} \, \frac{\partial u_k}{\partial
  \nu} \right)^2 ds}{\int_{\Gamma_\varrho} \varrho \left(\frac{\partial
u_k}{\partial \nu}\right)^2 ds} >0 \quad \; \mbox{for all}\;\;
0<\alpha<1.$$
 This implies that $\lambda_k$ is increasing with respect to $\alpha$ in
 $(0,1)$.  Note that
 if we change the $\alpha$ from $0$ to $1$, each individual Steklov
 eigenvalue $\lambda_k$ increase monotonically form the value
 $\varsigma_k$ which is the $k$-th Steklov eigenvalue of (\ref{3-30}) to
the value $\kappa_k$ which is the $k$-th Steklov eigenvalue
 (\ref{3-31}).
 Thus, we have that $ \varsigma_k \le \kappa_k$  for all $k$.   \ \ $\square$

\vskip 0.26 true cm

 Conversely, the following proposition shows that a sufficiently smooth function satisfying
(\ref{3-15}) (respectively, (\ref{3--30})) is an eigenfunction of
$G_{\mathcal{K}}^{(\star)}$ or $G_{\mathcal{N}}^{(\star)}$
(respectively, $G_{{\mathcal{K}}^d}^{(\star)}$).

 \vskip 0.20 true cm

 \noindent  {\bf Proposition 3.9.} \ \  {\it  Let $\bar D$ be bounded domain with  piecewise smooth
  boundary. Assume that $u$ belongs to $C^4(\bar D)$.

a) \ \  If $\Gamma_{\varrho}\ne \Gamma$ and $u$ satisfies
(\ref{3-15}), then $u\in \mathcal{K}$ and $u$ is an  eigenfunction
of $G_{\mathcal{K}}^{(\star)}$ with the eigenvalue
$\mu^\star=\gamma^{-1}$,
\begin{eqnarray}  \label{3-27} G_{\mathcal{K}}^{(\star)}u=\gamma^{-1} u.\end{eqnarray}

b) \ \  If $\Gamma_{\varrho}\ne \Gamma$ and $u$ satisfies
(\ref{3--30}), then $u\in {\mathcal{K}}^d$ and $u$ is an
eigenfunction of $G_{{\mathcal{K}}^d}^{(\star)}$ with the eigenvalue
$\mu^\star=\kappa^{-1}$,
\begin{eqnarray}  \label{3-28} G_{{\mathcal{K}}^d}^{(\star)}u=\kappa^{-1} u.\end{eqnarray}

c) \ \  If $\Gamma_{\varrho}= \Gamma$ and $u$ satisfies
(\ref{3-15}), then $u\in \mathcal{N}$ and $u$ is an  eigenfunction
of $G_{\mathcal{N}}^{(\star)}$ with the eigenvalue
$\mu^\star=\gamma^{-1}$,
\begin{eqnarray} \label {3-29} G_{\mathcal{N}}^{(\star)}u=\gamma^{-1} u.\end{eqnarray}}

 \noindent  {\bf Proof.} \ \  i) \ \  $\Gamma_\varrho \ne \Gamma$. We
 claim that there is no eigenvalue $\gamma=0$. Suppose by contradiction that there is a  function
 $u$ in $C^4(\bar D)$ satisfying
 \begin{eqnarray} \label {3-36}\left\{\begin{array} {ll} \triangle_g^2 u=0 \;\;\mbox{in}\;\; D, \; \quad \;\;
 u=0\;\;\mbox{on}\;\; \Gamma,\\
  \frac{\partial u}{\partial \nu}=0 \;\;\mbox{on}\;\;
  \Gamma_{00},  \quad \;\; \mbox{and}\;\; \triangle_g
  u=0 \;\; \mbox{on} \;\; \Gamma-\Gamma_{00}.\end{array} \right.\end{eqnarray}
  By multiplying the above equation by $u$, integrating the result over
  $D$, and using Green's formula, we derive
  \begin{eqnarray*}
  && 0=\int_D u(\triangle_g^2 u) dR= \int_D |\triangle_g u|^2 dR -  \int_\Gamma u
  \frac{\partial (\triangle_g u)}{\partial \nu} \,ds \\
  && \;\; \quad \; \; + \int_{\Gamma} (\triangle_g u)
 \frac{\partial u}{\partial \nu} ds=
 \int_D |\triangle_g u|^2 dR.  \end{eqnarray*}
This implies that $\triangle_g u=0$ in $D$. Since $u=0$ on $\Gamma$,
by the maximum principle we get that $u=0$ in $D$. The claim is
proved.

  In view of assumptions, we see that $u\in \mathcal{K}$.
  By (\ref{3-15}) and Green's formula, it follows that
  for an arbitrary $v\in K(D)$
  \begin{eqnarray*}  \langle G_{\mathcal{K}}^{(\star)}
   u, v\rangle^\star &=&[u,v] = \int_{\Gamma_\varrho}
  \varrho \, \frac{\partial u}{\partial \nu}\, \frac{\partial v}{\partial \nu}\, ds\\
    &=&- \gamma^{-1}\int_{\Gamma_\varrho} (\triangle_g u)\frac{\partial
  v} {\partial \nu} \, ds
  =  - \gamma^{-1}
  \int_{\Gamma} (\triangle_g u)\frac{\partial v}{\partial \nu}\, ds \\
  &=&  -\gamma^{-1}  \left[\int_{\Gamma}  \frac{\partial
  (\triangle_g u)} {\partial \nu}\, v\, ds
   -\int_D (\triangle_g u) (\triangle_g v) dR + \int_D
  v(\triangle_g^2 u)dR\right]\\
  &=& \gamma^{-1} \int_D (\triangle_g
  u)(\triangle_g v) dR = \gamma^{-1}\langle u, v \rangle^\star. \end{eqnarray*}
  Therefore,  \begin{eqnarray*} \langle G_{\mathcal{K}}^{(\star)}u
     - \gamma^{-1}u, v \rangle^\star=0 \quad\quad \mbox{for all } v\in K(D), \end{eqnarray*}
 which  implies (\ref{3-27}).
 By a similar way, we can prove b).

   ii) \ \  $\Gamma_\varrho = \Gamma$.  We claim that
    there is no eigenvalue $\gamma=0$. If it is not this case,  then there is a function $u$
 in $C^4(\bar D)$ satisfying
 \begin{eqnarray*} \left\{ \begin{array}{ll}  \triangle_g^2 u=0 \;\;\mbox{in}\;\; D,\\
  u=0\;\;\mbox{on}\;\; \Gamma,\\
   \triangle_g
  u=0\;\;\mbox{on}\;\; \Gamma.\end{array}\right.\end{eqnarray*}
  Setting $v:=\triangle_g u$ in $D$, we get
\begin{eqnarray*} \left\{\begin{array} {ll}\triangle_g v=0 \;\;\mbox{in}\;\; D, \\
 v=0\;\;\mbox{on}\;\; \Gamma.\end{array} \right.\end{eqnarray*}
   By the maximum principle it follows  that $v=0$ in $D$.
   Thus, we have
   \begin{eqnarray*} \left\{\begin{array} {ll} \triangle_g u=0 \;\;\mbox{in}\;\; D, \\
 u=0\;\;\mbox{on}\;\; \Gamma,
\end{array} \right.\end{eqnarray*}
 so that $u=0$ in $D$.

Now, if $u$ is a solution of (\ref{3-15}) with eigenvalue $\gamma>
0$, proceeding as in a), we can prove that $u\in \mathcal{N}$ and
(\ref{3-29}) holds. \ \ \ $\square$

\vskip 0.23 true cm

\noindent {\bf Remark 3.10.} \  Each of transformations
$G_{\mathcal{K}}^\star$, $G_{{\mathcal{K}}^d}^\star$ and
$G_{\mathcal{N}}^\star$ corresponds to a biharmonic Steklov problem
given by
  the quadratic forms
  \begin{eqnarray*}  \langle u, u\rangle^\star =\int_D
  |\triangle_g u|^2 dR \end{eqnarray*}
  and \begin{eqnarray*} [u,u] =\int_{\Gamma_\varrho}  \varrho
  \left(\frac{\partial u}{\partial \nu}\right)^2
  ds \end{eqnarray*}
  and the function classes of ${\mathcal{K}}^\star$, ${{\mathcal{K}}^d}^\star$
  and ${\mathcal{N}}^\star$, respectively.
  The eigenvalues $\gamma_k$ and $\kappa_k$ of these  biharmonic Steklov problems are given by
  \begin{eqnarray} \label {3-37} \gamma_k \;\; \mbox{and}\;\; \kappa_k=1/\mu^\star_k
  \;\quad k=1,2,3, \cdots.\end{eqnarray}
Since $0$ is the only limit point of $\mu_k^\star$, the only
possible limit points of $\gamma_k$ and $\kappa_k$ are $+\infty$.

 \vskip 1.39 true cm

\section{Biharmonic Steklov eigenvalues on an $n$-dimensional rectangular parallelepiped}

\vskip 0.45 true cm

  Let  $D=\{x\in {\Bbb R}^n \big|
 0\le x_i\le l_i, \, i=1, \cdots, n\}$ with boundary $\Gamma$, and let
 $\Gamma_\varrho^+ =\{x\in {\Bbb R}^n\big| 0\le x_i \le l_i \,\,
\mbox{when}\,\, i<n, \, x_n =0\}$. Let
 $\Gamma^{l_n} =\{x\in {\Bbb R}^n\big| 0\le x_i \le l_i \,\,
\mbox{when}\,\, i<n, \, x_n =l_n\}$.
 Our first purpose, in this section, is to discuss
 the biharmonic Steklov eigenvalue problem on
 $n$-dimensional rectangular parallelepiped $D$:
 \begin{eqnarray} \label {4-1}  \left\{ \begin{array}{ll} \triangle^2 u = 0\quad \;
   \mbox{in}\;\; D,\\
   u=0 \; \; \mbox{on}\;\; \Gamma,\\
   \frac{\partial u}{\partial \nu} =0 \;\; \mbox{on}\;\; \Gamma^{l_n},
\;\; \Delta u=0\;\; \mbox{on}\;\; \Gamma -(\Gamma_{\varrho}^+ \cup \Gamma^{l_n}),\\
   \triangle u+\gamma \varrho \frac{\partial u}{\partial \nu}=0 \;\;
    \mbox{on}\;\; \Gamma_{\varrho}^+, \quad \,  \varrho =constant>0\;\;\mbox{on}\;\;
     \Gamma_\varrho^+.
   \end{array}\right.\end{eqnarray}

   We first consider the special solution of (\ref{4-1}) which has the following form:
\begin{eqnarray*} u=X(x_1,\cdots, x_{n-1})\, Y(x_n).\end{eqnarray*}
 Since \begin{eqnarray*} \Delta u&=&
\big(\Delta_{n-1} X(x_1, \cdots, x_{n-1})\big)Y(x_n)+2\nabla X (x_1,
\cdots, x_{n-1}) \cdot \nabla Y(x_n) \\ && +\big(X(x_1, \cdots,
 x_{n-1})\big)Y''(x_n)=
  \big(\Delta_{n-1} X(x_1, \cdots, x_{n-1})\big)Y(x_n)\\&&
  +\big(X(x_1, \cdots,
 x_{n-1})\big)Y''(x_n) \end{eqnarray*}
 and
    \begin{eqnarray*}  \Delta^2 u &=&(\Delta^2_{n-1} X(x_1, \cdots, x_{n-1}))Y(x_n)
    + 2 \big(\Delta_{n-1} X(x_1, \cdots, x_{n-1})\big) Y'' (x_n)\\
    &&+ \big(X(x_1,\cdots, x_{n-1})\big)Y''''(x_n),\end{eqnarray*}
where $$\Delta_{n-1} X(x_1, \cdots, x_{n-1}) =\sum_{i=1}^{n-1}
\frac{\partial^2
 X}{\partial x_i^2},$$
 we find by  $\Delta^2 u=0$ that
  \begin{eqnarray*} & \big(\Delta^2_{n-1} X(x_1, \cdots, x_{n-1})\big) Y(x_n)
 +2 (\Delta_{n-1} X(x_1, \cdots, x_{n-1})\big) Y''(x_n)\\
 &+\big(X(x_1, \cdots, x_{n-1})\big)
 Y''''(x_n)=0,\end{eqnarray*}
  so that  \begin{eqnarray} \label {4..2}\quad\quad \quad \; \frac{\Delta^2_{n-1} X(x_1, \cdots,
 x_{n-1})}{X(x_1, \cdots, x_{n-1})} +2\,\frac{\Delta_{n-1} X(x_1, \cdots,
 x_{n-1})}{X(x_1, \cdots, x_{n-1})}\, \frac{Y''(x_n)}{Y(x_n)} +
 \frac{Y''''(x_n)}{Y(x_n)}=0.\end{eqnarray}
Differentiating (\ref{4..2}) with respect to $x_n$, we obtain that
\begin{eqnarray*} \label{4..3}  2 \,\frac{\Delta_{n-1} X(x_1, \cdots, x_{n-1})}{ X(x_1, \cdots,
x_{n-1})}\left[\frac{Y''(x_n)}{Y(x_n)}\right]'
+\left[\frac{Y''''(x_n)}{Y(x_n)}\right]'=0.\end{eqnarray*} The above
equation holds if and only if
\begin{eqnarray} \label{4..4}   \frac{\Delta_{n-1} X(x_1, \cdots, x_{n-1})}{ X(x_1, \cdots,
x_{n-1})} =-\frac{\left[\frac{Y''''(x_n)}{Y(x_n)}\right]'}{ 2
\left[\frac{Y''(x_n)}{Y(x_n)}\right]'}= -\eta^2,\end{eqnarray} where
$\eta^2$ is a constant.  Therefore, we have that
\begin{eqnarray}\label{4..5} \Delta_{n-1} X(x_1, \cdots, x_{n-1}) +\eta^2 X(x_1,
\cdots, x_{n-1}) =0 \end{eqnarray} and
$$ \left[\frac{Y''''(x_n)}{Y(x_n)}\right]' -2 \eta^2
\left[\frac{Y''(x_n)}{Y(x_n)}\right]'=0.$$
 From (\ref{4..5}), we get \begin{eqnarray} \qquad \quad \quad \;\Delta^2_{n-1}
 X(x_1,\cdots, x_{n-1})
=-\eta^2 \Delta_{n-1} X (x_1,\cdots, x_{n-1})=\eta^4 X(x_1,\cdots,
x_{n-1}).
\end{eqnarray}
  Substituting this in (\ref{4..2}), we obtain the following equation
\begin{eqnarray} \label{4..6} Y''''(x_n)- 2\eta^2 Y''(x_n) +\eta^4
Y(x_n)=0.\end{eqnarray}
   It is easy to verify  that the general solutions of (\ref{4..6}) have
   the form:
  \begin{eqnarray} \label{4.6'} Y(x_n) = A \cosh \, \eta x_n +B\sinh\, \eta
  x_n  +  C x_n \cosh\, \eta x_n +D x_n \sinh\, \eta x_n.\end{eqnarray}
 By setting   $Y(0)=Y(l_n)=0, \;\; Y'(0)=1, \; \; Y'(l_n)=0$, we get
 \begin{eqnarray}  \label{4--7} \quad \;  Y(x_n)&=& \left(\frac{-\eta l_n^2}{\sinh^2 \eta l_n
  - \eta^2 l_n^2}\right) \sinh \eta x_n +\left(\frac{\sinh^2 \eta l_n}{
 \sinh^2 \eta l_n -\eta^2 l_n^2} \right) x_n \cosh \eta x_n
 \\ && + \left(\frac{\eta l_n -(\sinh \eta l_n) \cosh \eta l_n}{
 \sinh^2 \eta l_n -\eta^2 l_n^2}\right) x_n\sinh \eta x_n.\nonumber
   \end{eqnarray}
     It is well-known that for the Dirichlet eigenvalue problem
     \begin{eqnarray}  \left\{ \begin{array} {ll} \Delta_{n-1} X (x_1,\cdots, x_{n-1})+ \eta^2
      X(x_1, \cdots, x_{n-1})=0
     \quad \; \mbox{in}\;\; \Omega,\\
     u=0 \quad \; \mbox{on}\; \; \partial \{(x_1,\cdots, x_{n-1}) \big|
  0\le x_i\le l_i, \,\, i=1, \cdots, n-1\},
     \end{array} \right. \end{eqnarray}
  there exist the eigenfunctions
  \begin{eqnarray} \label{4..10} X(x_1, \cdots, x_{n-1})= c \left(\sin \frac{m_1\pi}{l_1}x_1\right)\cdots
  \left(\sin
  \frac{m_{n-1} \pi } {l_{n-1}}x_{n-1}\right), \end{eqnarray}
  which correspond to the eigenvalues
   \begin{eqnarray*}  \eta^2=\sum_{i=1}^{n-1} \left(\frac{m_i \pi}{l_i}\right)^2,\quad \;\;
 \mbox{where}\;\; m_i=1, 2,3, \cdots. \end{eqnarray*}
    Therefore,
    \begin{eqnarray} \label{4.11}  \quad \;\,\; u&=& \big(X(x_1, \cdots, x_{n-1})\big)Y(x_n)
    \\
    &= & c \left(\sin \frac{m_1\pi}{l_1}x_1\right)\cdots
  \left(\sin \frac{m_{n-1} \pi } {l_{n-1}}x_{n-1}\right)
   \left[
   \left(\frac{-\eta l_n^2}{\sinh^2 \eta l_n
  - \eta^2 l_n^2}\right) \sinh \eta x_n \right.\nonumber \\&& \left. +\left(\frac{\sinh^2 \eta l_n}{
 \sinh^2 \eta l_n -\eta^2 l_n^2} \right) x_n \cosh \eta x_n \right.\nonumber\\
 &&  + \left.\left(\frac{\eta l_n -(\sinh \eta l_n) \cosh \eta l_n}{
 \sinh^2 \eta l_n -\eta^2 l_n^2}\right) x_n\sinh \eta x_n\right].\nonumber
  \end{eqnarray}
 Since $$Y''(0)= 2 \eta \left(
\frac{\eta l_n- (\sinh \eta l_n)\cosh \eta l_n}{\sinh^2 \eta l_n
-\eta^2 l_n^2}\right)\quad \;\mbox{and}\; \; Y'(0)=1,$$
 we obtain
 \begin{eqnarray*}(\triangle u)\big|_{x_n=0}&=&
  \big(\Delta_{n-1} X(x_1,\cdots, x_{n-1})\big)Y(0) +\big(X
  (x_1, \cdots, x_{n-1})\big)Y''(0)\\
    &=& 2\eta
   \left(\frac{\eta l_n -(\sinh \eta l_n) \cosh \eta l_n}{
 \sinh^2 \eta l_n -\eta^2 l_n^2}\right)X(x_1,\cdots, x_{n-1}),\end{eqnarray*}
 \begin{eqnarray*} \mbox{and}\qquad  \qquad \qquad \qquad\qquad \;  \frac{\partial u}{\partial
 \nu}\bigg|_{\Gamma_\rho^+}=X( x_1, \cdots, x_{n-1}), \qquad \qquad \qquad \qquad  \qquad\qquad \quad\end{eqnarray*}
so that $$\triangle u+\gamma \varrho \frac{\partial u}{\partial \nu}
=0\quad \; \mbox{on}\;\; \Gamma_\varrho^+$$ with
  $$\gamma = \frac{2\eta l_n}{\varrho l_n} \left(
\frac{(\sinh \eta l_n)\cosh \eta l_n -\eta l_n}{\sinh^2 \eta l_n
-\eta^2 l_n^2}\right).$$

\vskip 0.29 true cm

  Our second purpose is to discuss the biharmonic
 Steklov eigenvalue problem on the $n$-dimensional rectangular parallelepiped $D$:
\begin{eqnarray} \label {4-13}  \left\{ \begin{array}{ll} \triangle^2 u = 0\quad \;
   \mbox{in}\;\; D,\\
   u=0 \; \; \mbox{on}\;\; \Gamma_\varrho^+, \quad \,
   u=\frac{\partial u}{\partial \nu}=0 \;\;\mbox{on}\;\;
   \Gamma^{l_n},\\
      \frac{\partial u}{\partial \nu}= \frac{\partial (\Delta u)}{\partial \nu}= 0
      \;\; \mbox{on}\;\; \Gamma-(\Gamma_\varrho^+\cup \Gamma^{l_n}),\\
    \triangle u+\varsigma \varrho \frac{\partial u}{\partial \nu}=0 \;\;
    \mbox{on}\;\; \Gamma_\varrho^+, \quad \,  \varrho =constant>0\;\;\mbox{on}\;\;
     \Gamma_\varrho^+.
   \end{array}\right.\end{eqnarray}

\vskip 0.22 true cm

  Similarly, (\ref{4-13}) has the special solution
  $u=\big(X(x_1, \cdots, x_{n-1})\big)Z(x_n)$ with $Z(x_n)$ having form
 (\ref{4.6'}). According to the boundary conditions of (\ref{4-13}),
 we get that the problem (\ref{4-13}) has the solutions
 \begin{eqnarray*}  \label{4-19} u(x)&=&  u(x_1, \cdots, x_n)\\
   &=& c\left(\cos \frac{m_1 \pi}{l_1}\, x_1\right)\cdots
   \left(\cos \frac{m_{n-1} \pi}{l_{n-1}}\, x_{n-1}\right)
  Z(x_n),\nonumber \end{eqnarray*}
where $m_1, \cdots, m_{n-1}$ are whole numbers, and $Z (x_n)$ is
given by  \begin{eqnarray} \label {4-9} \quad \quad \quad  Z
(x_n)&=&
 \left(\frac{-\beta l_n^2}{\sinh^2 \beta l_n - \beta^2 l_n^2}\right) \sinh \beta x_n
+ \left(\frac{\sinh^2 \beta l_n}{
 \sinh^2 \beta l_n -\beta^2 l_n^2} \right) x_n \cosh \beta x_n\\
 &&  +  \left(\frac{\beta l_n -(\sinh \beta l_n) \cosh \beta l_n}{ \sinh^2 \beta l_n -\beta^2 l_n^2} \right)
 x_n\sinh \beta x_n,
 \nonumber\end{eqnarray}
 $\beta=\big[\sum_{i=1}^{n-1}
(m_i\pi/l_i)^2\big]^{1/2}$ with $\sum_{i=1}^{n-1} m_i\ne 0$. Since
  $\frac{\partial
 u}{\partial \nu}\big|_{\Gamma_\rho^+}=X(x_1, \cdots, x_{n-1})$,
  $\; (\Delta u)\big|_{x_n=0}=\big(X(x_1,\cdots, x_{n-1})\big)
 Z''(0)$ and
  $Z''(0)=2\eta \left(\frac{\beta l_n -(\sinh \beta l_n) \cosh \beta l_n}{
 \sinh^2 \beta l_n -\beta^2 l_n^2}\right)$,
  we get
$\Delta u+\varsigma \varrho \frac{\partial u}{\partial \nu}=0$ on
$\Gamma_\varrho^+$, where $$\varsigma = \frac{2\beta l_n}{\varrho
l_n} \left( \frac{(\sinh \beta l_n)\cosh \beta l_n -\beta
l_n}{\sinh^2 \beta l_n -\beta^2 l_n^2}\right).$$

\vskip 1.39 true cm

\section{Asymptotic distribution of eigenvalues on special domains}

\vskip 0.45 true cm

\noindent {\bf 5.1. Counting function $A(\tau)$.}

\vskip 0.2 true cm

In order to obtain our asymptotic formula, it is an effective way
 to investigate the
distribution of the eigenvalues of the transformation
$G_{\mathcal{K}}$ (respectively, $G_{{\mathcal{K}}^d}$,
$\,G_{\mathcal{N}}$) instead of the transformations
$G_{\mathcal{K}}^{(\star)}$ (respectively,
  $ G_{{\mathcal {K}}^d}^{(\star)}$, $\,G_{\mathcal{N}}^{(\star)}$).
 It follows from (\ref{3.10})---(\ref{3.12}) and (\ref{3-37}) we obtain
 \begin{eqnarray} \mu_k= (1+\lambda_k)^{-1}, \quad \; k=1, 2, 3,
 \cdots, \end{eqnarray}
  where $\mu_k$ denote the $k$-th eigenvalue of
 $G_{\mathcal{K}}$ or
 $G_{{\mathcal{K}}^d}$ or
 $G_{\mathcal{N}}$, and $\frac{1}{\lambda_k}$  is the
 $k$-th eigenvalue of $G_{\mathcal{K}}^{(\star)}$ or
 $G_{{\mathcal{K}}^d}^{(\star)}$
  or $G_{\mathcal{N}}^{(\star)}$) (More precisely, $\lambda_k=\gamma_k$ for $G_{\mathcal{K}}^{(\star)}$
  and $G_{\mathcal{N}}^{(\star)}$, and $\lambda_k=\kappa_k$ for
 $G_{{\mathcal{K}}^d}^{(\star)}$).)
 Since $A(\tau)= \sum_{\lambda_k\le \tau} 1$,  we have
 \begin{eqnarray} A(\tau)= \sum_{\mu_k \ge (1+\tau)^{-1}} 1.
 \end{eqnarray}

\vskip 0.48 true cm

\noindent {\bf 5.2. $\,D$ is an $n$-dimensional  rectangular
parallelepiped and
 $g_{ij}=\delta_{ij}$}.

\vskip 0.3 true cm

 Let $D$ be an $n$-dimensional rectangular parallelepiped,
 $g_{ij}=\delta_{ij}$ in the whole of $\bar D$,
  $\varrho=constant>0$ on one face $\Gamma_{\varrho}^+$ of the
  rectangular parallelepiped, i.e., $D=\{x\in {\Bbb R}^n\big|
 0\le x_i\le l_i, \, i=1, \cdots, n\}$,
 $\Gamma_\varrho^+ =\{x\in {\Bbb R}^n\big| 0\le x_i \le l_i \,\,
\mbox{when}\,\, i<n, \, x_n =0\}$) and  $\Gamma_{00}=\Gamma^{l_n}
=\{x\in {\Bbb R}^n\big| 0\le x_i \le l_i \,\, \mbox{when}\,\, i<n,
\, x_n =l_n\}$.
 Without loss of generality, we assume $l_i < l_n$ for all $i<n$.

 For the above domain $D$, except for the
$K(D)$ and $K^d (D)$ in Section 3, we introduce the linear space of
functions
     \begin{eqnarray*}  K^0(D)= \{u\big| u\in H_0^1(D) \cap
     H^2 (D)\cap C^\infty(\bar D), \;  \frac{\partial u}{\partial \nu}=0\;\;
  \mbox{on}\;\; \Gamma_{00}, \;\; \Delta u=0 \;\;
  \mbox{on}\;\; \Gamma-(\Gamma_\varrho^+\cup \Gamma_{00})\}.\end{eqnarray*}
 Clearly,
   \begin{eqnarray}  \label {5-3} K^0(D) \subset K(D) \subset K^d
  (D).\end{eqnarray}
     Closing $K^0$, $K$ and $K^d$ respect to the norm $\|u\|=
   \sqrt{\langle u, u\rangle}$, we obtain the Hilbert spaces
   ${\mathcal{K}}^0$, $\mathcal{K}$ and ${\mathcal{K}}^d$, and
  \begin{eqnarray} {\mathcal{K}}^0 \subset \mathcal{K} \subset
  {\mathcal{K}}^d.\end{eqnarray}
   According to Theorem 3.3, we see that the bilinear functional
   \begin{eqnarray} \label{5-4} [u,v] =\int_{\Gamma_\varrho^+} \varrho
   \frac{\partial u}{\partial \nu} \frac{\partial v}{\partial \nu}
   \, ds\end{eqnarray}
   defines self-adjoint, completely continuous transformations $G^0$,
   $G$ and $G^d$ on ${\mathcal{K}}^0$, $\mathcal{K}$ and
   ${\mathcal{K}}^d$, respectively (cf. Section 3).
    Obviously,   \begin{eqnarray*} \label{5-5}  \langle G^0 u, v\rangle = \langle Gu,
      v\rangle
   \quad \; \mbox{for all}\;\; u, v  \;\; \mbox{in}\;\;
   {\mathcal{K}}^0,\end{eqnarray*}
     \begin{eqnarray*}  \label{5-6} \langle G u, v\rangle = \langle G^d u,
     v\rangle
   \quad \; \mbox{for all}\;\; u, v  \;\; \mbox{in}\;\;
   {\mathcal{K}},\end{eqnarray*}
  from which we deduce immediately by Proposition 2.3 that
  \begin{eqnarray} \mu_k^0 \le \mu_k \le \mu_k^d, \quad k=1,2,3,
  \cdots, \end{eqnarray}
  where $\{\mu_k^0\}$ and $\{\mu_k^d\}$ are the eigenvalues of $G^0$
  and $G^d$, respectively. Hence
  \begin{eqnarray} \label{5-66} A^0(\tau) \le A(\tau) \le A^d(\tau) \quad \;\;
  \mbox{for all}\;\; \tau,\end{eqnarray}
  where \begin{eqnarray} \label {5-7} A^0(\tau)=\sum_{\mu_k^0\ge (1+\tau)^{-1}}
1\end{eqnarray}
 and \begin{eqnarray} \label {5-8} A^d(\tau)=\sum_{\mu_k^d\ge (1+\tau)^{-1}}
1.\end{eqnarray}
 We shall estimate the asymptotic behavior
  of $A^0(\tau)$ and $A^d(\tau)$.
 It is easy to verify (cf. Theorems 3.6, 3.7) that the eigenfunctions of the transformations $G^0$ and $G^d$,
  respectively, satisfy
 \begin{eqnarray} \label {5-13} \quad \;\; \left\{ \begin{array}{ll} \triangle^2 u = 0\quad \;
   \mbox{in}\;\; D,\\
   u=0 \; \; \mbox{on}\;\; \Gamma,\quad \\
   \frac{\partial u}{\partial \nu} =0 \;\; \mbox{on}\;\; \Gamma^{l_n},
\,\, \mbox{and}\;\;
\Delta u=0\;\;\mbox{on}\;\; \Gamma-(\Gamma_\varrho \cup \Gamma^{l_n}),\\
   \triangle u+\gamma \varrho \frac{\partial u}{\partial \nu}=0 \;\;
    \mbox{on}\;\; \Gamma_{\varrho}^+, \quad \,  \varrho =constant>0\;\;\mbox{on}\;\;
     \Gamma_\varrho^+.
   \end{array}\right.\end{eqnarray}
and
\begin{eqnarray} \label {5-14} \quad  \left\{ \begin{array}{ll} \triangle^2 u = 0\quad \;
   \mbox{in}\;\; D,\\
   u=0 \; \; \mbox{on}\;\; \Gamma_\varrho^+,\quad \;\; u=
   \frac{\partial u}{\partial \nu}=0 \;\; \mbox{on}\;\;
   \Gamma^{l_n},\\
    \frac{\partial (\Delta u)}{\partial \nu}=
    0 \;\; \mbox{and}\;\; \Delta u=0\;\; \mbox{on}\;\;
    \Gamma-(\Gamma_\varrho^+\cup \Gamma^{l_n}),\\
    \triangle u+\kappa \varrho \frac{\partial u}{\partial \nu}=0 \;\;
    \mbox{on}\;\; \Gamma_\varrho^+, \quad \,  \varrho =constant>0\;\;\mbox{on}\;\;
     \Gamma_\varrho^+.
   \end{array}\right.\end{eqnarray}

As being verified in Section 4, the functions of form
\begin{eqnarray} \label{5-15} u(x)=
c \left(\sin \frac{m_1\pi}{l_1}x_l\right)\cdots
  \left(\sin
  \frac{m_{n-1} \pi } {l_{n-1}}x_{n-1}\right)
 Y (x_n)\end{eqnarray}
 are the solutions of the problem
(\ref{5-13}), where $m_1, \cdots, m_{n-1}$ are positive integers,
and  $Y(x_n)$ is given by (\ref{4--7}). Since the functions in
(\ref{5-15})
 have derivatives of any order in $D$, it follows from Proposition 3.9 and Theorem 3.4 that they
  are eigenfunctions of the transformation $G^0$ with eigenvalues $(1+\gamma)^{-1}$,
 where
\begin{eqnarray} \label{5-16}   \quad \quad \; \gamma  =\frac{2\eta l_n}{\varrho l_n}
  \left(\frac{(\sinh \eta l_n)\cosh \eta l_n - \eta l_n}{\sinh^2 \eta \l_n -\eta^2
  l_n^2}\right),\quad \;\;  \eta= \left[ \sum_{i=1}^{n-1}
  \big(\frac{m_i\pi}{l_i}\big)^2 \right]^{1/2}.\end{eqnarray}
    Note that the normal derivatives
    \begin{eqnarray} \label{5;17} \frac{\partial u}{\partial \nu} =
   c\left(\sin \frac{m_1 \pi}{l_1}  x_1\right)
   \cdots \left(\sin \frac{m_{n-1} \pi}{l_{n-1}} x_{n-1}\right), \end{eqnarray}
 when $m_1, \cdots, m_{n-1}$ run through all positive integers (see, Section 4), form
  a complete system of orthogonal functions in $L^2_\varrho(\Gamma_\varrho^+)$.
It follows from Proposition 3.5 that if $m_1,
  \cdots, m_{n-1}$ run through all positive integers,
  then the functions (\ref{5-15})
 form an orthogonal basis of the subspace of ${\mathcal{K}}^0$, spanned by the
 eigenfunctions of $G^0$, corresponding to
  positive eigenvalues. That is, when
   $m_1, \cdots,  m_{n-1}$ run through all
  positive integers, then $(1+\gamma)^{-1}$,
   where $\gamma$ is given by (\ref{5-16}), runs through all positive eigenvalues of $G^0$.

 Similarly, for the problem (\ref{5-14}), the eigenfunctions $\{u_k\}$ of the operator
 $G^d$ on ${\mathcal{K}}^d$, corresponding to non-zero eigenvalues,
  form  an orthogonal basis
 of the subspace of ${\mathcal{K}}^d$.
 The non-zero eigenvalues  of $G^d$ are $\mu^d_k=(1+\kappa_k)^{-1}$, where $\kappa_k$
 is the $k$-th Steklov eigenvalue of (\ref{5-14}).

In order to give the upper bound estimate of $A^d(\tau)$, we further
introduce the following Steklov
 eigenvalue problem
\begin{eqnarray} \label {5.*.15}  \left\{ \begin{array}{ll} \triangle^2 u = 0\quad \;
   \mbox{in}\;\; D,\\
   u=0 \; \; \mbox{on}\;\; \Gamma_\varrho^+, \quad \,
   u=\frac{\partial u}{\partial \nu}=0 \;\;\mbox{on}\;\;
   \Gamma^{l_n},\\
      \frac{\partial u}{\partial \nu}=0 \;\; \mbox{and}\;\; \frac{\partial (\Delta u)}{\partial \nu}= 0
      \;\; \mbox{on}\;\; \Gamma-(\Gamma_\varrho^+\cup \Gamma^{l_n}),\\
    \triangle u+\varsigma \varrho \frac{\partial u}{\partial \nu}=0 \;\;
    \mbox{on}\;\; \Gamma_\varrho^+, \quad \,  \varrho =constant>0\;\;\mbox{on}\;\;
     \Gamma_\varrho^+.
   \end{array}\right.\end{eqnarray}
 Let $\varsigma_k$ be the $k$-th eigenvalue of (\ref{5.*.15}). By Theorem 3.8, we have
 \begin{eqnarray}\label{5.-.15}  \varsigma_k\le \kappa_k\quad \;\; \mbox{for all }\;\;
 k\ge 1.\end{eqnarray}
   We define \begin{eqnarray} \label{5..0}\mu_k^f=\frac{1}{1+\varsigma_k}, \quad \;\; A^f(\tau)
   = \sum_{\mu_k^f\ge (1+\tau)^{-1}} 1.\end{eqnarray}
It follows from (\ref{5.-.15}) and (\ref{5..0}) that
\begin{eqnarray} \label{5*1} A^d (\tau)\le A^f(\tau)\quad \; \mbox{for all}\;\; \tau.\end{eqnarray}

  We know (cf. Section 4) that the problem (\ref{5.*.15}) has the solutions of form
 \begin{eqnarray}  \label{5-21} u(x)=c \left(\cos \frac{m_1 \pi}{l_1} x_1\right)\cdots
 \left(\cos \frac{m_{n-1} \pi}{l_{n-1}} x_{n-1}\right)
  Z (x_n), \end{eqnarray}
where $m_1, \cdots, m_{n-1}$ are non-negative integers with
$\sum_{i=1}^{n-1}m_{i}\ne 0$, and $Z (x_n)$ is given by
(\ref{4-9}).
 This implies that if $m_1, \cdots, m_{n-1}$  run through all non-negative integers with
  $\sum_{i=1}^{n-1}m_{i} \ne 0$, then
\begin{eqnarray} \quad \quad \; \label{5--16}   \varsigma  =\frac{2\beta l_n}{\rho l_n}
  \left(\frac{(\sinh \beta l_n)\cosh \beta l_n - \beta l_n}{\sinh^2 \beta \l_n -\beta^2 l_n^2}\right),
  \quad \; \beta =  \left[ \sum_{i=1}^{n-1}
  \big(\frac{m_i\pi}{l_i}\big)^2 \right]^{1/2}\end{eqnarray}
 runs throughout
 all eigenvalues of problem (\ref{5.*.15}).

 We first compute the asymptotic behavior of $A^f (\tau)$.
By (\ref{5..0}), (\ref{5--16}) and the argument as in p.$\,$44 of
\cite{We4} or p.$\,$373 of \cite{CH} or p.$\,$51-53 of \cite{Sa},
$A^f(\tau)=$the number of $(n-1)$-tuples $(m_1, \cdots, m_{n-1})$
satisfying the inequality
\begin{eqnarray} \quad \quad \; \label{5-22'}  \frac{2\beta l_n}{\varrho l_n}
  \left(\frac{(\sinh \beta l_n)\cosh \beta l_n -
  \beta l_n}{\sinh^2 \beta \l_n -\beta^2 l_n^2}\right)\le \tau,
  \end{eqnarray}
  where $m_1, \cdots,
 m_{n-1}$ are non-negative integers with $\sum_{i=1}^{n-1} m_i \ne
 0$.
   By setting \begin{eqnarray} \label {005.10} t(s)=2s \left(
\frac{(\sinh s)\cosh s -s}{\sinh^2 s -s^2}\right),\end{eqnarray}
  we see that $$\lim_{s\to +\infty} t(s)/s= 2.$$
 We claim that
for all $s\ge 1$, \begin{eqnarray*} t'(s) =2\left[ \frac{-3 s(\sinh^2 s) +
3s^2 (\sinh s)\cosh s + (\sinh^3 s) \cosh s -s^3(\sinh^2 s +\cosh^2
s)}{(\sinh^2 s-s^2)^2}\right]>0.\end{eqnarray*}
  In fact, let \begin{eqnarray*}
\theta(s)=-3 s(\sinh^2 s) + 3s^3 (\sinh s)\cosh s + (\sinh^3 s) \cosh
s -s^3(\sinh^2 s +\cosh^2 s).\end{eqnarray*} Then
\begin{eqnarray*} &&\theta (1)>0, \,\, \mbox{and}\\
  &&\theta'(s)= -3 (\sinh^2 s) -4 s^3(\sinh s)\cosh s
+  3 (\sinh^2 s) \cosh^2 s +\sinh^4 s\\
&& \quad \quad \; =4(\sinh s)[\sinh^3 s -s^3 \cosh s]>0 \quad
\,\mbox{for}\;\; s\ge 1,\end{eqnarray*}
 This implies that $\theta(s)>0$ for $s\ge 1$.
  Thus, the function ${t}(s)$ is increasing in $[1, +\infty)$.
  Denote by $s={h}(t)$ the inverse of function
 ${t}(s)$ for $s\ge 1$. Then $$\lim_{t\to +\infty} \frac{{h}(t)}{t}=\frac{1}{2}.$$
Furthermore, we can easily check that
\begin{eqnarray} \label{5;;23} h(t)\sim \frac{t}{2} +O(1) \quad \,
\mbox{as}\,\, t\to +\infty.\end{eqnarray}
  Note that, for $s\ge 1$, the inequalities $t(s)\le t$ is equivalent to
  $s\le h(t)$.
 Hence (\ref{5-22'}) is equivalent to
\begin{eqnarray*}  \beta l_n  \le h(\varrho l_n \tau), \end{eqnarray*}
which can be written as
  \begin{eqnarray} \label{5-23}  \sum_{i=1}^{n-1}
  (m_i/l_i)^2 \le \left[\frac{1}{\pi l_n} h(\varrho l_n \tau)\right]^2,
  \quad \, m_i=0, 1,2,\cdots. \end{eqnarray}
 We consider the $(n-1)$-dimensional ellipsoid
\begin{eqnarray*} \label{5-24}  \sum_{i=1}^{n-1}
  (z_i/l_i)^2 \le \left[\frac{1}{\pi l_n} h(\varrho l_n \tau)\right]^2. \end{eqnarray*}
  Since $A^f(\tau)+1\,$ just is the number of those $(n-1)$-dimensional unit cubes of the
  $z$-space that have corners whose
  coordinates are non-negative integers in the ellipsoid (see, VI. \S4 of \cite{CH}).
 Hence $A^f(\tau)+1$ is the sum of the volumes of these cubes.
  Let $V(\tau)$ denote the volume and $T(\tau)$ the
  area of the part of the ellipsoid situated in
  the positive octant $z_i\ge 0,\, i=1,\cdots, n-1$. Then
\begin{eqnarray} \label {5;;25} V(\tau) \le A^f(\tau)+1 \le V(\tau)+ (n-1)^{\frac{1}{2}}
T(\tau),\end{eqnarray}  where $(n-1)^{\frac{1}{2}}$ is the diagonal
length of the unit cube (see, \cite{CH} or \cite{Sa}).  Since
  \begin{eqnarray*} \label{5-25}  V(\tau)= \omega_{n-1}
  2^{-(n-1)} l_1\cdots l_{n-1}\left[\frac{h(\varrho l_n
  \tau)}{\pi l_n}\right]^{(n-1)},
 \end{eqnarray*}
  by (\ref{5;;23}), we get that
  \begin{eqnarray}  \label{5;;26} \quad \qquad \;\; V(\tau)\sim \omega_{n-1}
 (4\pi)^{-(n-1)} l_1\cdots l_{n-1} \varrho^{n-1}
 \tau^{n-1} +O(\tau^{n-2}), \quad \; \mbox{as}\;\; \tau\to +\infty.\end{eqnarray}
   Note that
\begin{eqnarray} \label{5;;27}
  T(\tau)\sim \mbox{constant}\cdot \tau^{n-2}.\end{eqnarray}
    It follows that
      \begin{eqnarray*} \lim_{\tau\to +\infty} \frac{A^f (\tau)}{\tau^{n-1}}= \omega_{n-1}
  (4\pi)^{-(n-1)} l_1 \cdots l_{n-1} \varrho^{n-1}, \end{eqnarray*}
  i.e.,
  \begin{eqnarray} \label{5.2.8}  A^f(\tau) \sim \frac{\omega_{n-1}}{(4\pi)^{(n-1)}}
  |\Gamma_\varrho^+|\varrho^{n-1} \tau^{n-1},\;\; \mbox{as}\;\; \tau \to +\infty,\end{eqnarray}
 where $|\Gamma_\varrho^+|$ denotes the area of the face $\Gamma_\varrho^+$.

Next, we consider $A^0(\tau)$.
    Similarly,
\begin{eqnarray} \quad \quad \; \label{5-22}  \frac{2\eta l_n}{\varrho l_n}
  \left(\frac{(\sinh \eta l_n)\cosh \eta l_n - \eta l_n}{\sinh^2 \eta \l_n -\eta^2 l_n^2}\right)\le \tau,
  \end{eqnarray}
   is equivalent to
  \begin{eqnarray*} \eta l_n \le {h}(\varrho l_n \tau),\end{eqnarray*}
  i.e.,
  \begin{eqnarray*} \sum_{i=1}^{n-1} \big[m_i /l_i\big]^2 \le \left( \frac{h
  (\varrho l_n \tau)}{\pi l_n}\right)^2, \quad \,  m_i= 1,2, 3, \cdots.\end{eqnarray*}
    Similar to the argument for $A^f(\tau)$, we find  (see also, \cite{Liu} or \S4 of \cite{CH}) that
    \begin{eqnarray*}  \# \{ (m_1, \cdots, m_{n-1})\big|
  \sum_{i=1}^{n-1} \big(\frac{m_i}{l_i}\big)^2 \le \big( \frac{h
  (\varrho l_n \tau)}{\pi l_n}\big)^2,\;\, m_i= 1,2,3,\cdots \big\} \\
     \sim
   \frac{\omega_{n-1}}{(4\pi)^{(n-1)}}
  |\Gamma_\varrho^+|\varrho^{n-1} \tau^{n-1} \;\; \mbox{as}\;\; \tau\to +\infty.\qquad \quad\end{eqnarray*}
  i.e.,
    \begin{eqnarray} \label {5-30}   \lim_{\tau \to +\infty}
    \frac{A^0(\tau)}{\tau^{n-1}}= \frac{\omega_{n-1}}{(4\pi)^{(n-1)}}
  |\Gamma_\varrho^+|\varrho^{n-1}.\end{eqnarray}
  Noting that $\varrho=0$ on $\Gamma_\varrho -\Gamma_\varrho^+$,
   by (\ref{5-66}), (\ref{5*1}), (\ref{5.2.8}) and (\ref{5-30}), we have
  \begin{eqnarray} \label{5-31} A(\tau) \sim \frac{\omega_{n-1}\tau^{n-1}}{(4\pi)^{(n-1)}  }
  \int_{\Gamma_{\varrho}} \varrho^{n-1}ds \quad \;
   \mbox{as}\;\; \tau\to +\infty. \end{eqnarray}

\vskip 0.63 true cm

\noindent {\bf 5.3.  $\,$A cylinder $D$ whose base is an
 $n$-polyhedron of ${\Bbb R}^{n-1}$ having $n-1$ orthogonal plane surfaces and
  $g_{ij}=\delta_{ij}$}.

\vskip 0.38 true cm

 \noindent  {\bf Lemma 5.1.}  {\it Let $D^{(r)}
 =\Gamma^{(r)}_\varrho \times [0, l_n]$, $\,r=1,2$,
    where
    $\Gamma^{(1)}_\varrho=\{(x_1,\cdots, x_{n-1})\in {\Bbb R}^{n-1}\big|
    \, x_i\ge 0\;\;\mbox{for}\;\; 1\le i\le n-1,\,\,\mbox{and}\,\, \sum_{i=1}^{n-1} \frac{x_i}{l_i}\le 1\}$,
  and $\Gamma^{(2)}_\varrho$ is an $(n-1)$-dimensional
  cube with side length $l=\max_{1\le i\le n-1} l_i$.
     Assume that
    $\Gamma^{(r)}_{00}=\Gamma^{(r)}_\varrho\times \{l_n\}$, $\; r=1,2$.
     Assume also that  $\varrho$ is a positive constant
     on $\Gamma^{(r)}_\varrho$, $\; r=1,2$. If $l<l_n$, then
      \begin{eqnarray} \label {7-01}\varsigma_k^f(D^{(1)}) \ge  \varsigma_k^f(D^{(2)}) \quad \; \mbox{for}
    \;\;  k=1,2,3,\cdots,\end{eqnarray}
 where $\varsigma_k^f (D^{(r)})$ (similar to $\varsigma$ of (\ref{3-30}) in Theorem 3.8)
  is the $k$-th Steklov eigenvalue for
 the domain $D^{(r)}$ .}

\vskip 0.14 true cm

 \noindent  {\bf Proof.} \ \  Let $v^{(r)}_k$  be the $k$-th Neumann
 eigenfunction
 corresponding to $\alpha^{(r)}_k$  for the $(n-1)$-dimensional
 domain
 $\Gamma^{(r)}_\varrho$, $(r=1,2)$, \ i.e.,
 \begin{eqnarray} \label {7-02} \left\{ \begin{array} {ll}
   \triangle v^{(r)}_k +
 \alpha^{(r)}_k v^{(r)}_k =0 \;\; &\mbox{in}\;\; \Gamma^{(r)}_\varrho,\\
  \frac{\partial v^{(r)}_k}{\partial \nu} =0 \,\, &\mbox{on}\;\; \partial
  \Gamma^{(r)}_\varrho. \end{array} \right.\end{eqnarray}
    Put \begin{eqnarray*}u^{(r)}(x) = \big(v^{(r)}(x_1,
    \cdots, x_{n-1})\big)(Z^{(r)}(x_n)) \quad \, \mbox{in}\;\;
   D^{(r)},\end{eqnarray*}
   where $Z^{(r)}(x_n)$ is as in (\ref{4-9}) with $\beta$
   being replaced by
  $\sqrt{\alpha^{(r)}_k}$. It is easy to verify that $u^{(r)}_k(x)$ satisfies
  \begin{eqnarray}  \label{7-08} \left\{\begin{array}
   {ll} \triangle^2 u^{(r)}_k=0 \;\; &\mbox{in}\;\; D^{(r)},\\
      u^{(r)}_k=0 \; \; &\mbox{on}\;\; \Gamma^{(r)}_\varrho,\\
      u^{(r)}_k=\frac{\partial u^{(r)}_k}{\partial \nu}=0 \;\; & \mbox{on}\;\;
      \Gamma_{00}^{(r)}, \\
  \frac{\partial u^{(r)}_k}{\partial \nu}= \frac{\partial
   (\triangle u^{(r)}_k)}{\partial \nu}=0 \;\; &\mbox{on}\;\;
  \partial D^{(r)} -(\Gamma^{(r)}_\varrho \cup \Gamma^{(r)}_{00}),\\
   \triangle u^{(r)}_k+ \varsigma^f_k (D^{(r)}) \, \varrho \,
  \frac{\partial u_k^{(r)}}{\partial \nu}=0 \;\; &\mbox{on}\;\; \Gamma^{(r)}_{\varrho}.
  \end{array}\right.\end{eqnarray}
 with
\begin{eqnarray}  \label{5;;31} \quad \,\; \varsigma_k^f (D^{(r)}) =\frac{2\sqrt{\alpha_k^{(r)}} \, l_n}
{\varrho l_n}
  \left(\frac{(\sinh \sqrt{\alpha_k^{(r)}} \,l_n)\cosh \sqrt{\alpha_k^{(r)}}\, l_n -
  \sqrt{\alpha_k^{(r)}}\, l_n}{\sinh^2 \sqrt{\alpha_k^{(r)}} \,\l_n -\alpha_k^{(r)} l_n^2}\right).
 \end{eqnarray}
 It follows from p.$\,$437-438 of \cite{CH} that
 the $k$-th Neumann eigenvalue $\alpha^{(1)}_k$ for the domain $\Gamma^{(1)}_\varrho$ is
  at least as large as the $k$-th Neumann eigenvalue $\alpha_k^{(2)}$
  for the domain $\Gamma_\varrho^{(2)}$.
  Recalling that
  $2s \left( \frac{(\sinh s)\cosh s -s}{\sinh^2 s
-s^2}\right)$ is increasing when $s\ge 1$,
 we get  $$\varsigma_k^f (D^{(1)}) \ge \varsigma_k^f(D^{(2)}), \;\;
 k=1,2,3,\cdots$$
 if $l<l_n$. Here we have used the fact that $ \sqrt{\alpha_k^{(2)}}\,l_n
 \ge 1$ since any Neumann eigenvalue for $\Gamma_\varrho^{(2)}$ has the form
 $\sum_{i=1}^{n-1} \big(\frac{m_i\pi}{l}\big)^2$.
 In other words, if $l<l_n$, then the number
 $A^f(\tau)$ of eigenvalues less than or equal
 to a given bound $\tau$ for the domain $D^{(1)}$
 is at most equal to the corresponding
 number of eigenvalues for the domain $D^{(2)}$.
 \ \ $\square$

 Similarly, we can easily verify that
 the number $A^f(\tau)$
 of eigenvalues less than or equal
 to a given bound $\tau$
  for an arbitrary $n$-dimensional rectangular parallelepiped $D$
 is never larger than the corresponding number for
 an $n$-dimensional
 rectangular parallelepiped of the same height with its base an
 $(n-1)$-dimensional cube whose side length is at least equal to the largest side length of the
 base of $D$.

\vskip 0.52 true cm

\noindent {\bf 5.4. $\,D$ is a cylinder and $g_{ij}=\delta_{ij}$}.

\vskip 0.3 true cm

Let $D$ be an open $n$-dimensional cylinder in ${\Bbb R}^n$, whose
boundary consists of an $(n-1)$-dimensional cylindrical surface and
two parallel plane surfaces perpendicular to the cylindrical
surface.  Assume that
  $g_{ij}=\delta_{ij}$ in the whole of $\bar D$, that $\Gamma_\varrho$
 includes at least one of the plane surfaces, which we call
 $\Gamma_\varrho^+$, and that $\varrho$ is positive constant on
 $\Gamma_\varrho^+$ and vanishes on $\Gamma_\varrho -\Gamma_\varrho^+$.
  We let the plane surface $\Gamma_\varrho^+$ be situated in the
  plane $x_n=0$ and let another parallel surface $\Gamma^{l_n}$ be situated in the plane
  $\{x\in {\mathbb{R}}^n \big|x_n =l_n\}$. We now divide the
  plane $x_n=0$ into a net of $(n-1)$-dimensional
  cubes, whose faces are parallel to the coordinate-planes in $x_n=0$. Let $\Gamma_1, \cdots,
  \Gamma_p$ be those open cubes in the net, closure of which
  are entirely contained in $\Gamma_\varrho^+$,
   and let $Q_{p+1}, \cdots, Q_q$ be the remaining open cubes, whose closure
     intersect $\Gamma_\varrho^+$. We may let the subdivision  into
   cubes  be so fine that, for every piece of the boundary of $\Gamma_\varrho^+$
   which is contained  in one of the closure cubes, the direction of the normal varies by less than a
   given angle $\vartheta$, whose size will be determined  later.
 (This can be accomplished by repeated halving of the side of cube.)
         We can make the side length $l$ of each
   cube be less than $l_n$.
   Furthermore, let $D_j$, $(j=1, \cdots, p)$, be the open $n$-dimensional
   rectangular parallelepiped with the cube $\Gamma_j$
   as a base and otherwise bounded by the ``upper''
  plane surface $\Gamma^{l_n}$ of the cylinder $\bar D$ and planes parallel to the
  coordinate-planes $x_1=0, \cdots, x_{n-1}=0$ (cf.$\,$\cite{Sa}).

  We define the linear spaces of functions
\begin{eqnarray*} &K=\{u \big| u \in H_0^1 (D)\cap  H^2(D),
\;  \frac{\partial u}{\partial \nu}=0 \;\; \mbox{on}\;\;
\Gamma_{00}\},\\
& K_j^0 =\{u_j \big| u_j \in H_0^1 (D_j)\cap H^2(D_j)\cap
C^\infty(\bar D_j),\;
 \frac{\partial u_j}{\partial \nu}=0 \;\; \mbox{on}\;\;
 \Gamma_j^{l_n},\;\;\mbox{and} \qquad \quad \\
 & \qquad \qquad \qquad \quad \; \; \qquad  \Delta u_j=0\;\;\mbox{on}\;\;
 \partial D_j-(\Gamma_j \cup \Gamma_j^{l_n})
 \}, \quad\;  (j=1,\cdots, p)\end{eqnarray*}
 with the inner products
\begin{eqnarray*} & \langle u, v\rangle =
  \int_{D} (\triangle u)(\triangle v) dR  +\int_{\Gamma_\varrho} \varrho \frac{\partial u}{\partial \nu}
  \, \frac{\partial v}{\partial \nu}\,ds \quad \; \mbox{for}\;\; u,
  v\in K,\\
&\langle u_j, v_j\rangle_j =\langle u_j, v_j \rangle_j^{\star} +[u_j
, v_j]_j  =
  \int_{D_j} (\triangle u_j)(\triangle v_j)dR +\int_{\Gamma_j} \varrho \frac{\partial u_j}{\partial \nu}
  \, \frac{\partial v_j}{\partial \nu}ds \quad \; \mbox{for} \;\; u_j, v_j\in K_j^0,\end{eqnarray*}
   respectively.  Closing $K$ and $K_j^0$
   with respect to the norms $\|u\| =
  \sqrt{ \langle u, u\rangle}$ and $\|u_j\|_j =
  \sqrt{ \langle u_j, u_j\rangle_j}$, we obtain the
   Hilbert spaces ${\mathcal{K}}$ and  ${\mathcal{K}}_j^0\,$ ($j=1, \cdots, p$), respectively.
 Clearly, the bilinear functional
\begin{eqnarray*}  &&[u,v] = \int_{\Gamma_\varrho} \varrho \,
\frac{\partial u}{\partial \nu}\, \frac{\partial v}{\partial \nu}
\, ds\\
 &&  [u_j, v_j]_j = \int_{\Gamma_j} \varrho \,\frac{\partial
u_j}{\partial \nu} \,\frac{\partial u_j}{\partial \nu}\, ds, \quad
\;\; (j=1, \cdots, p), \end{eqnarray*}
  define  self-adjoint, completely continuous transformations $G$ and $G^0_j$ on
  ${\mathcal{K}}$ and ${\mathcal{K}}_j^0$
  by \begin{eqnarray} \label{5-2-1}
   \langle Gu, v\rangle&=&[u,v] \quad \;\; \mbox{for}\;\; u,v
\;\;\mbox{in}\;\; {\mathcal{K}},\\
  \label {5-2/1} {\langle G_j^0 u_j, v_j\rangle}_j &=&[u_j,v_j]_j
  \quad \;\; \mbox{for}\;\;
  u_j \,\,\mbox{and}\;\; v_j\,\, \mbox{in}\;\; {\mathcal{K}}_j^0,\end{eqnarray}
 respectively.
 By defining a space \begin{eqnarray*} {\mathcal{K}}^0
 =\sum_{j=1}^p \oplus {\mathcal{K}}_j^0 =\{u^0 \big| u^0=
u_1+\cdots + u_p, \,\, u_j\in {\mathcal{K}}_j^0\}\end{eqnarray*}
    with its inner product
   \begin{eqnarray}  \label{5-2-2}\langle u^0, v^0\rangle =\sum_{j=1}^p
   \langle u_j, v_j\rangle_j, \;\;
   \end{eqnarray}
 we find that the space ${\mathcal{K}}^0$ becomes a Hilbert space.
 If  we define the transformation $G^0$ on ${\mathcal{K}}^0$ by
 \begin{eqnarray} \label{5-2/2} G^0 u^0=G_1^0 u_1+ \cdots+ G^0_p u_p\quad \;\; \mbox{for}\;\;
 u^0 =u_1+\cdots+ u_p\,\, \mbox{in}\,\,
  {\mathcal{K}}^0,\end{eqnarray}
  we see  that $G^0$ is a self-adjoint, completely continuous
  transformation on ${\mathcal{K}}^0$. If we put
  \begin{eqnarray}\label {5-2-3} [u^0, v^0]= \sum_{j=1}^p [u_j,v_j]_j, \end{eqnarray}
  we find by (\ref{5-2/1})---(\ref{5-2-3}) that
  \begin{eqnarray} \label{5-2-4} \langle G^0 u^0, v^0\rangle
  = [ u^0, v^0] \quad \; \,\mbox{for all}\;\;
  u^0\;\;\mbox{and}\;\; v^0 \;\; \mbox{in}\;\; {\mathcal{K}}^0.\end{eqnarray}
   Let us define a mapping of ${\mathcal{K}}^0$ into $\mathcal{K}$.
    Let $u^0=u_1+ \cdots+ u_p, \, u_j\in H_j^0$, be an element of ${\mathcal{K}}^0$
     and define
     \begin{eqnarray} \label{5-2-5} u=\Pi^0 u^0,\end{eqnarray}
     where $u(x)= u_j(x)$, when $x\in \bar D_j$, and $u(x)=0$,
     when $x\in \bar D -\cup_{j=1}^p \bar D_j$. Clearly $u\in \mathcal{K}$
     and thus (\ref{5-2-5}) defines a transformation $\Pi^0$ of ${\mathcal{K}}_1^0 \oplus \cdots \oplus
     {\mathcal{K}}_p^0$ into
     $\mathcal{K}$.
  It is readily seen that
\begin{eqnarray}\label{5-2-6}
      [\Pi^0 u^0, \Pi^0 v^0] =[u^0, v^0] \quad \; \mbox{for all}\;\;
      u^0 \;\; \mbox{and}\;\; v^0 \;\; \mbox{in}\;\; {\mathcal{K}}^0.
      \end{eqnarray}
and
\begin{eqnarray}
      \label{5-2-7}\langle G^0 u^0,  v^0\rangle
      =\langle G\Pi^0 u^0, \Pi^0 v^0\rangle  \quad \; \mbox{for all}\;\;
      u^0 \;\; \mbox{and}\;\; v^0 \;\; \mbox{in}\;\; {\mathcal{K}}^0.
      \end{eqnarray}
  By (\ref{5-2-6}) and (\ref{5-2-7}), we find by applying Proposition 2.3 that
 \begin{eqnarray*} \mu_k^0 \le \mu_k \quad \; \mbox{for}\;\; k=1,2,3, \cdots.\end{eqnarray*}
 Therefore \begin{eqnarray} \label {5-38} A^0(\tau)\le A(\tau).\end{eqnarray}
 The definition of $G^0$ implies that
 \begin{eqnarray} \label{5-2-8}  G^0 {\mathcal{K}}_j^0
 \subset {\mathcal{K}}_j^0, \quad \; (j=1, \cdots, p),\end{eqnarray}
  and
\begin{eqnarray} \label{5-2-9} G^0 u^0 = G_j^0 u^0,
 \quad \, \mbox{when}\;\; u^0\in  {\mathcal{K}}_j^0.\end{eqnarray}
 From (\ref{5-2-3}), (\ref{5-2-4}), (\ref{5-2-8}), (\ref{5-2-9}) and Proposition 2.4, we obtain
\begin{eqnarray} \label {5-40} A^0(\tau) =\sum_{j=1}^p A^0_j (\tau),\end{eqnarray}
  where $A_j^0 (\tau)$ is the number of eigenvalues
  of the transformation
   $G_j^0$ on ${\mathcal{K}}_j^0$
 which are greater or equal to $(1+\tau)^{-1}$.
 Because  $\bar D_j, (j=1, \cdots, p)$, is an $n$-dimensional
    rectangular parallelepiped
   we find  by (\ref{5-30}) that
\begin{eqnarray} \label{5-2-10} A^0_j (\tau)  \sim \omega_{n-1}(4\pi)^{-(n-1)}
 |\Gamma_j| \varrho^{n-1}  \tau^{n-1}\quad \;\; \mbox{as}\;\; \tau\to +\infty,\end{eqnarray}
where $|\Gamma_j|$ denotes the area of the face $\Gamma_j$ of $D_j$.
By (\ref{5-40}) and (\ref{5-2-10}) we infer that
\begin{eqnarray} \label{5-2-11} A^0 (\tau) \sim \omega_{n-1}(4\pi)^{-(n-1)}
  \sum_{j=1}^p |\Gamma_j| \varrho^{n-1}  \tau^{n-1}\quad \, \mbox{as}\;\; \tau\to +\infty.\end{eqnarray}

\vskip 0.12 true cm

 Next, we shall calculate the upper estimate of $A(\tau)$.
  Let $\bar P_j, (j=p+1, \cdots,q)$, be the $n$-dimensional rectangular parallelepiped
   with the cube $\bar Q_j$ as a base and otherwise
 bounded by the ``upper'' plane surface $\Gamma^{l_n}$ of the cylinder $\bar D$ and planes parallel to the
 coordinate-planes $x_1=0, \cdots, x_{n-1}=0$.
  The intersection $\bar P_j \cap \bar D$ is a cylinder $\bar D_j, (j=p+1, \cdots, q)$,
  with $\bar \Gamma_j:=\bar Q_j \cap \bar \Gamma_\varrho^+$  as a base. Obviously
  \begin{eqnarray} \bar D=\sum_{j=1}^q \bar D_j.\end{eqnarray}

  We first define the linear spaces of functions
  \begin{eqnarray*}& K^d =\{u\big| u\in  H^2(D),
  \; u=0\,\, \mbox{on} \;\; \Gamma_\varrho, \; u=\frac{\partial u}{\partial \nu}
  =0 \,\, \mbox{on}\,\,
   \Gamma^{l_n})\},\\
 & K_j^d =\{u_j\big| u_j\in  H^2(D_j),
  \; u_j=0\,\, \mbox{on} \;\; \Gamma_j, \,\,  u_j=\frac{\partial u_j}{\partial \nu}
  =0 \,\, \mbox{on}\,\,
   \Gamma^{l_n}_j\},\; \; (j=1, \cdots, q)\end{eqnarray*}
  with the inner products
    \begin{eqnarray} \label{5-3-1}  & \langle u, v\rangle = \int_{D} (\triangle u)(\triangle v) dR
  +\int_{\Gamma_\varrho}  \varrho\, \frac{\partial u}{\partial \nu}\, \frac{\partial v}{\partial \nu}
  \, ds, \end{eqnarray}
  and
  \begin{eqnarray} \label {5-3.2}  \langle u_j, v_j\rangle_j = \int_{D_j}  (\triangle u_j)(\triangle v_j) dR
  +\int_{\Gamma_j}  \varrho \,\frac{\partial u_j}{\partial \nu}\, \frac{\partial v_j}{\partial \nu}
  \, ds,\end{eqnarray}
  respectively.
  Closing $K^d$ and $K_j^d$ with respect to the norms $\|u\|=\sqrt{\langle u, u\rangle}$
  and $\|u_j\|_j
  =\sqrt{\langle u_j, u_j\rangle_j}$,
 we obtain Hilbert spaces
 ${\mathcal{K}}^d$ and ${\mathcal{K}}_j^d, (j=1, \cdots, q)$, and then we define the Hilbert
 space
\begin{eqnarray} \label{5-3-2} {\mathcal{K}}^d = \sum_{j=1}^q \oplus {\mathcal{K}}^d_j=\{
  u^d\big| u^d =u_1+\cdots+ u_q, \; \; u_j \in
{\mathcal{K}}_j^d\}\end{eqnarray}
 with its inner product
  \begin{eqnarray} \label{5-3-3}\langle u^d,v^d \rangle
   = \sum_{j=1}^q \langle u_j, v_j\rangle_j.    \end{eqnarray}
  The bilinear functional
\begin{eqnarray} \label{5-3-4} [u_j, v_j]_j = \int_{\Gamma_j}  \varrho \frac{\partial u_j}{\partial \nu}\,
\frac{\partial v_j}{\partial \nu}\, ds, \quad \;\; (j=1, \cdots, q),\end{eqnarray}
 define a self-adjoint, completely continuous transformation $G_j^d$
 on ${\mathcal{K}}_j^d$ given by
\begin{eqnarray} \label{5-3-5}
      \langle G^d_j u_j,  v_j\rangle_j =[u_j, v_j]_j \quad \; \mbox{for all}\;\;
      u_j \;\; \mbox{and}\;\; v_j \;\; \mbox{in}\;\; {\mathcal{K}}_j^d.
      \end{eqnarray}
       The self-adjoint, completely continuous transformation $G^d$ on ${\mathcal{K}}^d$ is defined by
 \begin{eqnarray} \label{5-3-6}G^d u^d=G_1^d u_1+\cdots+ G^d_q u_q\quad \;\; \mbox{for}\;\;
 u^d =u_1+\cdots+ u_q\,\, \mbox{in}\,\,
  {\mathcal{K}}^d.\end{eqnarray}
  With
    \begin{eqnarray} \label{5-3-7} [u^d, v^d]= \sum_{j=1}^q [u_j,v_j]_j, \end{eqnarray}
   we find by (\ref{5-3-3}), (\ref{5-3-5})---(\ref{5-3-7}) that
  \begin{eqnarray} \label{5-3-8} \langle G^d u^d, v^d\rangle = [u^d, v^d] \quad \; \,\mbox{for all}\;\;
  u^d \,\, \mbox{and}\,\, v^d \;\; \mbox{in}\;\; {\mathcal{K}}^d.\end{eqnarray}
  Let us define a mapping  $\Pi$ of $\mathcal{K}$ into ${\mathcal{K}}^d$.
 Let $u\in K(D)$, and put
 \begin{eqnarray*} u^d =\Pi u=u_1+ \cdots + u_q,\end{eqnarray*}
where $u_j(x)= u(x)$, when $x\in \bar D_j$. It can be easily
 verified that
  \begin{eqnarray} \label{5-3-9} \langle \Pi u, \Pi v\rangle =\langle u, v\rangle
      \quad \, \mbox{for all}\;\; u \;\; \mbox{and}\;\; v \;\; \mbox{in}\;\;
      {\mathcal{K}}. \end{eqnarray}
and
 \begin{eqnarray} \label{5-3-10} \langle Gu, v\rangle =\langle G^d \Pi u, \Pi v\rangle
        \quad \, \mbox{for all}\;\; u \;\; \mbox{and}\;\; v \;\; \mbox{in}\;\;
            \mathcal{K}.\end{eqnarray}
 From (\ref{5-3-9})---(\ref{5-3-10}), with the aid of Proposition 2.3, we obtain
  \begin{eqnarray*} \mu_k \le \mu_k^d \quad \; \mbox{for}\;\; k=1,2,3, \cdots,\end{eqnarray*}
 and hence
  \begin{eqnarray} \label {5..56} A(\tau)\le A^d(\tau).\end{eqnarray}
 By $G^d {\mathcal{K}}_j^d\subset {\mathcal{K}}_j^d$, $\;(j=1, \cdots,
 q)$, and $G^d u^d =G^d_j u^d$ when $u^d\in {\mathcal{K}}_j^d$,  we
 get
\begin{eqnarray} A^d(\tau) =\sum_{j=1}^q A^d_j (\tau),\end{eqnarray}
  where $A_j^d (\tau)$ is the number of eigenvalues
   of the transformation
   $G_j^d$ on ${\mathcal{K}}_j^d$
   which are
   greater than or equal to $(1+\tau)^{-1}$.
 Further, we define $A_j^f (\tau)$ similar to (\ref{5.*.15}) and (\ref{5..0}), i.e.,
 \begin{eqnarray*}  A^f_j(\tau)
   = \sum_{\mu_k^f\ge (1+\tau)^{-1}} 1\; \quad \;\mbox{with}\;\;
   \mu_k^f= \frac{1}{1+\varsigma_k},\end{eqnarray*}
 where  $\varsigma_k$ is  the $k$-th Steklov eigenvalue of the following problem
  \begin{eqnarray*} \label {5..15}  \left\{ \begin{array}{ll} \triangle^2 u_j = 0\quad \;
   \mbox{in}\;\; D_j,\\
   u_j=0 \; \; \mbox{on}\;\; \Gamma_j, \quad \,
   u_j=\frac{\partial u_j}{\partial \nu}=0 \;\;\mbox{on}\;\;
   \Gamma^{l_n}_j,\\
     \frac{\partial u}{\partial \nu}= \frac{\partial (\Delta u_j)}{\partial \nu}= 0
      \;\; \mbox{on}\;\; \partial D_j-(\Gamma_j\cup \Gamma^{l_n}_j),\\
    \triangle u_j+\varsigma \varrho \frac{\partial u_j}{\partial \nu}=0 \;\;
    \mbox{on}\;\; \Gamma_j, \quad \,  \varrho =constant>0\;\;\mbox{on}\;\;
     \Gamma_\varrho^+.
   \end{array}\right.\end{eqnarray*}
    From Theorem 3.8, it follows that
 \begin{eqnarray*} \varsigma_k\le \kappa_k\quad \;\; \mbox{for all }\;\;
 k\ge 1,\end{eqnarray*}
 and hence
 \begin{eqnarray}  \label{5---1} A^d_j(\tau)\le A^f_j(\tau)\quad \;\mbox{for all}\;\;
   \tau\;\;\mbox{and}\;\;j=1,\cdots, q,\end{eqnarray}
  where $\frac{1}{1+\kappa_k}$ is the $k$-th eigenvalue of the transformation
 $G^d_j$.
   Since $\bar D_j, (j=1, \cdots, p)$, is an $n$-dimensional rectangular
   parallelepiped,
  we find from
 (\ref{5.2.8}) that
 \begin{eqnarray}\label{5..58}  A^f_j (\tau) \sim \omega_{n-1}(4\pi)^{-(n-1)}
 |\Gamma_j| \varrho^{n-1}  \tau^{n-1}, \quad \;(j=1, \cdots, p). \end{eqnarray}

 It remains to estimate $A^f_j(\tau), (j\ge p+1)$.
 According to the argument in p.$\,$438-440 of \cite{CH},
  each of the $(n-1)$-dimensional domains $\Gamma_j$ is
 bounded either by $n-1$ orthogonal plane surfaces of the partition
 (the diameter of the intersection  of any two plane surfaces lies
   between $l$ and $3l$), and an $(n-2)$-dimensional
  surface of the boundary (see, in two dimensional case, Figure 5 of p.$\,$439 of \cite{CH}),
 or by $2n-3$ orthogonal plane surfaces of the partition (the diameter of the intersection of any
 two plane surfaces lies between $l$ and $3l$),
 and a surface of the boundary $\partial \Gamma_\varrho$ (see, in two dimensional case,
   Figure 6 of p.$\,$439 of \cite{CH}).
The number $q-p$ is evidently smaller than a constant $C/l^{n-2}$,
where $C$ is independent of $l$ and depends essentially on the area
of
 the boundary $\partial \Gamma_\varrho$.
  Now, we take any point on the  boundary surface of $\Gamma_j$
  and take the tangent plane through it.
   This tangent plane together with the plane parts of $\partial \Gamma_j$
   bounds an $n$-polyhedron of ${\Bbb R}^{n-1}$ with a vertex at which $n-1$
   orthogonal plane surfaces meet  (see,
    Figure 5 of p.$\,$439 of \cite{CH} in two dimensions),
   e.g., if $\vartheta$ is sufficiently small it forms an $(n-1)$-dimensional
   $n$-polyhedron of ${\Bbb R}^n$
   with a vertex having $n-1$ orthogonal plane surfaces (the diameter of the intersection  of any two
   plane surfaces is also smaller than $4l$), or else an  $(n-1)$-dimensional $2(n-1)$-polyhedron
   of ${\Bbb R}^{n-1}$ (see, Figure 6 of p.$\,$439 of \cite{CH} in two dimensional case),
      the diameter of the intersection of any two plane surfaces (except for the top inclined plane surface)
      of the $2(n-1)$-polyhedron is also smaller than $4l$;
     The shape of the result domain depends on the type to which  $\bar \Gamma_j$ belongs.
   We shall denote the result domains by $S'_j$.  The domain $\Gamma_j$
   can always be deformed into the domain $S'_j$ by a transformation
   of the form (\ref{02-001}), as defined in Section 2.  In the case of
   domains of the first  type,  let the intersection point of $n-1$ orthogonal plane surfaces
    be the pole of a system  of pole coordinates
   $r$, $\theta_1$,  $\theta_2, \cdots, \theta_{n-2}$, and let $r=f(\theta_1,
   \theta_2, \cdots, \theta_{n-2})$ be the equation of the boundary
   surface of $\Gamma_\varrho$, $r=h(\theta_1,
   \theta_2, \cdots, \theta_{n-2})$
  the equation of the inclined plane surface of the $n$-polyhedron of
  ${\Bbb R}^{n-1}$ having a vertex
  of $n-1$ orthogonal plane surfaces. Then the
  equations
   $$\theta_1'=\theta_1, \quad  \theta'_2=\theta_2, \;\cdots,\;\,
      \theta'_{n-2}=\theta_ {n-2},  \quad r'= r\,\frac{h(\theta_1,
   \theta_2, \cdots, \theta_{n-2})}{f(\theta_1,
   \theta_2, \cdots, \theta_{n-2})}$$
  represents a transformation of the domain  $\Gamma_j$ into
  the $n$-polyhedron  $S'_j$ of ${\Bbb R}^{n-1}$.
     For a domain of the second type, let $x_{n-1}=
  h (x_1, \cdots, x_{n-2})$  be the equation of top plane surface of the
  $2(n-1)$-polyhedron  and let $x_{n-1}=f(x_1, \cdots, x_{n-2})$
  be the  equation of the boundary  surface of $\Gamma_\varrho$.
  We then consider the  transformation
  $$ x'_1=x_1, \;\; \cdots, \quad x'_{n-2}=x_{n-2}, \quad
   x'_{n-1}= x_{n-1}\,\frac{h(x_1, \cdots, x_{n-2})}
   {f(x_1, \cdots, x_{n-2})}. $$
If we assume that the side length $l$ of cube in the partition is
sufficiently small, and therefore the rotation of the normal on the
 boundary surface  is taken sufficiently small, then the
transformations considered here evidently have precise the form
(\ref{02-001}),  and the quantity denoted by $\epsilon$ in
(\ref{02-001}) is arbitrarily small.
  From Corollary to Theorem 10 of p.$\,$423 of \cite{CH}, we know that
  there exists a number $\delta>0$ depending on $\epsilon$
 and approaching zero with $\epsilon$, such that
  \begin{eqnarray*} \bigg| \frac{\alpha_k(S'_j)}{\alpha_k(\Gamma_j)}
  -1\bigg|<\delta \quad \, \mbox{uniformly for all}\;\; k,
  \end{eqnarray*}
  where $\alpha_k (\Gamma_j)$ and $\alpha_k (S'_j)$  are the $k$-th Neumann
  eigenvalues of $\Gamma_j$ and $S'_j$, respectively.
  According to the argument as in the proof of Lemma 5.1 (i.e., (\ref{5;;31})), we see
  that
  \begin{eqnarray*}
 \varsigma_k^f (E_j)=\frac{1}{\varrho l_n}\, t(
  l_n \, \sqrt{\alpha_k(\Gamma_j)}), \quad \;  \varsigma_k^f (E'_j)=\frac{1}{\varrho l_n}\, t(
 l_n \, \sqrt{\alpha_k(S'_j)}),\end{eqnarray*}
 where  $t(s)$ is given by (\ref{005.10}), and $\varsigma_k^f (E_j)$
  and $\varsigma_k^f(E'_j)$ (similar to $\varsigma$ of (\ref{3-30})) are the $k$-th Steklov eigenvalue
   for the $n$-dimensional  domains
  $E_j=\Gamma_j\times [0, l_n]$
  and $E'_j= S'_j \times [0, l_n]$, respectively.  Recalling
  that the function
  $t=t(s)$ is continuous and increasing for $s\ge 1$, we get
 that there exists a constant $\delta'>0$ depending on $\epsilon$
 approaching
 zero with $\epsilon$, such that
\begin{eqnarray*} \bigg|\frac{  \varsigma_k^f
(E'_j)}{\varsigma_k^f (E_j)}-1\bigg|<\delta'.\end{eqnarray*}
 In other words, the corresponding
$k$-th Steklov eigenvalues for the $n$-dimensional
 domains $E_j=\Gamma_j\times [0, l_n]$ and $E'_j=S'_j\times
[0,l_n]$ differ only by a  factor which itself differs by a small
amount from $1$, uniformly  for all $k$. Therefore, the same is true
also for the corresponding
 numbers $A^f_{E_j} (\tau)$ and $A^f_{E'_j} (\tau)$ of the eigenvalues less or equal to the bound
 $\tau$.

The  domain $E'_j$ is either a cylinder whose base is an
$n$-polyhedron of ${\Bbb R}^{n-1}$ having $(n-1)$ orthogonal plane
surfaces with its largest side length smaller than $4l$ or a
cylinder whose base is a
 combination of such an $n$-polyhedron of ${\Bbb R}^{n-1}$
 and an $(n-1)$-dimensional cube with
sides smaller than $3l$; it follows
 from the estimates for $E'_j$ (cf. (\ref{5;;25})---(\ref{5;;27})) and Lemma
 5.1 that if $l$ is taken sufficiently small,  the number
$A^f_{E_j}(\tau)$ from some $\tau$ on satisfies the inequality
  \begin{eqnarray*} \label{07-00}  A^f_{E_j}(\tau) < C_1 l^{n-1}
  \tau^{n-1} +C_2 l^{n-2} \tau^{n-2} \end{eqnarray*}
where $C_1$, $C_2$ are constants, to be chosen suitably.
  Thus, $A_{E_j}^f (\tau)$ can be  written as
  $A_{E_j}^f (\tau)=\theta (C_3 l^{n-1} \tau^{n-1}+C_4
  l^{n-2}\tau^{n-2})$, where $\theta$ denotes a number between $-1$
  and $+1$ and $C_3, C_4$ are constants independent of $l,j$ and
  $\tau$. It follows that $$\sum_{j=p+1}^q A_{E_j}^f(\tau)=\tau^{n-1}
  \big[ \theta C_3 (q-p)l^{n-1} +\theta C_4 (q-p)l^{n-2}
  \,\frac{1}{\tau}\big].$$
  As pointed out before, $(q-p)l^{n-2}<C$; therefore, for sufficiently small
  $l$, $\,(q-p)l^{n-1}$ is arbitrarily
 small and we have  the asymptotic relation
 \begin{eqnarray} \label{5.-01} \lim_{\tau\to +\infty}
 \sum_{j=p+1}^q \frac{A_{E_j}^f (\tau)}{\tau^{n-1}}
 =\varpi (l),  \end{eqnarray}
 where $\varpi (l)\to 0$ as $l\to 0$.
 For, we may choose the quantity $l$ arbitrarily, and  by taking
   a sufficiently small fixed $l$, make the factors of $\tau^{n-1}$ in
  the previous equalities arbitrarily close to zero for sufficiently
  large $\tau$.
 Since \begin{eqnarray}   \label{5.-02}  A^d_{E_j}(\tau)\le A^f_{E_j}(\tau)
 \quad \, \mbox{for}\;\; j=p+1, \cdots, q,\end{eqnarray}
 we get
 \begin{eqnarray} \label{5.-03} \lim_{\tau\to +\infty}
 \frac{ \sum_{j=p+1}^q A_{E_j}^d (\tau)}{\tau^{n-1}}\le
 \lim_{\tau\to +\infty}
 \frac{ \sum_{j=p+1}^q A_{E_j}^f (\tau)}{\tau^{n-1}}=\varpi (l).  \end{eqnarray}
  From (\ref{5-38}), (\ref{5-2-11}), (\ref{5..56}), (\ref{5---1}),
  (\ref{5..58}), (\ref{5.-01}),  (\ref{5.-02}) and
  (\ref{5.-03}), we obtain
  \begin{eqnarray} \label {5.63} &&
 \omega_{n-1}(4\pi)^{-(n-1)} \varrho^{n-1}
  \sum_{j=1}^p |\Gamma_j| \le
 \underset {\tau\to \infty} {\underline{\lim}}\, \frac{A(\tau)}{\tau^{n-1}}
\le \overline{\lim_{\tau\to \infty}} \,\frac{A(\tau)}{\tau^{n-1}}
\\
 && \;\;\quad \qquad \quad\;\;  \le
 \bigg(\omega_{n-1}(4\pi)^{-(n-1)}
   \varrho^{n-1}\sum_{j=1}^p |\Gamma_j| \bigg)+\varpi (l) \nonumber.
   \end{eqnarray}
 Letting $l\to 0$, we immediately see that
$\sum_{j=1}^p |\Gamma_j|$ tends to the area $|\Gamma_\varrho|$ of
 $\Gamma_\varrho$ and $\lim_{l\to 0} \varpi (l)=0$.
 Therefore,
 (\ref{5.63}) gives
\begin{eqnarray} \label{5;;66}A (\tau) \sim
\frac{\omega_{n-1}}{(4\pi)^{(n-1)}}
 |\Gamma_\varrho^+| \varrho^{n-1}  \tau^{n-1} \quad \; \mbox{as}\;\; \tau\to +\infty,\end{eqnarray}
or
\begin{eqnarray}\label{5-90} A (\tau) \sim \frac{\omega_{n-1}\tau^{n-1}}{(4\pi)^{(n-1)}}
   \int_{\Gamma_\varrho}  \varrho^{n-1} ds \quad \mbox{as} \;\; \tau\to +\infty.\end{eqnarray}

\vskip 0.26 true cm

\noindent {\bf Remark 5.2.}  \ \  In the above argument, we first
made the assumption that the boundary $\partial \Gamma_\varrho$ of
$\Gamma_\varrho$  was smooth. However,  the corresponding discussion
and result remain essentially valid if $\partial \Gamma_\varrho$ is
composed of a finite number of $(n-2)$ dimensional smooth surfaces.

 \vskip 1.39 true cm

\section{Proofs of main results}

\vskip 0.45 true cm

 \noindent  {\bf Lemma 6.1.} \ \  {\it  Let $g_{ij}$ and $g'_{ij}$
  be two metric tensors on manifold $\mathcal{M}$ such that
    \begin{eqnarray} \label {006-1}
 \big|g^{ij} -g'^{ij}\big|< \epsilon, \quad \, i,j=1, \cdots, n
    \end{eqnarray}
and
 \begin{eqnarray} \label {006-2} \bigg| \frac{1}{\sqrt{|g|}} \,\frac{\partial}{\partial x_i}
   \big(\sqrt{|g|} g^{ij}\big)\,- \, \frac{1}{\sqrt{|g'|}} \,\frac{\partial}{\partial x_i}
   \big(\sqrt{|g'|} g'^{ij}\big)\bigg|\le \epsilon, \quad \,
   i,j=1,\cdots,n   \end{eqnarray}
   for all points in $\bar D$, where $D$ is a bounded domain in $\mathcal{M}$.
    Let \begin{eqnarray*} \mu_1\ge \mu_2\ge \cdots \ge \mu_n \ge
  \cdots >0 \;\; \mbox{and}\;\;  \mu'_1\ge \mu'_2\ge \cdots \ge \mu'_n \ge
  \cdots >0\end{eqnarray*}
  be positive eigenvalues of $G$ and $G'$, respectively,
  where $G$ and $G'$ are given by \begin{eqnarray*} \langle Gu,
v\rangle =\int_{\Gamma_\varrho} \varrho \, \frac{\partial
u}{\partial\nu} \,\frac{\partial v}{\partial \nu}\, ds, \quad\;
\mbox{for}\;\; u \,\, \mbox{and}\;\; v
\,\, \mbox{in}\;\; \mathcal{K}, \\
\langle G'u, v\rangle' =\int_{\Gamma_\varrho} \varrho \,
\frac{\partial u}{\partial \nu} \,\frac{\partial v}{\partial\nu}\,
ds', \quad\; \mbox{for}\;\; u \,\, \mbox{and}\;\; v \,\,
\mbox{in}\;\; {\mathcal{K}}'.\end{eqnarray*}
  Then
  \begin{eqnarray}\label {006-3} && (1+{\tilde M}\epsilon)^{-(n+1)/2}
  \left( \max\{(1+\epsilon M), (1+{\tilde M}\epsilon)^{(n+1)/2}\}\right)^{-1} \mu_k \le \mu'_k
 \\
  && \quad \;\;\le
 (1+{\tilde M}\epsilon)^{(n+1)/2} \left( \min\{ (1-\epsilon M),
 (1+{\tilde M}\epsilon)^{-(n+1)/2}\} \right)^{-1} \mu_k, \nonumber
\\ && \;\;
 \qquad
 \qquad \;\;\mbox{for}\;\;
 k=1,2,3,\cdots,\nonumber \end{eqnarray}
where $\tilde M$ and $M$ are constants depending only on $g$, $g'$,
$\frac{\partial g_{ij}}{\partial x_l}$, $\frac{\partial
g^{ij}}{\partial x_l}$, $\frac{\partial g'_{ij}}{\partial x_l}$,
$\frac{\partial g'^{ij}}{\partial x_l}$
  and $\bar D$.}

 \vskip 0.24 true cm

\noindent  {\bf Proof.} \
 It follows from (\ref{006-1})  that
 there exists a positive constant $\tilde M$ independent of
 $\epsilon$ and depending only on $g^{ij}$, $g'^{ij}$ and $\bar D$ such that
 \begin{eqnarray*} (1+\epsilon \tilde M)^{-1} \sum_{i,j=1}^{n}
 g^{ij} t_it_j \le \sum_{i,j=1}^{n}
 g'^{ij} t_it_j \le (1+\epsilon \tilde M) \sum_{i,j=1}^{n}
 g^{ij} t_it_j  \end{eqnarray*}
 for all points in $\bar D$ and all real numbers $t_1, \cdots, t_n$.
 Thus we have
 \begin{eqnarray*} (1+{\tilde M}\epsilon)^{-n/2} \sqrt{|g|}\le \sqrt{|g'|}
\le (1+{\tilde M}\epsilon)^{n/2} \sqrt{|g|}, \end{eqnarray*}
   which implies (see
  p.$\,$64-65 of \cite{Sa}) that \begin{eqnarray*} (1+{\tilde M}\epsilon)^{-n/2} dR \le dR'
  \le (1+{\tilde M}\epsilon)^{n/2} dR
\end{eqnarray*}
 and \begin{eqnarray*}   (1+\tilde{M}\epsilon)^{-(n+1)/2} ds \le ds' \le
 (1+{\tilde M}\epsilon)^{(n+1)/2} ds.\end{eqnarray*}
 Thus \begin{eqnarray} \label {006-4} (1+{\tilde M}\epsilon)^{-(n+1)/2} [u,u]\le [u, u]' \le
 (1+{\tilde M}\epsilon)^{(n+1)/2} [u, u].\end{eqnarray}
  Putting   $$\omega_{ij}=g'^{ij}- g^{ij}, \quad \; \theta_{ij}
  =\frac{1}{\sqrt{|g'|}} \,\frac{\partial}{\partial x_i}
   \big(\sqrt{|g'|} g'^{ij}\big)\,- \, \frac{1}{\sqrt{|g|}} \,\frac{\partial}{\partial x_i}
   \big(\sqrt{|g|} g^{ij}\big),$$
  we  immediately see that \begin{eqnarray*}
   \max_{x\in \bar D}|\omega_{ij}|\le \epsilon \quad \;
    \mbox{and}\;\; \max_{x\in \bar D} |\theta_{ij}|\le
   \epsilon.\end{eqnarray*}
 Thus, for any $u\in {\mathcal{K}}^0 (D)$ or $u\in {\mathcal{K}}^d
 (D)$,
we have \begin{eqnarray*} \quad \triangle_{g'} u
   = \sum_{i,j=1}^n (\omega_{ij} +g^{ij} ) \frac{\partial^2 u}{\partial
x_i\partial x_j} + \sum_{i,j=1}^n \left[\theta_{ij}+
  \frac{1}{\sqrt{|g|}} \,\frac{\partial}{\partial x_i}
   \big(\sqrt{|g|} g^{ij}\big)\right] \frac{\partial
   u}{\partial x_j}, \end{eqnarray*}
so that \begin{eqnarray*}\label{006-5}  \triangle_{g'} u -
 \triangle_g u = \sum_{i,j=1}^n \left[ \omega_{ij}  \frac{\partial^2 u}{\partial
  x_i\,\partial x_j} +\theta_{ij} \,
   \frac{\partial u}{\partial
   x_j}\right]\nonumber\end{eqnarray*}
  It follows that
\begin{eqnarray*} | \triangle_{g'} u -
 \triangle_g u| \le \epsilon \big(M_1|\nabla_g^2 u|+M_2|\nabla_g
 u|\big),\end{eqnarray*}
 where $|\nabla^2_g u|^2$ is defined in an invariant ways as
  \begin{eqnarray*} |\nabla_g^2 u|^2 =\nabla^l \nabla^k u \,
  \nabla_l \nabla_k u = g^{pl} g^{kq} \left(\frac{\partial^2
  u}{\partial x_k \partial x_l} -\Gamma_{kl}^m \frac{\partial
  u}{\partial x_m}\right)\left(\frac{\partial^2
  u}{\partial x_p \partial x_q} -\Gamma_{pq}^r \frac{\partial
  u}{\partial x_r}\right),\end{eqnarray*}
 and $M_1$ and $M_2$ are constants depending only on $g$, $g'$,
$\frac{\partial g_{ij}}{\partial x_l}$, $\frac{\partial
g^{ij}}{\partial x_l}$, $\frac{\partial g'_{ij}}{\partial x_l}$,
$\frac{\partial g'^{ij}}{\partial x_l}$ and
 $\bar D$.
 Thus, \begin{eqnarray} \label {006-7}
 \int_{D} | \triangle_{g'} u -
 \triangle_g u|^2 dR  \le 2 \epsilon^2 \left(M_1^2\int_D |\nabla_g^2 u|^2dR+M_2^2\int_D |\nabla_g
 u|^2dR   \right).\end{eqnarray}

Recall also (see, Section 2) that $\Lambda_1^0(D)\ge
\Lambda_1^d(D)$, where
\begin{eqnarray} \label {006-10}
  \Lambda_1^0 (D) =\inf_{v\in K^0(D), \;\int_D |\nabla_g v|^2 dR =1}
 \,\frac{\int_D|\triangle_g v|^2 dR }{\int_D |\nabla_g v|^2
dR},\end{eqnarray}
  \begin{eqnarray} \label {006-11}
  \Lambda_1^d (D) =\inf_{v\in K^d(D), \;\int_D |\nabla_g v|^2 dR =1}
\, \frac{\int_D|\triangle_g v|^2 dR }{\int_D |\nabla_g v|^2
dR},\end{eqnarray}
 and $K^0(D)$ and  $K^d(D)$ are  as in  Section 3.
  Let \begin{eqnarray} \label {006-8}
  \Theta_1^0 (D) =\inf_{v\in K^0(D), \;\int_D |\nabla_g^2 v|^2 dR =1}
 \,\frac{\int_D|\triangle_g v|^2 dR }{\int_D |\nabla_g^2 v|^2
dR},\end{eqnarray}
  \begin{eqnarray} \label {006-9}
  \Theta_1^d (D) =\inf_{v\in K^d(D), \;\int_D |\nabla_g^2 v|^2 dR =1}
 \,\frac{\int_D|\triangle_g v|^2 dR }{\int_D |\nabla_g^2 v|^2
dR}.\end{eqnarray}
  Clearly,  $\Theta_1^0(D)\ge \Theta_1^d
  (D)$.
  As in the proofs of Lemmas 2.1, 2.2, it is easy to prove that
  the existence of the minimizers to (\ref{006-11}) and (\ref{006-9}),
  respectively. Therefore, we have that $\Lambda_1^d(D)>0$ and $\Theta_1^d
  (D)>0$ (Suppose by contradiction that $\Lambda_1^d(D)=0$ and $\Theta_1^d
  (D)=0$.
  Then $\triangle_g u=0$ in $D$ for the corresponding minimizer $u\in K^d(D)$ in both cases.
 By applying Holmgren's uniqueness theorem for the minimizer $u\in K^d(D)$
   in both cases, we immediately see that
   $u\equiv 0$ in $D$.
 This contradicts the assumption $\int_D |\nabla_g u|^2
 dR=1$ or $\int_D |\nabla_g^2 u|^2 dR =1$ for the minimizer $u\in K^d(D)$ in the corresponding cases).
  Combining these inequalities, we obtain
\begin{eqnarray*} \label {006-12}
 \int_{D} | \triangle_{g'} u -
 \triangle_g u|^2 dR  \le 2\epsilon^2 \left(\frac{M^2_1}{\Theta_1^0(D)}  +\frac{M^2_2}{\Lambda_1^0
 (D)}\right)
  \int_D |\triangle_g
 u|^2, \quad \mbox{for}\;\; u\in K^0(D)\end{eqnarray*}
and \begin{eqnarray*} \label {006-13}
 \int_{D} | \triangle_{g'} u -
 \triangle_g u|^2 dR  \le 2 \epsilon^2 \left(\frac{M^2_1}{\Theta_1^d(D)}  +\frac{M^2_2}{\Lambda_1^d (D)}
  \right) \int_D |\triangle_g
 u|^2 \quad \mbox{for}\;\; u\in K^d(D).\end{eqnarray*}
  Thus we have  that, for all $u\in K^0(D)$ or $u\in
 K^d(D)$,
 \begin{eqnarray*} \label {006-14}
    (1-\epsilon M)\int_D | \triangle_g u|^2 dR
  \le    \int_{D} |\triangle_{g'} u|^2 dR
   \le (1+\epsilon M) \int_D |\triangle_{g} u|^2 dR,\end{eqnarray*}
where $M$ is a constant depending only $g$, $g'$, $\frac{\partial
g_{ij}}{\partial x_l}$, $\frac{\partial g^{ij}}{\partial x_l}$,
$\frac{\partial g'_{ij}}{\partial x_l}$, $\frac{\partial
g'^{ij}}{\partial x_l}$ and $\bar D$. That is,
\begin{eqnarray} \label {006-15} (1-\epsilon M) \langle u,u\rangle^\star \le \langle u, u\rangle'^\star \le
 (1+\epsilon M)  \langle u, u\rangle^\star.\end{eqnarray}
By (\ref{006-4}) and (\ref{006-15}) we obtain that, for all $u\in
K^0 (D)$ or $u\in K^d(D)$,
\begin{eqnarray*} \label {006-16} &&\frac{(1+{\tilde M}\epsilon)^{-(n+1)/2}[u,u]}{
 \big(\max \{ (1+\epsilon M), (1+{\tilde M}\epsilon)^{(n+1)/2}\}\big) \big(
  \langle u, u\rangle^* +[u,u]\big)} \le
   \frac{[u,u]'}{\langle u,u\rangle'^\star
   +[u,u]'}\\
&&\quad \;\;\le
 \frac{(1+{\tilde M}\epsilon)^{(n+1)/2} [u,u]}
 {\big(\min \{(1-\epsilon M), (1+{\tilde M}\epsilon)^{-(n+1)/2}\}\big) \big( \langle u,
 u\rangle^\star  +[u,u]\big)},
 \end{eqnarray*}
 which implies (\ref{006-3}). \ \  $\square$

\vskip 0.23 true cm
 \noindent  {\bf Remark 6.2.} \  Let $\tilde \Gamma$ and $\Gamma$ be two bounded domains
 in ${\Bbb R}^{n-1}$, let $\tilde \Gamma$ is similar to $\Gamma$ (in
 the  elementary sense of the term; the length of any line in $\tilde \Gamma$
  is to the corresponding length in $\Gamma$ as $h$  to $1$), and
  let $\Gamma_{00} =\Gamma\times \{\sigma\}$ and $\tilde \Gamma_{00} =\tilde \Gamma\times \{h\sigma \}$.
   It is easy to verify that \begin{eqnarray*} \Lambda_1^d(\tilde D) = h^{-2} \Lambda_1^d(D), \quad
   \;\, \Theta_1^d(\tilde D)=\Theta_1^d (D), \end{eqnarray*}
where $D=\Gamma_i\times [0, \sigma]$, $\tilde D=\tilde \Gamma \times
[0, h\sigma]$, and $\Lambda_1^d (D)$ and $\Theta_1^d(D)$ are defined
as in (\ref{006-11}) and
 (\ref{006-9}), respectively.

\vskip 0.26 true cm

 \noindent  {\bf Lemma 6.3.} \ \  {\it
  Let $G$ and  $G'$ be the continuous linear transformations defined by
\begin{eqnarray*} \langle Gu, v\rangle = \int_{\Gamma_\varrho}
\varrho\,\frac{\partial u}{\partial \nu}\, \frac{\partial
v}{\partial \nu}\, ds \quad \; \; \mbox{for} \;\; u\;\;
\mbox{and}\;\;
 v \;\; \mbox{in}\;\; K^0(D)\;\;\mbox{or}\;\; K^d (D)\end{eqnarray*}
 and \begin{eqnarray*} \langle G'u, v\rangle' = \int_{\Gamma_\varrho}
 \varrho'\,
 \frac{\partial u}{\partial \nu}\,\frac{\partial v}{\partial \nu}\, ds
\quad \; \; \mbox{for} \;\; u\;\; \mbox{and}\;\;
 v \;\; \mbox{in}\;\; K^0(D)\;\;\mbox{or}\;\; K^d
 (D),\end{eqnarray*}
 respectively.
 Let \begin{eqnarray*}
 \mu_1\ge \mu_2\ge\cdots\ge \mu_k \ge \cdots >0\quad \,
 \mbox{and}\;\;
\mu'_1\ge \mu'_2\ge\cdots\ge \mu'_k \ge \cdots >0\end{eqnarray*} be
the positive eigenvalues of $G$ and $G'$, respectively.
 If $\varrho\le \varrho'$, then
 \begin{eqnarray} \label{a-1} \mu_k\le \mu'_k \quad \;\mbox{for}\;\;
 k=1,2,3,\cdots.\end{eqnarray}}

\vskip 0.22 true cm

\noindent  {\bf Proof.} \ Since $\varrho\le \varrho'$, we see that
for any $u\in K^0(D)$ or $K^d(D)$,
\begin{eqnarray*} \frac{\langle Gu, u\rangle}{\langle u, u\rangle}
=\frac{\int_{\Gamma_\varrho} \varrho \left(\frac{\partial
u}{\partial \nu}\right)^2 \, ds }{\langle u, u\rangle^*
+\int_{\Gamma_\varrho} \varrho \left(\frac{\partial u}{\partial
\nu}\right)^2\, ds}\le \frac{\int_{\Gamma_\varrho} \varrho'
\left(\frac{\partial u}{\partial \nu}\right)^2 \, ds }{\langle u,
u\rangle^* +\int_{\Gamma_\varrho} \varrho' \left(\frac{\partial
u}{\partial \nu}\right)^2\, ds} =\frac{\langle G'u,
u\rangle'}{\langle u, u\rangle'},\end{eqnarray*} which
 implies (\ref{a-1}).  \ \ $\square$

 \vskip
0.586 true cm

\noindent  {\bf Proof  of Theorem 1.1.} \ \  a)  \ \  First, let
$(\mathcal{M}, g)$ be a real analytic
 Riemannian manifold, and let the boundary $\partial \Omega$ of $\Omega$ be real
 analytic.
    We divide the domain $\bar \Omega$
 into subdomains in the following manner.
  It is clear that the boundary $\partial \Omega$ of the domain $\Omega$
  is the union of a finite number of closed
 pieces $\bar \Gamma_1, \cdots, \bar \Gamma_p$ (without common
 inner point on the surface).
    Let $U$ be a coordinate neighborhood which contains
 $\bar \Gamma_j$, let $x_i=x_i (Q)$ and $a_i= a_i(\nu_Q)$ be the coordinates
 of a point $Q$ in $\bar \Gamma_j$
  and the interior Riemannian normal $\nu_Q$ at $Q$, respectively.
   We define the subdomain  $D_j$ and surface $\Gamma_j^\sigma$ by
  \begin{eqnarray*} \label{6-1}   D_j = \{P\big| x(P)=x(Q)+\xi_n a(\nu_Q),
  \quad  Q\in \Gamma_j, \,\, 0<\xi_n< \sigma\}\end{eqnarray*}
  and \begin{eqnarray*} \Gamma_j^\sigma =\{ P\big| x(P)=x(Q)+\sigma \,
    a(\nu_Q),\; Q\in \Gamma_j\},\nonumber
  \end{eqnarray*}
 where $\sigma$ is a positive constant. The closure of  $D_j$ is
 \begin{eqnarray} \bar D_j = \{P\big| x(P) =x(Q)+\xi_n a(\nu_Q),
 \; Q\in \bar \Gamma_j, \; 0\le \xi_n \le \sigma\}.\end{eqnarray}
By the assumption,
 each $\bar \Gamma_j$, which is contained in a
 coordinate neighborhood, can be represented by equations
  \begin{eqnarray} \label {0006-1} x_i =\psi_i (\xi_1, \cdots,
  \xi_{n-1})\end{eqnarray}
  with real analytic functions $\psi_i$, i.e., it is the imagine of the closure
  $\bar \Upsilon_j$ of an open domain $\Upsilon_j$ of ${\Bbb R}^{n-1}$.
Hence, if $\sigma$ is sufficiently small, the definitions have a
sense  and the formula \begin{eqnarray}  x(P)=x(Q) +\xi_n a(\nu_Q),
\quad \,
  Q\in \bar \Gamma_j, \;\; 0\le \xi_n \le \sigma\end{eqnarray}
  defines a real analytic homeomorphism of a neighborhood of the
  image of $\bar D_j$ in ${\Bbb R}^n$  given by the coordinates $x$
  and a neighborhood $U_j$ of the  closed cylinder
  $\bar F_j$ in ${\Bbb R}^n$ defined by $\bar F_j= \{\xi \big| (\xi_1,
  \cdots, \xi_{n-1})\in \bar \Upsilon_j, \, 0\le \xi_n \le \sigma\}$.
    Moreover, the domains $\bar D_1, \cdots, \bar D_p$
 have no common inner points and the remainder
 $D_0=\Omega -\cup_{j=1}^p \bar D_j$ of $\Omega$ has a finite number of connected parts.
 Note that  the boundary of $\bar D_0$ contains
 no part of $\partial \Omega$.

 Let us define the spaces $K=K(\Omega)$, $\mathcal{K}$ and the transformation
$G$ as in Section 5.
   We shall investigate the asymptotic behavior of $A(\tau)$ with
  regard to transformation $G$ on space $\mathcal{K}$.
 Moreover, we define the function spaces
   \begin{eqnarray*} & K_j^0 =\{u_j\big|u_j \in H_0^1 (D_j)\cap
   H^2(D_j)\cap C^\infty(\bar D_j),\,
\frac{\partial u_j}{\partial \nu}=0 \;\; \mbox{on}\;\;
\Gamma_j^{\sigma},\;\; \mbox{and} \quad \quad \; \\
&\qquad \qquad \qquad \qquad  \quad\; \Delta u=0\;\;\mbox{on}\;\;
\partial D_j -(\Gamma_j \cup \Gamma_j^{\sigma})\},\\
   & H_0^0 = \{u_0 \big|u_0\in H_0^1 (D_0)\cap H^2(D_0), \;
    \frac{\partial u_0}{\partial \nu} =0 \;\; \mbox{on}\;\;
   \partial D_0\},\\
   & K_j^d =
  \{u_j\big|u_j\in   H^2(D_j), \, u_j=0 \;\; \mbox{on}\;\;
  \Gamma_j, \, u_j=\frac{\partial u_j}{\partial \nu}=0\;\;\mbox{on}\;\;
  \Gamma_j^\sigma\},\\
 & \, \quad (j=0, 1, \cdots, p), \end{eqnarray*}
     and the bilinear functionals
     \begin{eqnarray}  \langle u_j, u_j \rangle_j^{\star} =\int_{D_j} |\triangle_g u_j |^2 dR,
     \quad \, (j=0, \cdots, p), \end{eqnarray}
      \begin{eqnarray} [u_j, v_j]_j = \int_{\Gamma_j} \varrho \frac{\partial u_j}{\partial \nu} \,
   \frac{\partial u_j}{\partial \nu} \, ds,  \quad (j=1, \cdots, p),  \quad \, [u_0, v_0]=0, \end{eqnarray}
   and
  \begin{eqnarray} \langle u_j, v_j\rangle_j=  \langle u_j,
  v_j\rangle_j^{\star} +[u_j, v_j]_j, \quad (j=0, \cdots, p),\end{eqnarray}
  where $u_j, v_j \in K_j^0$ or $K_j^d$. Closing $K_j^0$ and $K^d_j$ with respect to
 the norm $|u_j|_j=\sqrt{\langle u_j, u_j\rangle_j}$,
 we get the Hilbert spaces ${\mathcal{K}}_j^0$ and ${\mathcal{K}}_j^d$,
  $\,(j=0, \cdots, p)$.
 In the same manner as in Section 5 we can define the  Hilbert
  spaces ${\mathcal{K}}^0$ and  ${\mathcal{K}}^d$, and then define the
 positive, completely continuous transformations
  $G^0$, $G^d$, $G_j^0$
 and $G_j^d$ on  ${\mathcal{K}}^0$,
    ${\mathcal{K}}^d$, ${\mathcal{K}}_j^0$
    and ${\mathcal{K}}_j^d$, respectively. Consequently, we can prove
 \begin{eqnarray} \label{6-100} A^0(\tau) \le A(\tau)\le A^d(\tau)
 \quad \, \mbox{for all}\;\; \tau,\end{eqnarray}
and \begin{eqnarray} \label{6-101} A^0(\tau) =\sum_{j=0}^p
A_j^0(\tau),\quad\;
  A^d(\tau)=\sum_{j=0}^p A_j^d(\tau),\end{eqnarray}
 where $A^0(\tau)$, $A^d(\tau)$, $A_j^0(\tau)$ and $A_j^d(\tau)$
 are the numbers of eigenvalues
 of the transformations  $G^0$, $G^d$, $G_j^0$ and $G^d_j$
 on ${\mathcal{K}}^0$, ${\mathcal{K}}^d$, ${\mathcal{K}}_j^0$ and
 ${\mathcal{K}}_j^d$
  which are greater than or equal to $(1+\tau)^{-1}$,
  respectively.

  Since  $[u_0, u_0]_0=0$ for all $u_0\in K_0^0$ or $K_0^d$  and  $\langle
  G_0^0 u_0, u_0\rangle_0 =\langle G_0^d u_0, u_0\rangle_0 =[u_0, u_0]_0$, we immediately find
   that $G_0^0=G_0^d=0$,
   so that  $A_0^0(\tau) =A_0^d(\tau) =0$, ($\tau\ge 0$).
  Thus we need estimate $A_j^0(\tau)$ and $A_j^d(\tau)$ for those domains $D_j$, where
  $\int_{\Gamma_j} \varrho \, ds>0$.

 We can choose a finer subdivision of $\partial \Omega$ by
subdividing the domains $\bar \Upsilon_j$ into smaller  ones, e.g.
 by means of a cubical net in the coordinates $\xi$.
 According to p.$\,$71 of \cite{Sa}, by performing a linear transformation
 $\Phi$ of the coordinates
 we can choose a new coordinate system $(\eta)$ such that
 \begin{eqnarray*}
 g^{il} (\bar \eta)=\delta^{il},\;\quad (i,l=1, \cdots, n),
\end{eqnarray*}
 for one point $\bar \eta\in T_j$, where $T_j:= \Phi(\Upsilon_j)$.
 Setting
   $\phi_i=\psi_i\circ \Phi^{-1}$ and $\tilde a_i= a_i \circ
   \Phi^{-1}$, we see that  \begin{eqnarray}
 \label {0.6-0.1}   & x_i(P)= \phi_i (\eta_1, \cdots, \eta_{n-1})+\eta_n
 \,\tilde a_i(\nu (\eta_1, \cdots, \eta_{n-1})),\\
   & \mbox{for}\;\; (\eta_1, \cdots, \eta_{n-1})\in \bar T_j,\;\,
    0\le \eta_n \le \sigma  \nonumber \end{eqnarray}
 defines a real analytic homeomorphism from  $\bar E_j$ to the image of $\bar
 D_j$, where $\bar E_j=\{\eta=(\eta_1$,$\cdots$, $\eta_n)\big| (\eta_1$,
 $\cdots$, $\eta_{n-1})\in \bar T_j, \,
  0\le \eta_n\le \sigma\}$ is a  cylinder
in ${\Bbb R}^n$ (This can also be realized by choosing a
(Riemannian) normal
 coordinates system at the point $\bar \eta\in T_j$ for the manifold
 $(\mathcal{M},g)$ (see, for example, p.$\,$77 of \cite{Le})
  such that $a(\nu(\eta))=(0, \cdots, 0, 1)$ and by using the mapping
  (\ref{0.6-0.1}).)
 If we denote the new
 subdomains of
 $\partial \Omega$ by $\bar \Gamma_j$ as before, it is clear that we can always choose
  them and $\sigma$ (i.e., by letting  $\sigma$ sufficiently small
  and further making a finer subdivision of
   $\partial \Omega$, see p.$\,$71 of \cite{Sa}), so that,
\begin{eqnarray} \label{0.6-0.2} |g'^{il} (\eta')-g^{il}(\bar
\eta)|<\epsilon, \quad \, i,l=1,\cdots, n,\end{eqnarray}
 \begin{eqnarray} \label {0006-2}  &&\bigg| \frac{1}{\sqrt{|g(\eta')|}} \,\frac{\partial}{\partial x_i}
   \big(\sqrt{|g(\eta')|} g^{il}(\eta')\big)\,- \,
   \frac{1}{\sqrt{|g(\bar \eta)|}} \,\frac{\partial}{\partial x_i}
   \big(\sqrt{|g(\bar \eta)|} g^{il}(\bar \eta)\big)\bigg|<\epsilon, \\
   &&\,\qquad \;\;\qquad  \; \qquad \, i,l=1,\cdots,n,\nonumber
   \end{eqnarray}
for any given $\epsilon>0$, and all points $\eta' \in {\bar E}_j$.
 The inequalities (\ref{0.6-0.2}) imply that
\begin{eqnarray} \label {0006-1} (1+{\tilde M_j}\epsilon)^{-1}\sum_{i=1}^n  t_i^2  \le \sum_{i,l=1}^n g^{il}
    (\eta') t_i t_l \le (1+{\tilde M_j}\epsilon)\sum_{i=1}^n
   t_i^2\end{eqnarray}
  for all points $\eta'\in {\bar E}_j$ and all real numbers $t_1, \cdots,
  t_n$, where ${{\tilde M}}_j$ is a positive constant depending only on
  $g^{il}$ and $\bar E_j$ (cf. Lemma 6.1).  This and formula (128) of \cite {Sa} say that
  \begin{eqnarray} \label{6>>1} (1+\epsilon{\tilde M}_j)^{-n/2} |T_j|\le |\Gamma_j|
  \le (1+\epsilon{\tilde M}_j)^{n/2} |T_j|,\end{eqnarray}
  where
  \begin{eqnarray*}  |\Gamma_j|=
 \int_{T_j} \sqrt{g(\eta)}\, d\eta_1 \cdots d\eta_{n-1},
 \quad \, |T_j|=\int_{T_j} d\eta \cdots d\eta_{n-1}\end{eqnarray*}
 are the Riemannian and Euclidean areas of $\Gamma_j$ and
  $|T_j|$, respectively,
  \vskip 0.22 true cm

  Next, we consider the Hilbert spaces ${\mathcal{K}}_j^0$ and
  ${\mathcal{K}}_j^d$. When transported to $\bar E_j$,
   the underlying incomplete function spaces
  $K_j^0$ and $K_j^d$ are
  \begin{eqnarray*} K_j^0 = \{ u\big| u\in H_0^1(E_j)\cap H^2(E_j)\cap C^\infty(\bar E_j),  \,
    \, \frac{\partial u}{\partial \nu}=0 \;\;\mbox{on}\;\; T_j^\sigma, \,\mbox{and}\;\,
    \Delta u=0\;\;\mbox{on}\;\; \partial E_j -(T_j \cup T_j^\sigma)\}
  \end{eqnarray*}
   and \begin{eqnarray*}  K_j^d =\{u_j\big| u_j\in
    H^2(E_j),\; u_j=0 \;\;\mbox{on}\;\; T_j,\; u_j=\frac{\partial u_j}{\partial \nu}=0\;\;
  \mbox{on}\;\; T_j^{\sigma}\}, \end{eqnarray*}
   respectively.
    The inner product, which is similar to Section 5, is defined  by
    \begin{eqnarray*} \langle u, v\rangle_j = \int_{E_j} (\triangle_g u)(\triangle_g v)
     \sqrt{g(\eta)} d\eta_1 \cdots d\eta_n +\int_{T_j}
  \varrho \, \frac{\partial u}{\partial \nu}\, \frac{\partial u}{
  \partial \nu} \sqrt{g(\eta)}\, d\eta_1\cdots \eta_{n-1} \end{eqnarray*}
  and the transformations $G_j^0$ and $G_j^d$ are defined by
  \begin{eqnarray*} \langle G_j^0u, v\rangle_j = \int_{T_j} \varrho
   \frac{\partial u}{\partial \nu}\, \frac{\partial v}{
  \partial \nu} \sqrt{g(\eta)}\, d\eta_1\cdots \eta_{n-1}, \quad \; \mbox{for}\;\; u,v \;\;
  \mbox{in}\;\; {\mathcal{K}}_j^0, \end{eqnarray*}
and
 \begin{eqnarray*} \langle G_j^d u, v\rangle_j = \int_{T_j} \varrho
   \frac{\partial u}{\partial \nu}\, \frac{\partial v}{
  \partial \nu} \sqrt{g(\eta)}\, d\eta_1\cdots \eta_{n-1}, \quad \; \mbox{for}\;\; u,v \;\;
  \mbox{in}\;\; {\mathcal{K}}_j^d, \end{eqnarray*}
respectively.

 Put
\begin{eqnarray}  {\underline{\varrho}}_j =\inf_{\bar \Gamma_j} \varrho \quad\; \mbox{and}\;\;
 \bar \varrho_j =\sup_{\bar \Gamma_j} \varrho,\end{eqnarray}
 and let us introduce the inner products
  \begin{eqnarray*} {\underline{\langle u, v\rangle}}_j =\int_{E_j}
  (\triangle u)(\triangle_g v) d\eta_1 \cdots d \eta_n +
  \int_{T_j} {\underline{\varrho}}_j \frac{\partial u}{\partial \nu}
 \,\frac{\partial v}{\partial \nu}\, d\eta_1\cdots d\eta_{n-1} \end{eqnarray*}
 and
 \begin{eqnarray*} {\overline{\langle u, v\rangle}}_j =\int_{E_j}
  (\triangle u)(\triangle_g v) d\eta_1 \cdots d \eta_n +
  \int_{T_j} {\overline{\varrho}}_j \frac{\partial u}{\partial \nu}
 \,\frac{\partial v}{\partial \nu}\, d\eta_1\cdots d\eta_{n-1} \end{eqnarray*}
in the spaces $K^0_j$ and $K_j^d$, respectively. By closing these
spaces in the corresponding norms, we get Hilbert spaces
${\underline{\mathcal{K}}}_j^0$
 and ${\overline{\mathcal{K}}}_j^d$. Furthermore, we  obtain the
  positive, completely continuous transformations ${\underline{G}}_j^0$ and
  ${\overline{G}}_j^d$ on
 ${\underline{\mathcal{K}}}_j^0$ and
${\overline{\mathcal{K}}}_j^d$, which are given by
\begin{eqnarray} \underline{\langle {\underline{G}}_j^0 u, v\rangle_j}=\int_{T_j}
{\underline{\varrho}}_j \frac{\partial u}{\partial \nu} \, \frac{\partial u}{\partial \nu}
 \, d\eta_1 \cdots d\eta_{n-1}, \quad \; \mbox{for}\;\; u\;\; \mbox{and}\;\; v\;\; \mbox{in}\;\;
 {\underline{\mathcal{K}}}_j^0\end{eqnarray}
 and
\begin{eqnarray} \overline{\langle {\overline{G}}_j^d u, v\rangle_j}=\int_{T_j}
 {\overline{\varrho}}_j \frac{\partial u}{\partial \nu} \, \frac{\partial u}{\partial \nu}
 \, d\eta_1 \cdots d\eta_{n-1}, \quad \; \mbox{for}\;\; u\;\; \mbox{and}\;\; v\;\; \mbox{in}\;\;
 {\overline{\mathcal{K}}}_j^d,\end{eqnarray}  respectively.

  Let $\mu_k(G_j^0)$ be the $k$-th positive eigenvalue of
  $G_j^0$ and so on. According to Lemma 6.1 and Remark 6.2,
  $\Lambda_1^d(D_j)$ and $\Theta_1^d(D_j)$ have uniformly positive lower bound when
  repeated taking finer division of $D$
  (In fact, by repeated halving the side length of
  every rectangular parallelepiped
  in the partition net of the coordinates $\eta$ for each cylinder $E_j$,
  we see that $\Lambda_1^d(D_j)$ will tend to $+\infty$, and
  that  $\Theta_1^d(D_j)$ will have a positive lower bound).
 This implies that the corresponding positive constants ${\tilde M}_j$ and
 $M_j$ have uniformly upper bound when
 we further divide the domain $D$ into finer a division, where
 $\tilde M_j$ is defined as before, and $M_j$ is a constant
  independent of $\epsilon$ and depending only on $g$, $\frac{\partial
 g_{im}}{\partial x_l}$, $\,\frac{\partial g^{im}}{\partial x_l}$
   and $\bar E_j$ as in Lemma 6.1.
  Denote by $c_j(\epsilon)$ the maximum value of
 $(1+\epsilon {\tilde M}_j)^{(n+1)/2}
  \big( \max\{(1+\epsilon M_j)$, $(1+\epsilon {\tilde M}_j )^{(n+1)/2}\}\big)$ and $
   (1+\epsilon{\tilde M}_j)^{(n+1)/2} \left( \min\{ (1-\epsilon M_j),
 (1+\epsilon{\tilde M}_j)^{-(n+1)/2}\} \right)^{-1}$.
Obviously, $c_j(\epsilon)\to 1$ as $\epsilon\to 0$.
  By virtue of (\ref{0.6-0.2}) and (\ref{0006-2}), it follows from Lemmas 6.1 and 6.3 that
\begin{eqnarray} \mu_k(G^d_j)\le c_j(\epsilon)\, \mu_k(\bar G_j^d)\end{eqnarray}
and
\begin{eqnarray} \mu_k(G^0_j)\ge c_j(\epsilon)^{-1}\mu_k({\underline{G}}_j^0),\end{eqnarray}
 so that
\begin{eqnarray*} A_j^d (\tau)\le {\bar A}_j^d
\big( c_j(\epsilon)\tau+c_j(\epsilon)-1\big)\end{eqnarray*} and
\begin{eqnarray} A_j^0 (\tau)\ge {\underline{A}}_j^0
\big(c_j(\epsilon)^{-1}\tau+
  c_j(\epsilon)^{-1} -1\big)\end{eqnarray}
 where ${\bar
 A}^d_j (\tau)$ and
${\underline{A}}^0_j (\tau)$ are
 the numbers of eigenvalues  of the transformation $\bar G_j^d$
 $(1+\tau)^{-1}$, respectively.
 and ${\underline{G}}_j^0$
  which are greater than or equal to
 By (\ref{6-100}) and (\ref{6-101}), we obtain
\begin{eqnarray} \label{6>1} && \sum_j {\underline{A}}_j^0
\big(c_j(\epsilon)^{-1} \tau+ c_j(\epsilon)^{-1}
-1\big)\le A(\tau) \\
 && \; \quad \qquad \le \sum_j {\bar A}_j^d
\big(c_j(\epsilon)\tau +
c_j(\epsilon)-1\big).\nonumber\end{eqnarray}

 Finally,   we shall apply the results of Section 5 to
estimate ${\underline{A}}_j^0(\tau)$
 and ${\overline{A}}_j^d(\tau)$.
Note that \begin{eqnarray} \label{5:0} {\overline{A}}_j^d(\tau)\le
{\overline{A}}_j^f(\tau) \quad\;\,\mbox{for all}\;\;
\tau>0,\end{eqnarray} where ${\bar A}_j^f$ is defined similarly to
(\ref{5.*.15})---(\ref{5..0}).
 It follows from (\ref{5-2-10}), (\ref{5..58}), (\ref{5:0}) and (\ref{5.-03})-(\ref{5;;66}) that
\begin{eqnarray}\label{6>2} \underset {\tau\to
+\infty}{\underline{\lim}}
 \frac{{\underline{A}}_j^0 (\tau)}{\tau^{n-1}}\ge  \omega_{n-1}(4\pi)^{-(n-1)}
 |{T}_j| {\underline{\varrho}}_j^{n-1}    \end{eqnarray}
and \begin{eqnarray} \label{6>3} \underset {\tau\to
+\infty}{\overline{\lim}} \frac{{\bar A}_j^f (\tau)}{\tau^{n-1}} \le
\omega_{n-1}(4\pi)^{-(n-1)}
  |{T}_j| {\bar \varrho}_j^{n-1},\end{eqnarray}
where $|{T_j}|$ is the Euclidean area of $T_j$. By (\ref{6>1}),
(\ref{5:0}), (\ref{6>2}), (\ref{6>3}), (\ref{6-100}), (\ref{6-101})
 and (\ref{6>>1}), we find that
\begin{eqnarray*}  \underset{\tau\to \infty} {\overline{\lim}} A(\tau)\, \tau^{-(n-1)} \le
  \omega_{n-1}(4\pi)^{-(n-1)}{\tilde c}_j(\epsilon) \sum_j
   {\overline{\varrho}}_j^{n-1} |\Gamma_j|,\end{eqnarray*}
and \begin{eqnarray*} \underset{\tau\to \infty} {\underline{\lim}}
A(\tau)
 \,\tau^{-(n-1)} \ge
  \omega_{n-1}(4\pi)^{-(n-1)}{\tilde c}_j(\epsilon)^{-1} \sum_j
   {\underline{\varrho}}_j^{n-1} |\Gamma_j|,\end{eqnarray*}
where ${\tilde c}_j(\epsilon) =(1+\epsilon{\tilde M}_j )^{n/2}
 {c}_j(\epsilon)^{n-1}$. Note that $\varrho$ is Riemannian integrable
since
 it is non-negative bounded measurable function on $\Gamma_\varrho$.
   Therefore, letting $\epsilon\to 0$, we obtain the desired result that
\begin{eqnarray}  \label {6..0} A(\tau)\sim \frac{\omega_{n-1}\tau^{n-1}}
{(4\pi)^{(n-1)}}\int_{\partial \Omega}
  \varrho^{n-1} ds \quad \; \mbox{as}\;\; \tau\to +\infty. \end{eqnarray}

  b) \ \  Next, since a $C^2$-smooth metric $g$
     can be approximated in $C^2$ by a metric $g'_\epsilon$
  which is $C^2$-smooth on $\mathcal{M}$ and piecewise real analytic
  (i.e., $g'_\epsilon$ is $C^2$-smooth and $g'_\epsilon$ is composed of a finite number
  of real analytic functions) in any compact submanifold of $(\mathcal{M}, g)$ such that
 \begin{eqnarray*} |g'^{il}_\epsilon - g^{il} |<
 \epsilon, \quad i,l=1,\cdots, n, \end{eqnarray*}
   \begin{eqnarray*} \bigg| \frac{1}{\sqrt{|g'_\epsilon|}} \frac{\partial}{\partial x_i}
 \big(\sqrt{|g'_\epsilon|} \,g'^{il}_\epsilon \big)-\frac{1}{\sqrt{|g|}} \frac{\partial}{\partial x_i}
 \big(\sqrt{|g|}\, g^{il}\big)\bigg|< \epsilon, \quad \, i,l=1,\cdots, n, \end{eqnarray*}
  for all points in  $\bar D$, with any given $\epsilon>0$.
 In addition, any bounded domain $D$ with $C^{2}$-smooth  boundary
 can also be approximated (see, the definition
  in Section 2) by
  domain $D'_\epsilon$ with $C^2$-smooth and piecewise real analytic
 boundary. Thus, the methods of Lemma 6.1 and a) still work in this case, so that
  we can estimate the eigenvalues
  for $g'^{il}_\epsilon$ in $D'_\epsilon$. But for these eigenvalues (\ref{6..0}) is true.
  Therefore, letting $\epsilon \to 0$ and noticing that $ds'_\epsilon\to ds$, we get
  that (\ref{6..0})  also holds for the $C^2$-smooth metric $g^{il}$ and $D$.
   \  \  $\quad \;\; \square$

 c) \ \  With the same arguments as in the case b), we immediately see that the formula (\ref{1-6}) is
 still true for a bounded domain with a piecewise $C^2$-smooth
 boundary in a $C^1$-smooth Riemannian manifold.

\vskip 0.28 true cm

\noindent {\bf Remark 6.4.}  \ \   Our method in the proof of
Theorem 1.1 is new and significantly different from that of
\cite{Sa}. In \cite{Sa}
   Sandgren used a technique of Lipschitz image of a convex subset for the harmonic Steklov problem.
    In our proof, $\Gamma_j$ needn't be the imagine of a convex subset.
    Next, in order to estimate $A^d(\tau)$, we introduce a new counting function $A^f(\tau)$
   as done in Section 5.
  In addition, we use the uniform boundedness of the constants $M_j$ and ${\tilde M}_j$
   to estimate the asymptotic behavior
   for any finer division according to Lemma 6.1.

 \vskip 0.398 true cm

\noindent  {\bf Proof of Corollary 1.2.} \ \  By (\ref{1-6}), we
have
\begin{eqnarray} \label{07-1}  A(\lambda_k)\sim
\frac{\omega_{n-1}\lambda_k^{n-1}}{(4\pi)^{(n-1)}}
   \big(\mbox{vol}(\partial \Omega)\big), \quad \, \; \mbox{as}\;\; k\to +\infty.\end{eqnarray}
 Since  $A(\lambda_k)=k$, we obtain (\ref{1-7}), which completes the proof.
 \  \ $\; \square$

\vskip 1.68  true cm

\centerline {\bf  Acknowledgments}

\vskip 0.38 true cm

\vskip 0.38 true cm
 I wish to express my sincere gratitude to Professor
 L. Nirenberg, Professor Fang-Hua Lin and Professor J. Shatah for their
  support during my visit at Courant Institute.
  This research was also supported by China
  Scholarship Council (No.2004307D01) and NSF of China.

  \vskip 1.65 true cm

\vskip 0.32 true cm

\end{document}